\documentclass[pdflatex,referee,sn-mathphys-num]{sn-jnl}

\usepackage{latexsym,amssymb,amsmath,amsbsy,amstext} 
\usepackage{mathrsfs}
\usepackage{amsfonts}
\usepackage{textcomp}
\usepackage{graphics,graphicx}
\usepackage{multirow}
\usepackage{subcaption,caption}
\usepackage{multirow}
\usepackage{array}
\usepackage{color,colordvi}
\usepackage{enumerate}
\usepackage{bm}
\usepackage{algpseudocode}
\usepackage{algorithm}
\usepackage{siunitx}
\usepackage{hyphenat}
\usepackage{cite}
\usepackage{mathtools}
\DeclareMathOperator*{\argmin}{arg\,min}
\usepackage{natbib} 
\usepackage[normalem]{ulem} % Load the ulem package for sout
\usepackage{bclogo} % tetrahedron
\usepackage{xcolor}
\usepackage[nobottomtitles*]{titlesec} 

\usepackage{tikz}
\usetikzlibrary{decorations.markings}
\usetikzlibrary{shapes.geometric}
\usetikzlibrary{arrows.meta,calc,3d,positioning}

\usepackage[utf8]{inputenc}
\usepackage[T1]{fontenc}
\usepackage{lmodern} % provides the "minted-style" Latin modern font
\usepackage{listings}

\definecolor{commentred}{rgb}{0.7,0,0}      % Deep Red
\definecolor{keywordblue}{rgb}{0,0,0.8}     % def, for, in, return
\definecolor{builtinteal}{rgb}{0,0.4,0.4}   % range, print, len, np.zeros
\definecolor{funcpurple}{rgb}{0.5,0,0.5}    % compute_g_function, boundary_function

\definecolor{stringgreen}{rgb}{0,0.5,0}     % Dark Green for strings/docstrings
\definecolor{bglight}{rgb}{0.95,0.95,0.95}  % Light gray background
\definecolor{codegreen}{rgb}{0,0.5,0}
\definecolor{codedarkgreen}{rgb}{1,0.3,0}
\definecolor{codeblue}{rgb}{0,0,1}
\definecolor{codegray}{rgb}{0.5,0.5,0.5}
\definecolor{codepurple}{rgb}{0.8,0,0.2}

\lstdefinestyle{mystyle}{
     backgroundcolor=\color{bglight},
     basicstyle=\fontfamily{lmtt}\selectfont\footnotesize, % Latin Modern Typewriter
     language=Python,
     keywords={}, 
     emph={np, numpy, compute_wachspress_square, compute_g_function, boundary_function},
     emphstyle=\color{codeblue}\bfseries,
     morekeywords={import, def, for, range, return},
     keywordstyle=\color{codegreen}\bfseries,
     moreemph=[2]{in},
     emphstyle=[2]\color{codepurple}\bfseries,
     moreemph=[3]{as},
     emphstyle=[3]\color{codegreen}\bfseries,
     commentstyle=\color{commentred}\itshape,
     stringstyle=\color{commentred}\itshape,
     numbers=left,                    % Position of the numbers
     numberstyle=\tiny\color{codegray}, % Size and color of the numbers
     stepnumber=1,                    % Number every line
     frame=single,
     breakatwhitespace=false,
     breaklines=true,
     captionpos=t,
     abovecaptionskip=10pt,
     belowcaptionskip=12pt,
     numbers=left,
     numbersep=10pt,
     columns=fullflexible,
     keepspaces=false,
     showspaces=false,
     showstringspaces=false
     showtabs=false,
     tabsize=4
}
\newcommand{\fref}[1]{Fig.~\ref{#1}}

\newcommand{\sref}[1]{Section~\ref{#1}}
\renewcommand{\Re}{{\mathbb{R}}}
\newcommand{\vm}[1]{\bm{#1}}
\newcommand{\vx}{\vm{x}}

\renewcommand{\it}[1]{{\textit{#1}}}

\newcommand{\talpha}{\tilde{\alpha}}

\newcommand{\liftTFI}{\mathscr{L}}
\newcommand{\proj}{\mathbb{P}}
\newcommand{\B}{\mathcal{B}}

\newcommand{\N}{{N_\theta}}

\newcommand{\Nxt}{\N (\vx;\vm{\theta})}
\newcommand{\NLt}{\N (\vm{\lambda};\vm{\theta})}
\renewcommand{\u}{{u_\theta}}

\newcommand{\uTFI}{u_\theta^{\mathrm{TFI}}}

\newcommand{\uxtTFI}{u_\theta^{\mathrm{TFI}}(\vx;\vm{\theta})}

\newcommand{\uLtTFI}{u_\theta^{\mathrm{TFI}}(\vm{\lambda};\vm{\theta})}
\newcommand{\uADF}{u_\theta^{\mathrm{ADF}}}

\graphicspath{ {./} {./figs/} 
               {./figs/ADF/} 
               {./figs/Transfinite/} 
               {./figs/GBC/}
               {./figs/Numerical_Examples/Harmonic/}       
               {./figs/Numerical_Examples/HarmonicOsc/} 
               {./figs/Numerical_Examples/Inverse/} 
               {./figs/Numerical_Examples/NonLinearPoisson/} 
               {./figs/Numerical_Examples/ParametricQuad/} 
               {./figs/Numerical_Examples/Quadrilateral/} 
               {./figs/Numerical_Examples/Eikonal/} 
               }

\makeatletter
\def\input@path{ {./}{./Python/} }
\makeatother

\begin{document}

\title{A Wachspress-based transfinite formulation 
for exactly enforcing Dirichlet boundary conditions on convex polygonal domains in physics-informed neural networks}

\author*[1]{\fnm{N.} \sur{Sukumar}}\email{nsukumar@ucdavis.edu}

\author[2]{\fnm{Ritwick} \sur{Roy}}\email{Ritwick.ROY@3ds.com}

\affil[1]{\orgdiv{Department of Civil and Environmental Engineering}, \orgname{University of California}, \orgaddress{\street{One Shields Avenue}}, \city{Davis}, \postcode{95616, \state{CA}, \country{USA}}}

\affil[2]{\orgdiv{3DS Simulia}, \orgname{Dassault Syst{\`e}mes Inc.}, \orgaddress{\street{1301 Atwood Avenue}, \city{Johnston}, \postcode{02919}, \state{RI}, \country{USA}}}

\abstract{ \unboldmath
In this paper, we present a Wachspress-based transfinite formulation on convex polygonal 
domains for exact enforcement of
Dirichlet boundary conditions in physics-informed  neural networks. This approach 
leverages prior advances in geometric design
such as blending functions and transfinite interpolation over convex 
domains. For prescribed Dirichlet 
boundary function $\B$, the transfinite interpolant of $\B$, 
$g : \bar P \to C^0(\bar P)$,  
\emph{lifts} functions from the boundary
of a two-dimensional polygonal domain to its interior. 
The transfinite
trial function is expressed as the difference between the neural network's output and
the extension of its boundary restriction into the interior of the domain, with $g$ added to it.
This ensures kinematic admissibility of the trial function in the deep Ritz method. 
Wachspress coordinates for an $n$-gon are 
used in the transfinite formula, which 
generalizes bilinear Coons transfinite interpolation on rectangles 
to convex polygons. Since Wachspress coordinates
are smooth, the neural network
trial function has a bounded Laplacian,
thereby overcoming a limitation in a previous contribution where approximate distance
functions were used to exactly enforce Dirichlet 
boundary conditions.
For a point $\vx \in \bar{P}$,  
Wachspress coordinates
$\vm{\lambda} : \bar P \to [0,1]^n$ serve as
a geometric feature map for the neural network: $\vm{\lambda}$ 
encodes the boundary edges of the polygonal
domain.
This offers a framework for solving problems 
on parametrized convex geometries using neural networks. The
accuracy of physics-informed 
neural networks is successfully
assessed on forward problems (linear and nonlinear), an inverse heat conduction problem, and a
parametrized geometric Poisson boundary-value problem.}


\maketitle

\section{Introduction}\label{sec:intro}
As a computational paradigm, physics-informed neural networks 
(PINNs)~\citep{Raissi:2019:PIN} provide
new pathways to solve forward, inverse and parametric design problems. 
In addition, they not only allow the ability to incorporate data
and physics into the solution procedure, but also to 
solve for nonlinear differential 
operators~\citep{Lu:2021:LNO,Li:2020:FNO,Lu:2021:DLL}.
For the forward problem, the simplest view of PINNs is as a
collocation-based meshfree computational method on the strong form~\citep{Dissanayake:1994:NNB,Lagaris:1998:ANN}.
A neural network approximation 
is formed via composition of nonlinear functions,
with unknown parameters residing in the contribution of each neuron, which renders the trial function 
(ansatz) to not be known a priori.  To solve
a partial differential equation (PDE) over a bounded domain, a set of collocation points is chosen in 
its interior and another
set of collocation points on the boundary (\emph{soft} imposition of boundary conditions).
The objective (loss) function is written as the 
sum of the mean squared error of the PDE in the interior of the domain (PDE loss)
and the mean squared error associated with
the prescribed boundary conditions (boundary loss). 
In addition, if the solution is also provided at specific 
collocation points in the interior of the domain, then the 
mean squared
error associated with this data (data loss) is also included 
in the loss function.
The resulting loss function
is minimized via model training 
to determine the optimal set of parameters.
An overview of recent advances in scientific machine learning is presented in~\citep{Toscano:2025:PIN}.

The minimization of the loss function is a highly nonlinear, nonconvex
optimization problem. It has been broadly appreciated
that the presence of PDE loss and boundary loss terms in the
objective function adversely
affects model training~\citep{Fuks:2020:LPI,Krishnapriyan:2021:CPF,Wang:2021:UMG,Wang:2022:WWP}, 
and the know-how that this
issue is related to finding
pareto-optimal solutions
in multiobjective optimization~\citep{Rohrhofer:2023:DVP}. 
For boundary-value problems, satisfying boundary conditions exactly (\emph{hard} imposition) 
through the neural network ansatz is desirable, as it 
simplifies the optimization landscape by reducing the 
objective to a single PDE loss term. 
The main contribution of this paper is a novel construction of a transfinite neural network trial function that exactly enforces Dirichlet boundary conditions on convex polygonal domains. 
This formulation adopts Wachspress coordinates~\citep{Wachspress:2016:RBA} as blending functions in a transfinite interpolant on convex domains~\citep{Randrianarivony:2011:OTI}.
The proposed approach 
ensures that the trial function has a bounded Laplacian, thereby overcoming a  key limitation from
previous work~\citep{Sukumar:2022:EIB} in 
which the Laplacian of the trial function 
was unbounded at the vertices of a polygonal 
domain. The detailed contributions in this work with key features of the  new formulation are described 
in~\sref{subsec:contribution}.

\subsection{Related work}\label{subsec:related}
Early attempts to strongly enforce Dirichlet boundary conditions
in neural networks can be found in the works of~\citet{Dissanayake:1994:NNB},
\citet{Lagaris:1998:ANN} and~\citet{McFall:2009:ANN}, where simple
univariate polynomials in one dimension and bivariate polynomials over two-dimensional Cartesian domains were adequate. 
Let $\bar \Omega = [0,1]$ be a closed
set in one dimension and
$\N(x;\vm{\theta})$ represents the neural network's output, where the unknown weights and
biases reside in the vector
$\vm{\theta}$. Assume that homogeneous Dirichlet boundary conditions are prescribed at $x = 0$ and $x = 1$. If so, we
write the neural network approximation as: $\u(x;\vm{\theta})
= x(1-x) \N(x;\vm{\theta}$), which
exactly satisfies the 
vanishing boundary conditions at the two ends;
similarly, the initial condition $u(0) = u_0$ in an 
initial-value problem 
is met on using the ansatz $ \u(t;\vm{\theta}) = u_0 + t \N(t;\vm{\theta}$). 
This approach has also been adopted over 
tensor-product Cartesian geometries in the early papers on PINNs~\citep{Lu:2021:DLL,Lu:2021:PIN,Yu:2022:GEP} and deep Ritz~\citep{E:2018:DRM}. As a method for transfinite
interpolation over Cartesian domains, the theory of
functional connections (TFC)~\citep{Mortari:2017:TOC,Mortari:2019:MTC} has been applied in PINNs to solve PDEs with exact enforcement of Dirichlet boundary conditions~\citep{Leake:2020:DTF,Schiassi:2021:ETF}.  Over the rectangle, the method arising from 
TFC~\citep{Mortari:2019:MTC} is identical to 
bilinear Coons transfinite interpolation~\citep{Coons:1967:SCA}. A
comprehensive and detailed discussion of Coons transfinite interpolation is presented 
in~\citet{Provatidis:2019:PIA}.

For complex geometries, using an independent low-capacity 
neural network to approximate the Dirichlet
boundary conditions was proposed in~\citet{Berg:2018:UDA}, which
has been adopted in many other studies with PINNs. 
A more general framework for the hard imposition of Dirichlet and Robin boundary conditions was introduced 
in~\citet{Sukumar:2022:EIB}, utilizing approximate distance functions (ADFs) based on the theory of 
R-functions~\citep{Rvachev:1995:RBV,Shapiro:2007:SAG}.
A limitation that
was noted in~\citep{Sukumar:2022:EIB} was that the Laplacian of the ADF became
unbounded at the vertices of a polygonal domain. Hence, to
ensure accurate model training, collocation
points could not be chosen very close to these vertices.  
The benefits in accuracy when using
\emph{hard} imposition of Dirichlet
boundary conditions over \emph{soft} imposition of
Dirichlet boundary condition (penalizing terms in the loss 
function or using Nitsche's
method) have been demonstrated in~\citep{Berrone:2023:EDB,Barschkis:2023:ESB,Cooley:2024:FPS,Sokratis:2024:RBA,Toscano:2025:PIN,Gladstone:2025:FOP,Deguchi:2025:REI}. It has been pointed out
that accuracy is affected by the 
form of the ansatz to strongly enforce Dirichlet boundary conditions. It has also been shown
in~\citep{Zeinhofer:2025:UFE} that
imposing Dirichlet boundary conditions via
penalty methods (soft imposition)
compromises the theoretical $L^2$ error estimates. 

\subsection{Contributions}\label{subsec:contribution}
For representing complex objects, transfinite interpolation over boundary curves and smoothly 
connecting surface patches are of fundamental
importance and interest in computer-aided geometry 
design. \citet{Gordon:1973:TEM} introduced the concept of transfinite 
interpolation, which refers to interpolation methods that match a given 
function on a nondenumerable set of points. 
The well-known bilinear and bicubic Coons patch use Boolean sum to map a square to a curved quadrilateral with given boundary curves~\citep{Coons:1967:SCA}. The analogous scheme on
triangles (triangular Coons interpolation) 
was proposed by~\citet{Barnhill:1973:SIT}. 
For $n$-sided polygons, V{\'a}rady et al.~\citep{Varady:2011:TSI} devised side blending and
corner blending functions that use distance functions 
to construct a transfinite interpolant. 
Related earlier work is due to Kato~\citep{Kato:1991:GSS},
who used special side blending functions
(identical to inverse-distance
based Shepard approximation~\citep{Shepard:1968:ATD}).
A different (topological) path is taken in~\citet{Randrianarivony:2011:OTI}, where
projection operations over the faces of a convex polytope
are used to develop the formula for a 
transfinite interpolant.  Our work
adopts the construction in~\citep{Randrianarivony:2011:OTI}, 
which generalizes bilinear Coons 
interpolation over the rectangle to convex polygons.
A review on blending surfaces for polygonal patches is presented in~\citet{Varady:2024:GMP}.

Generalized barycentric coordinates~\citep{Floater:2015:GBC,Hormann:2017:GBC} are
an extension of barycentric coordinates 
on simplices to polygons and polyhedra. 
In this work, 
Wachspress coordinates (nonnegative, smooth rational 
functions)~\citep{Wachspress:2016:RBA}
 over convex polygons
are used in the formula proposed 
in~\citep{Randrianarivony:2011:OTI} to
construct a transfinite interpolant on convex
domains. Consider a polygon $P$ and 
a piecewise continuous Dirichlet boundary function
$\B: \partial P \to \Re$. The
transfinite interpolant in~\citep{Randrianarivony:2011:OTI} is constructed
(hereafter identified by the function $g$) such 
that its restriction on the boundary is precisely 
$\B$. We
view this formula through the lens of 
a \emph{lifting operator}---$g$ extends $\B$ 
from the boundary to the interior of the polygon.
We adopt this lifting operator to
form the transfinite trial function in PINNs, which is expressed
as the difference between the neural network's 
output 
and the extension of its boundary restriction into the interior of 
the polygonal domain,  
with $g$ added to it. Since the
Dirichlet boundary conditions are strongly
enforced, this ensures that the transfinite trial function is kinematically admissible in the deep Ritz method for second-order boundary-value problems. As a consequence of 
the smoothness of
Wachspress coordinates, the Laplacian 
of the trial function 
is bounded, and therefore the
proposed approach is well-suited to solve PDEs with Dirichlet boundary conditions; in addition, it overcomes the previously stated limitation of  
ADFs~\citep{Sukumar:2022:EIB}.
The Wachpress-based transfinite formulation that is proposed is agnostic to the specific mathematical
expression of the 
Dirichlet function $\B$ that is prescribed on $\partial P$.

For a
convex polygon with vertices $\{v_i\}_{i=1}^n$ that are located at coordinates
$\{ \vx_i\}_{i=1}^n := (x_i,y_i)_{i=1}^n$, 
Wachspress coordinates 
are rational $C^\infty$ functions that are nonnegative, 
form a partition of unity, and
satisfy the linear reproducing conditions:
\begin{equation}\label{eq:reproducing}
\lambda_i (\vx) \ge 0, \quad
\sum_{i=1}^n \lambda_i(\vx) = 1, \quad
\sum_{i=1}^n \lambda_i(\vx) \vx_i = \vx .
\end{equation}
Equation~\eqref{eq:reproducing} reveals that Wachspress
coordinates are a convex combination and 
possess affine invariance. In addition,
$\lambda_i(\vx_j) = \delta_{ij} $ (Kronecker-delta property),
and on any boundary edge only two basis functions 
are nonzero (piecewise affine). These properties reveal 
that as barycentric coordinates on a polygon that 
encode the boundary edges, it is appealing to use 
$\vm{\lambda}$ as a geometric feature map (a new contribution) in a neural network.
In doing so, we show that a framework emerges for solving problems on parametrized convex geometries using neural networks.

\subsection{Outline}
The remainder of this paper is structured as follows. \sref{sec:unbounded_Laplacian} describes the issue of the unbounded Laplacian of the ADF at the vertices of a polygonal domain~\citep{Sukumar:2022:EIB}.
\sref{sec:wsp} presents the construction of 
Wachspress coordinates over polygons and their numerical computations. \sref{sec:formulation} outlines the 
transfinite formulation. First, we introduce the bilinear Coons interpolant on the square. Then, projection operations on the faces
of a convex polygonal domain are used to form the 
transfinite interpolant~\citep{Randrianarivony:2011:OTI}. 
Subsequently, we present the 
construction of a kinematically admissible transfinite
trial function using neural networks.
Lastly in~\sref{subsec:python}, we provide a Python implementation for the transfinite interpolant of the Dirichlet boundary 
function.  \sref{sec:model_training} presents details on the network architecture and model
training.
Numerical experiments in~\sref{sec:numerical_experiments} demonstrate the application of PINNs and deep Ritz to PDEs on square, quadrilateral, and pentagonal domains.
Several Poisson (linear and one nonlinear)
problems including one
with oscillatory boundary conditions, a parametric
geometric problem, an inverse heat
conduction problem
to determine the heat source, and the
Eikonal equation for the distance function to an interface
are considered to 
showcase the versatility and promise of the proposed
transfinite formulation to exactly enforce
Dirichlet boundary conditions.
Errors are assessed with
respect to either the exact solution (when available) or
to reference finite element solutions on highly
refined meshes using the Abaqus\texttrademark\ 
finite element
software package~\citep{Abaqus:2025}.
\sref{sec:conclusions} summarizes
the main findings from this study and provides a few 
perspectives for future work.

\section{Unbounded Laplacian of the ADF at a boundary vertex} \label{sec:unbounded_Laplacian}
\begin{figure}[!bht]
\centering
\begin{subfigure}{0.48\textwidth}
\includegraphics[width=\textwidth]{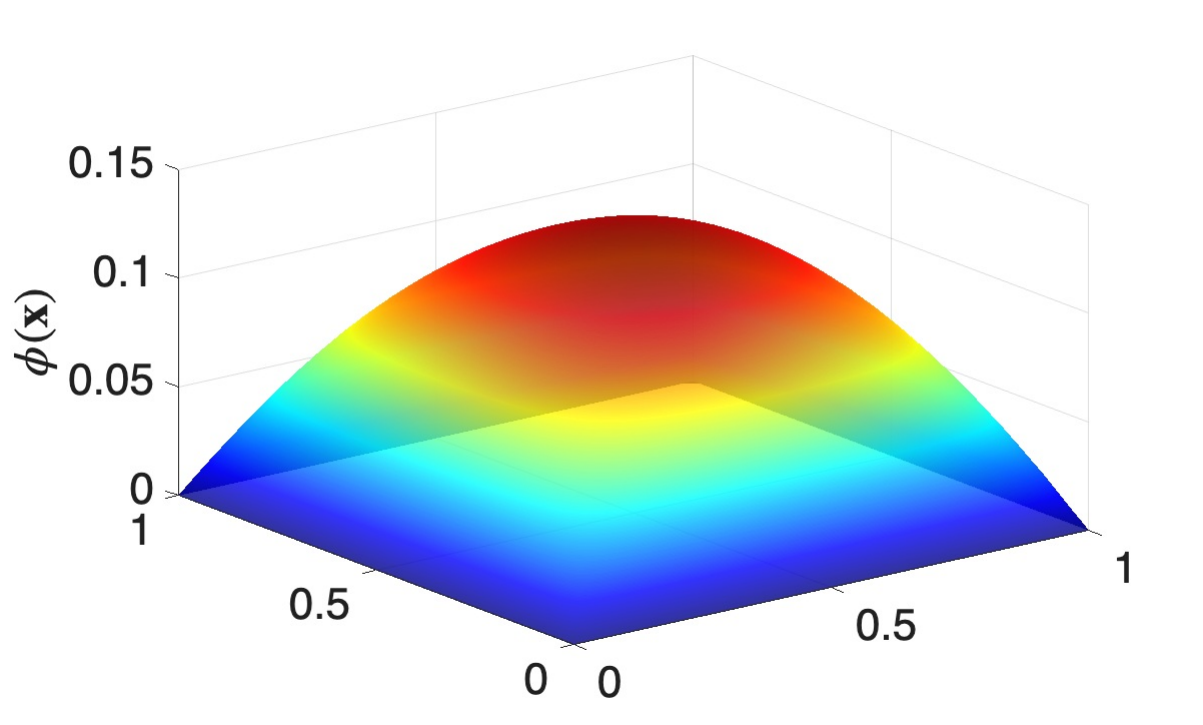}
\caption{}\label{fig:Laplacian-a}
\end{subfigure} \hfill
\begin{subfigure}{0.48\textwidth}
\includegraphics[width=\textwidth]{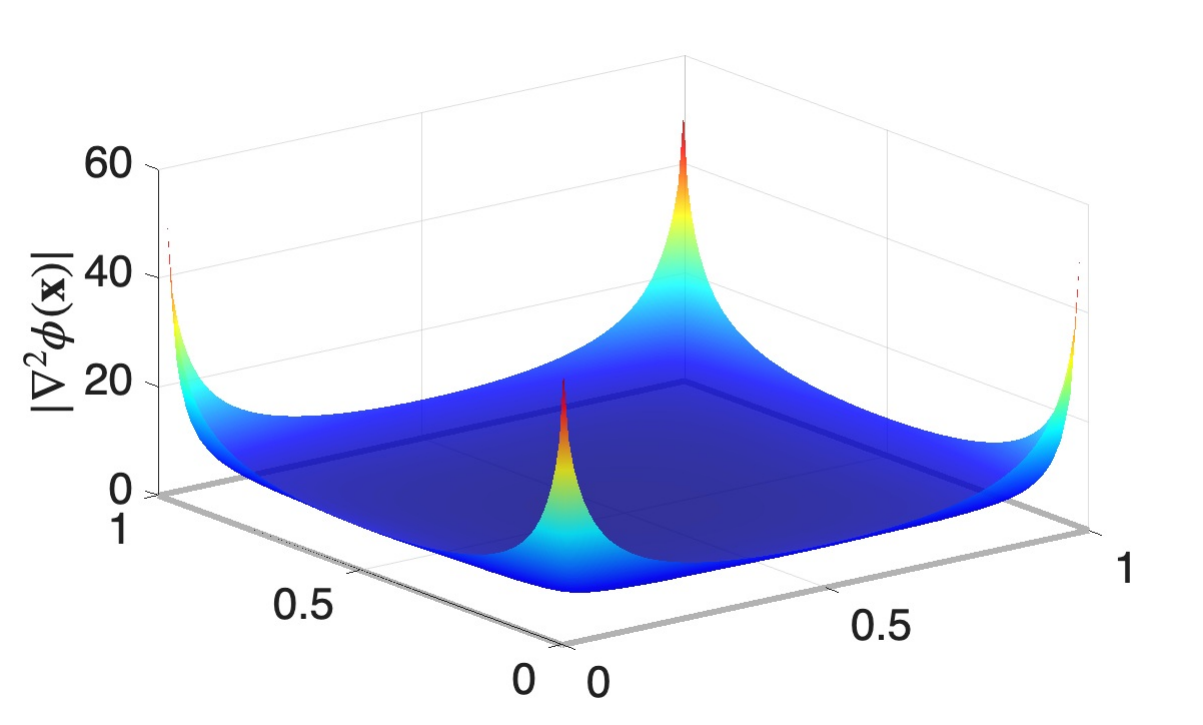}
\caption{}\label{fig:Laplacian-b}
\end{subfigure}
\medskip
\caption{Computation of $\phi$ and $\nabla^2\phi$  
         over the unit 
         square using R-equivalence with 
         order of normalization $m$ = 1
         (see~\protect\citet{Sukumar:2022:EIB}). 
         (a) 
         $\phi$; and 
         (b) $| \nabla^2 \phi |$.  Note that  
         $| \nabla^2 \phi |  \to \infty$ at the vertices.}\label{fig:Laplacian}
\end{figure}

Consider the open, bounded unit square domain, 
$\Omega = (0,1)^2$.   
The boundary $\partial \Omega$ consists of four edges that we label as
$\{e_i\}_{i=1}^4$.
For a second-order Poisson problem with homogeneous Dirichlet boundary conditions on $\partial \Omega$, the ansatz in 
PINNs is written as~\citep{Sukumar:2022:EIB}:
\begin{equation*}
u_\theta^{\mathrm{ADF}}(\vx;\vm{\theta}) 
= \phi(\vx) \N(\vx;\vm{\theta}), \quad
\phi(\vx) = \left( 
\frac{1}{\phi_1(\vx)} + 
\frac{1}{\phi_2(\vx)} + 
\frac{1}{\phi_3(\vx)} + 
\frac{1}{\phi_4(\vx)} 
\right)^{-1}, 
\end{equation*}
where R-equivalence (order of
normalization, $m = 1$) is used to form the approximate distance function to the square, $\phi(\vx)$. 
In addition, $\phi_i(\vx) = 0$ on $e_i$ and its inward
normal derivative on $e_i$ is unity. By construction, $\phi(\vx) = 0$ on 
$\partial \Omega$ and its inward normal derivative on $e_i$ is unity. 
In~\fref{fig:Laplacian}, $\phi$ and its Laplacian over the unit square are presented for $m = 1$. We note that $\phi$ is zero 
on the entire boundary and monotonic (concave) inside the domain. 
From~\fref{fig:Laplacian-b}, we observe that the
Laplacian of $\phi$ dramatically 
increases in magnitude proximal to 
the vertices of the square. In fact, it is 
known that 
$\nabla^2\phi$ is singular at the vertices 
of a polygonal domain, and therefore it is very large in magnitude near any of its vertices. 

The pathology of $\nabla^2\phi$ near corners fundamentally emerges from the simultaneous requirement that both
$\phi=0$ and ${\partial \phi} / {\partial \nu} = 1$ ($\vm{\nu}$ is the inward normal vector)
hold on each edge. Since two edges meet at a boundary vertex, these requirements are inconsistent as one approaches the vertex from each edge. The large values of the Laplacian near the corners in~\fref{fig:Laplacian-b} lead to 
a large weighting of the contributions to the loss function 
from the collocation points in the vicinity of the corners. This adversely affects model training if collocation points are sampled very close to these corners. 
To mitigate this issue, 
collocation
points were chosen in~\citep{Sukumar:2022:EIB} within a smaller 
square $[ \delta, 1 - \delta ]^2$, with
$\delta = 10^{-2}$.

\section{Wachspress coordinates on a convex polygon}\label{sec:wsp}
\begin{figure}
\centering
\begin{subfigure}{0.48\textwidth}
\includegraphics[width=\textwidth]{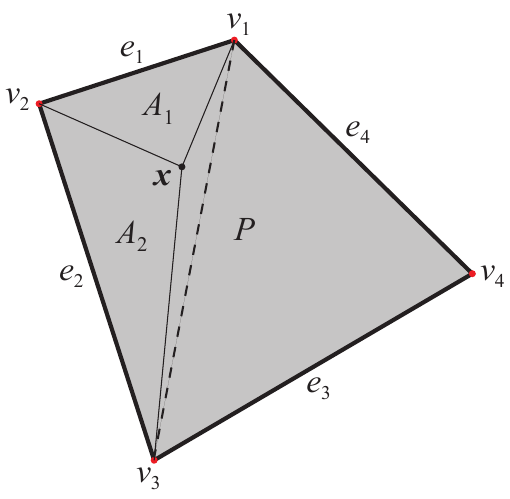}
\subcaption{}\label{fig:polygon_wsp_a}
\end{subfigure} \hfill
\begin{subfigure}{0.48\textwidth}
\includegraphics[width=\textwidth]{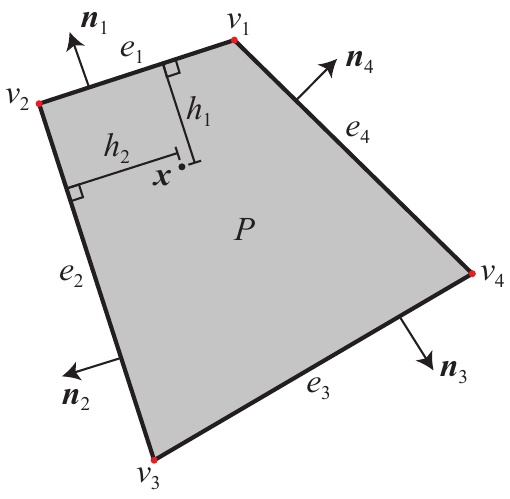}
\subcaption{}\label{fig:polygon_wsp_b}
\end{subfigure} 
\caption{Computation of Wachspress coordinates on a quadrilateral ($n = 4$) based on (a) area measures~\citep{Meyer:2002:GBC}, and
(b) edge normals and distance measures~\citep{Floater:2014:GBW}.
}\label{fig:polygon_wsp}
\end{figure}

Consider the quadrilateral ($n = 4$) shown
in~\fref{fig:polygon_wsp_a}.  
On a convex polygon $P$,
Wachspress coordinates~\citep{Wachspress:2016:RBA} are 
rational functions whose numerator and denominator are polynomials of degree
$n-2$ and $n-3$, respectively.  \citet{Meyer:2002:GBC} 
presented a simple three-point formula for Wachspress
coordinates:
\begin{equation}\label{eq:wsp_Meyer}
w_i(\vx) = \dfrac{ B_i  } 
{ A_{i-1}(\vx) A_i(\vx) },
\quad
    \lambda_i(\vx) = \dfrac{w_i(\vx)}{\sum\limits_{j=1}^n w_j(\vx)}, \ \ (i = 1,2,
    \ldots,n) ,
\end{equation}
where $\vx_i$ is the spatial coordinate of vertex $v_i$,
$B_i = A(\vx_{i-1},\vx_i,\vx_{i+1})$ is the area of the
triangle with vertices $(v_{i-1},v_i,v_{i+1})$, and
$A_i(\vx) := A(\vx_{i},\vx_{i+1},\vx)$. 
Referring to~\fref{fig:polygon_wsp_a},
$B_2 = A(\vx_1,\vx_2,\vx_3)$ is the area of 
the triangle with vertices
$(v_1,v_2,v_3$).
The
vertices of the $n$-gon are in counterclockwise orientation and cyclic ordering is assumed, i.e., 
$\vx_{n+1} := \vx_1$ and
$\vx_0 := \vx_n$.
A stable and efficient means to
compute these
coordinates on a convex
polygon is~\citep{Floater:2014:GBW}:
\begin{equation}\label{eq:wsp}
    w_i (\vx) := \frac{ \det( \vm{n}_{i-1}, \vm{n}_i) }
    { h_{i-1}(\vx) h_i(\vx) } , \quad
 \lambda_i(\vx) = \frac{w_i(\vx)}{\sum\limits_{j=1}^n w_j(\vx)} 
\ \ (i = 1,2,\dots,n), 
\end{equation}
where $\det(\cdot,\cdot)$ is the
two-dimensional scalar cross product, 
$\vm{n}_i$ is the unit outward normal to edge $e_i$, and
$h_i(\vx)$ is the distance from the point $\vx$ to the edge $e_i$
(see~\fref{fig:polygon_wsp_b}).

Both~\eqref{eq:wsp_Meyer} 
and~\eqref{eq:wsp} are not valid if the point $\vx$ lies 
on the boundary. For both formulas, an
expression for $\vm{\lambda}$ that is also valid if 
the point $\vx \in \partial P$ is available  (global form)~\citep{Floater:2015:GBC}:
\begin{equation}\label{eq:wsp_global}
w_i(\vx) = B_i
\prod_{ \underset{j \neq i-1,i}{j=1}  }^n A_j(\vx), \quad
w_i(\vx) = \det( \vm{n}_{i-1}, \vm{n}_i)
\prod_{ \underset{j \neq i-1,i}{j=1}  }^n h_j(\vx),
\end{equation}
is the weight function associated with the $i$-th vertex, which
when normalized yields $\lambda_i(\vx)$.
For the numerical experiments
in~\sref{sec:numerical_experiments},
we use the global form in~\eqref{eq:wsp_global} to compute Wachspress coordinates for $n = 5$.
On a quadrilateral,
an exact solution for $\vm{\lambda}(\vx)$ is used, which 
is presented next.

On a convex quadrilateral, moment coordinates
are introduced in~\citet{Dieci:2023:MC}, which provide 
an analytical formula for Wachspress coordinates~\citep{Dieci:2023:MC}. They are also given in~\citep{Dieci:2023:MC} for mean value coordinates~\citep{Floater:2003:MVC}
over convex or nonconvex quadrilaterals. Due to the 
availability of automatic differentiation, such formulas become attractive for applications in scientific machine learning. Consider a convex quadrilateral element. Define 
$\rho_i(\vx) := A_{i-1}(\vx) A_i(\vx) = 
\frac{1}{4} \ell_{i-1} \ell_{i} h_{i-1}(\vx) h_{i}(\vx) $, 
where
$\ell_i := |e_i| $ is the length of the $i$-th edge and $h_i(\vx)$
is the distance from $\vx$ to the $i$-th edge (see~\fref{fig:polygon_wsp_b}). Wachspress 
coordinates, 
$\vm{\lambda}(\vx) := \{\lambda_1(\vx),\, \lambda_2(\vx),\,
\lambda_3(\vx),\, \lambda_4(\vx)\}^\top$ 
are the solution to the linear system~\citep{Dieci:2023:MC}:  
\vspace*{0.5pt} 
\begin{equation}\label{eq:wsp_quad}
\begin{bmatrix}
  1 & 1 & 1 & 1 \\
  x_1 & x_2 & x_3 & x_4 \\
  y_1 & y_2 & y_3 & y_4 \\
  \rho_1(\vx) & -\rho_2(\vx) & \rho_3(\vx) & -\rho_4(\vx)
\end{bmatrix}
\begin{Bmatrix}
  \lambda_1 \\ \lambda_2 \\ \lambda_3 \\ \lambda_4
\end{Bmatrix}
= 
\begin{Bmatrix}
1 \\ x \\ y \\ 0
\end{Bmatrix} .
\end{equation}

\begin{figure}[!tbh]
\centering
\mbox{
\includegraphics[width=0.24\textwidth]{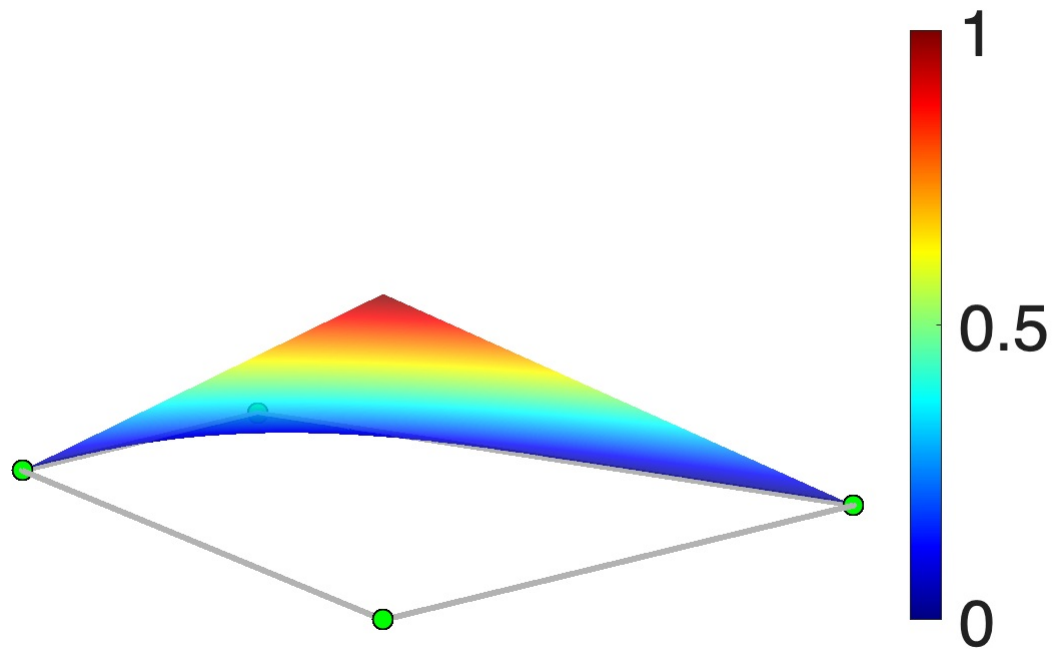}
\includegraphics[width=0.24\textwidth]{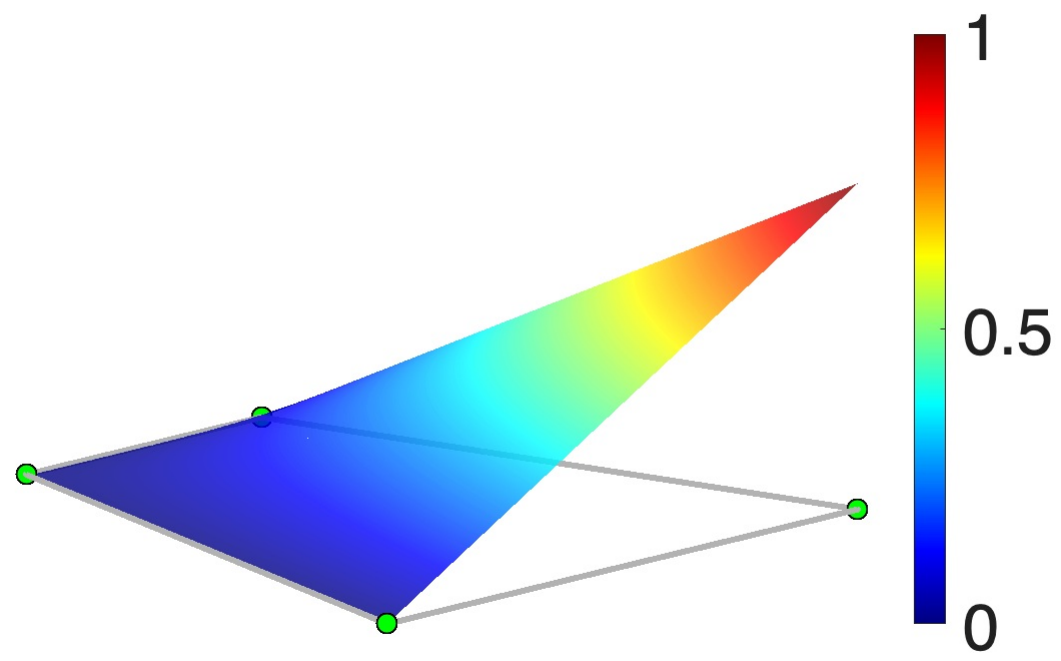}
\includegraphics[width=0.24\textwidth]{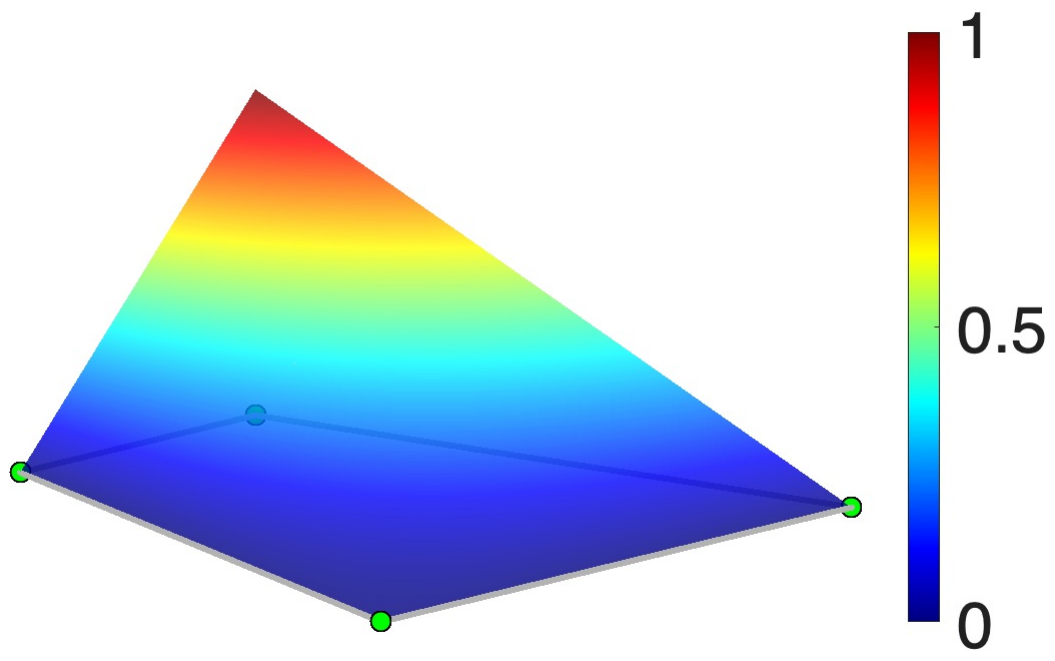}
\includegraphics[width=0.24\textwidth]{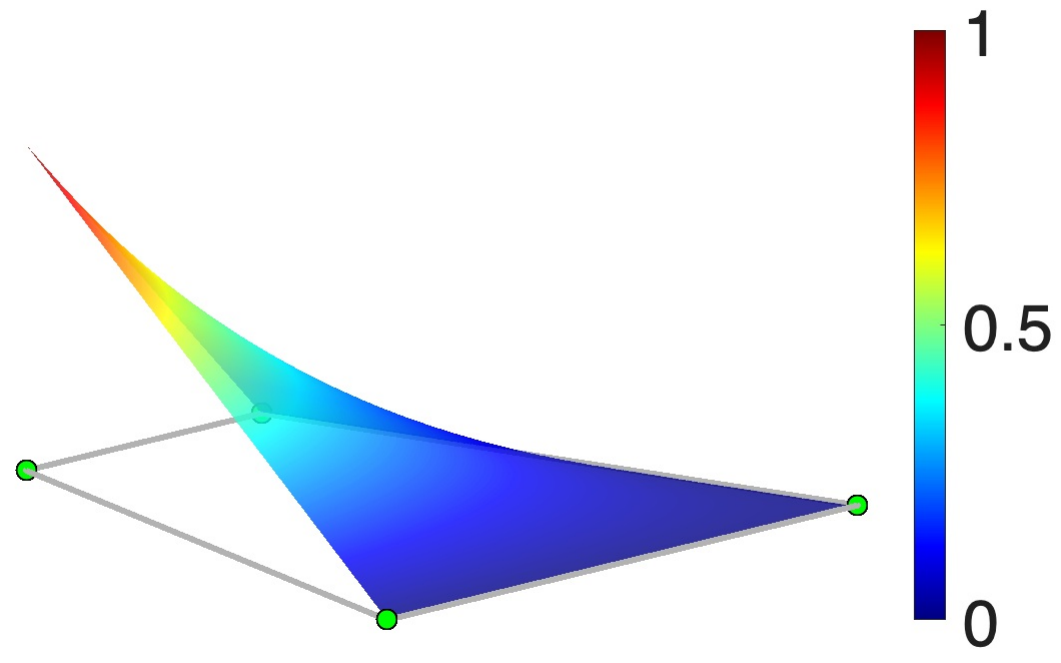}
}
\caption{Plots of Wachspress coordinates for a convex quadrilateral.} \label{fig:wsp_quad_1234}
\end{figure}

\begin{figure}[!tbh]
\centering
\mbox{
\includegraphics[width=0.19\textwidth]{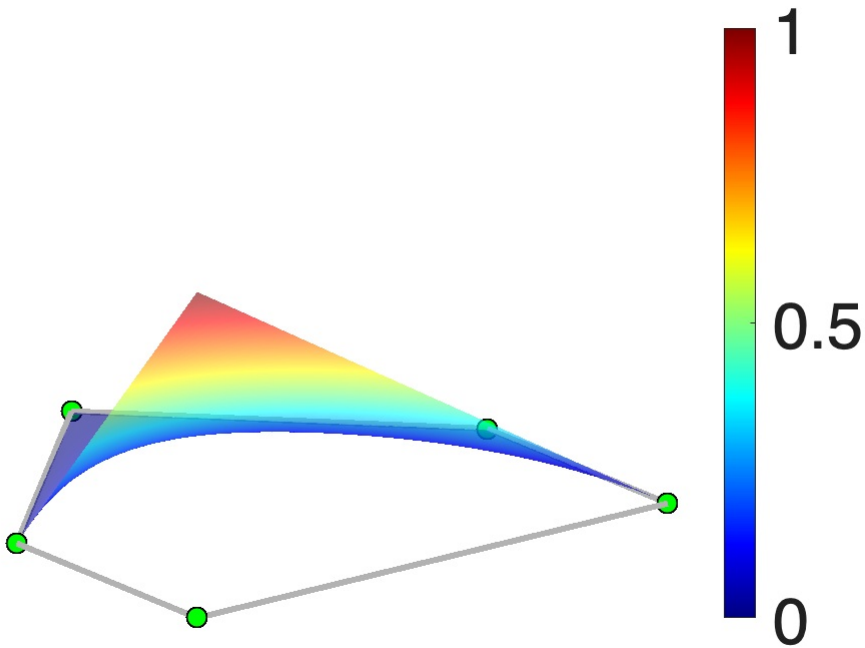}
\includegraphics[width=0.19\textwidth]{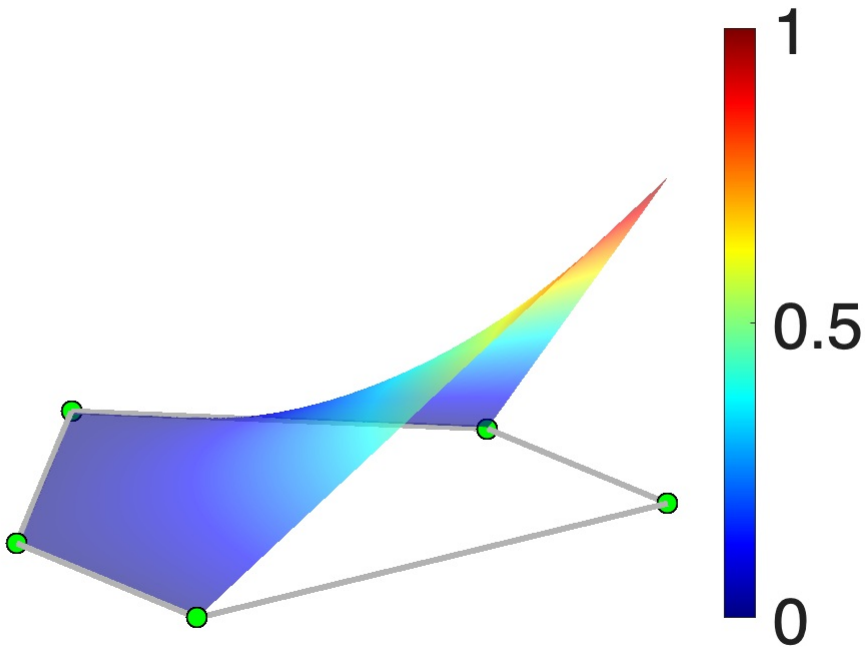}
\includegraphics[width=0.19\textwidth]{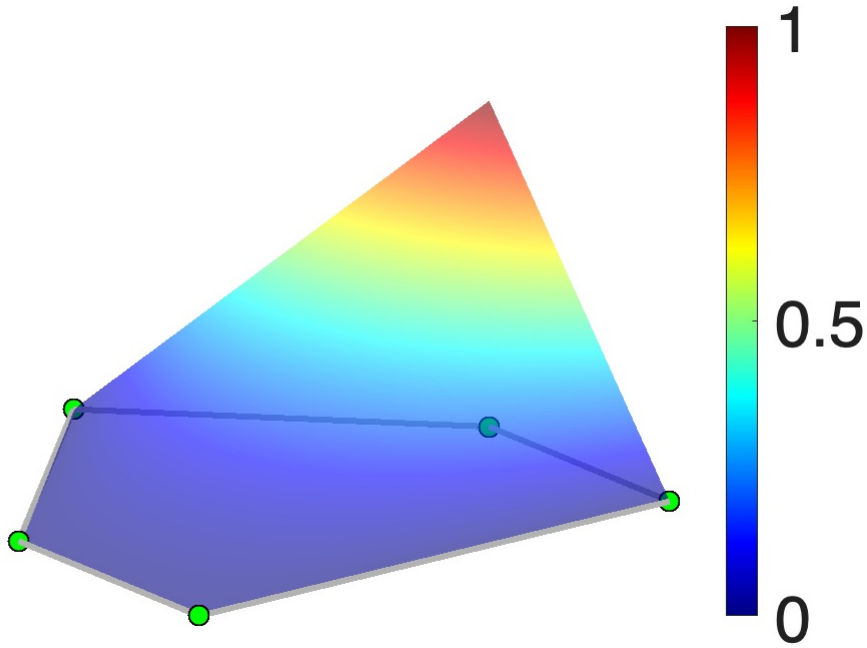}
\includegraphics[width=0.19\textwidth]{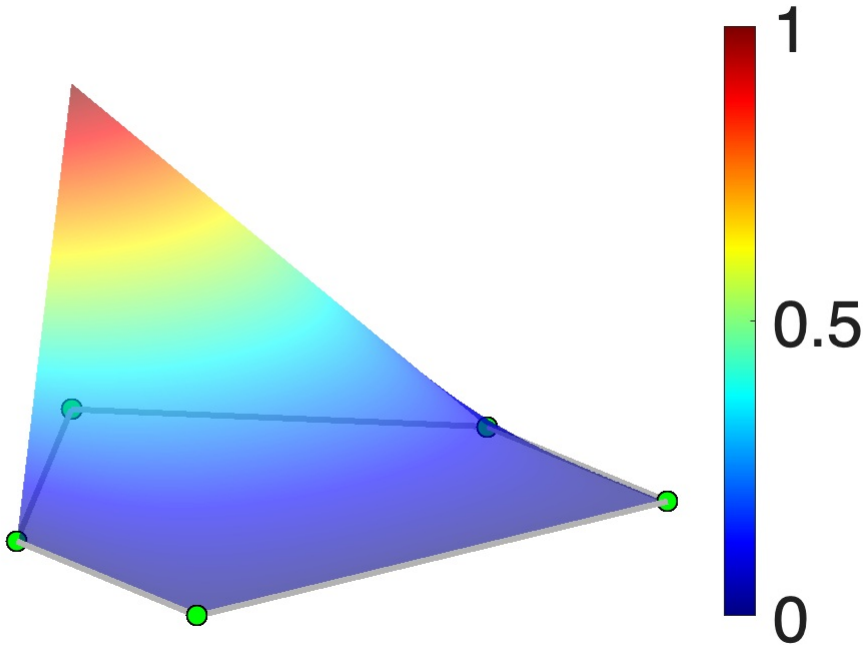}
\includegraphics[width=0.19\textwidth]{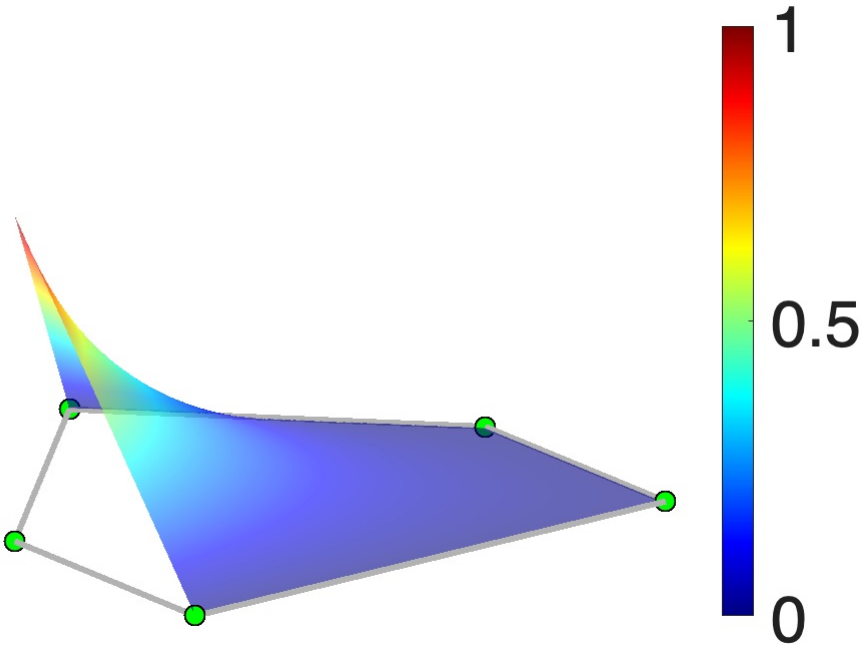}
}
\caption{Plots of Wachspress coordinates for a convex pentagon.} \label{fig:wsp_pentagon_12345}
\end{figure}

\begin{figure}[!tbh]
\centering
\mbox{
\includegraphics[width=0.24\textwidth]{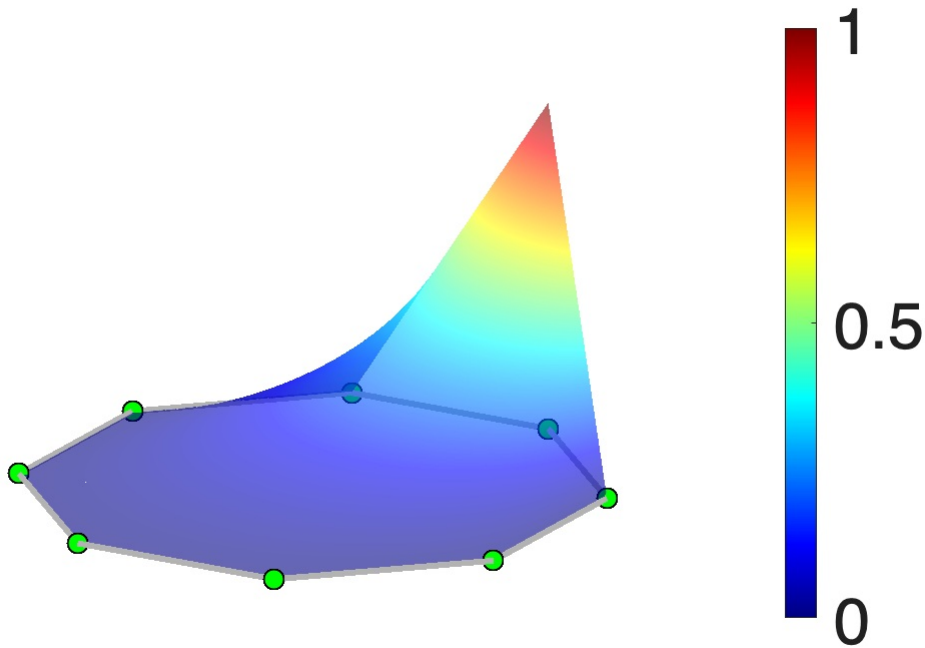}
\includegraphics[width=0.24\textwidth]{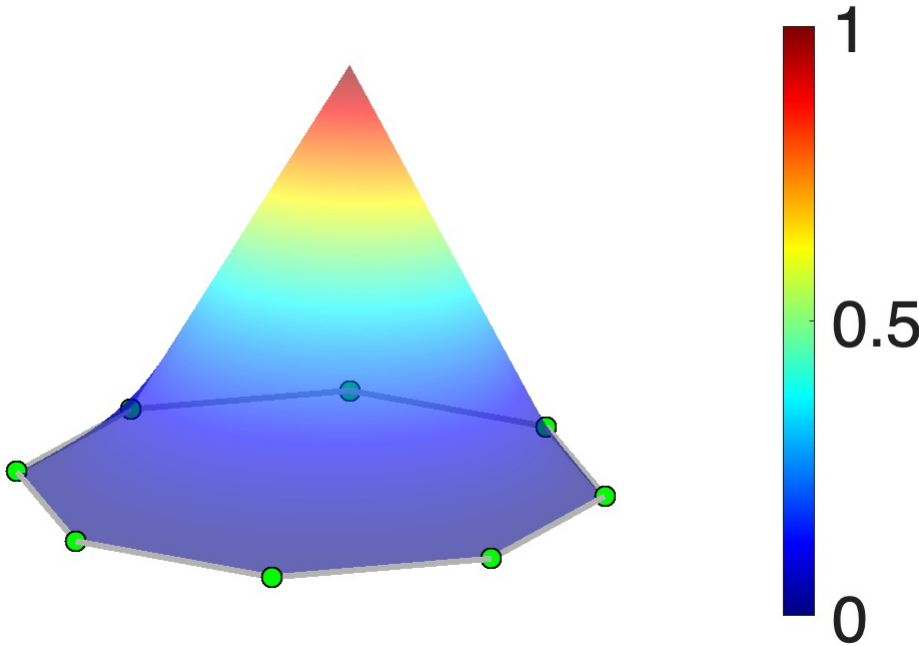}
\includegraphics[width=0.24\textwidth]{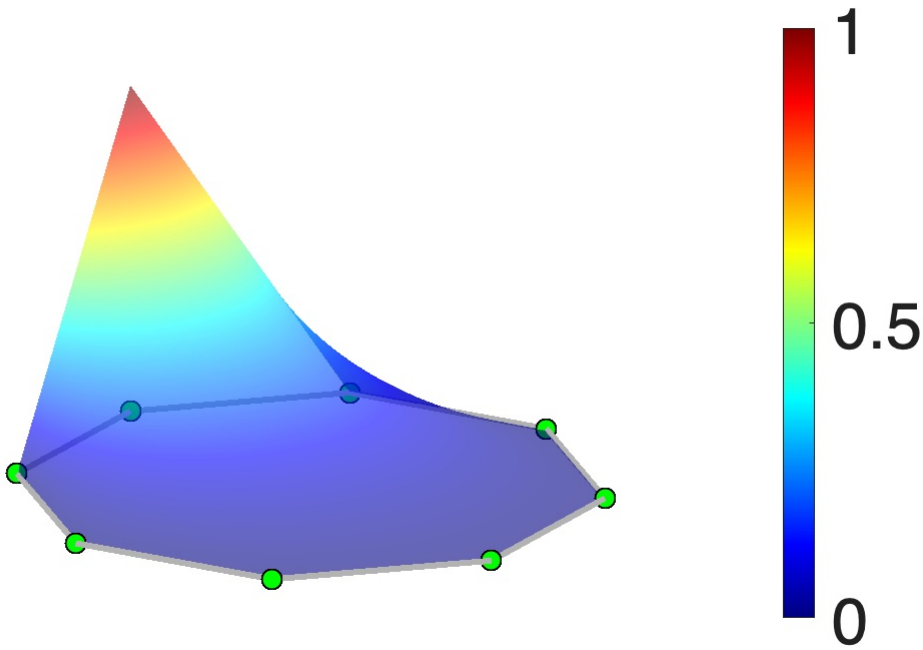}
\includegraphics[width=0.24\textwidth]{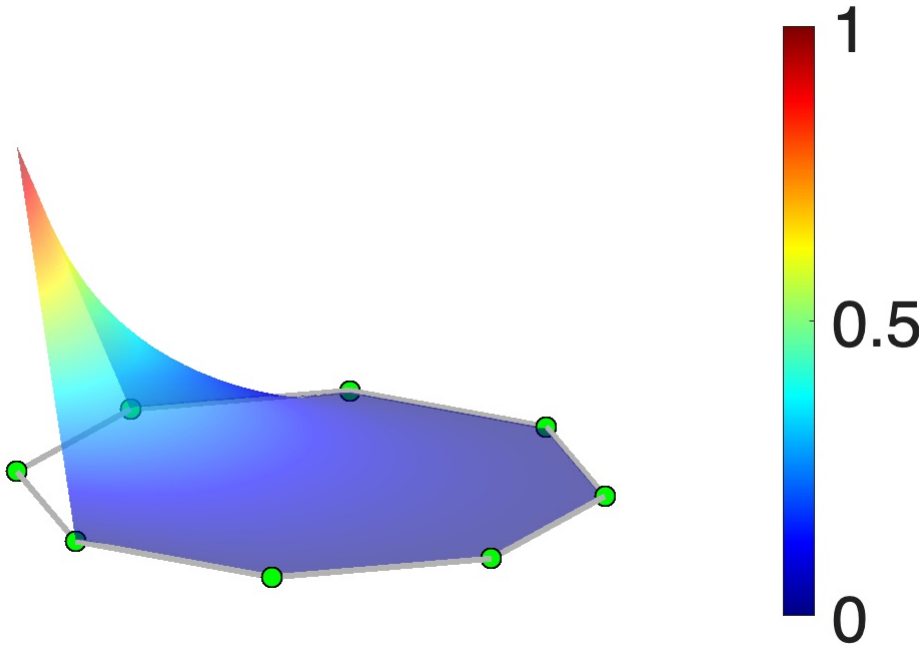}
}
\mbox{
\includegraphics[width=0.24\textwidth]{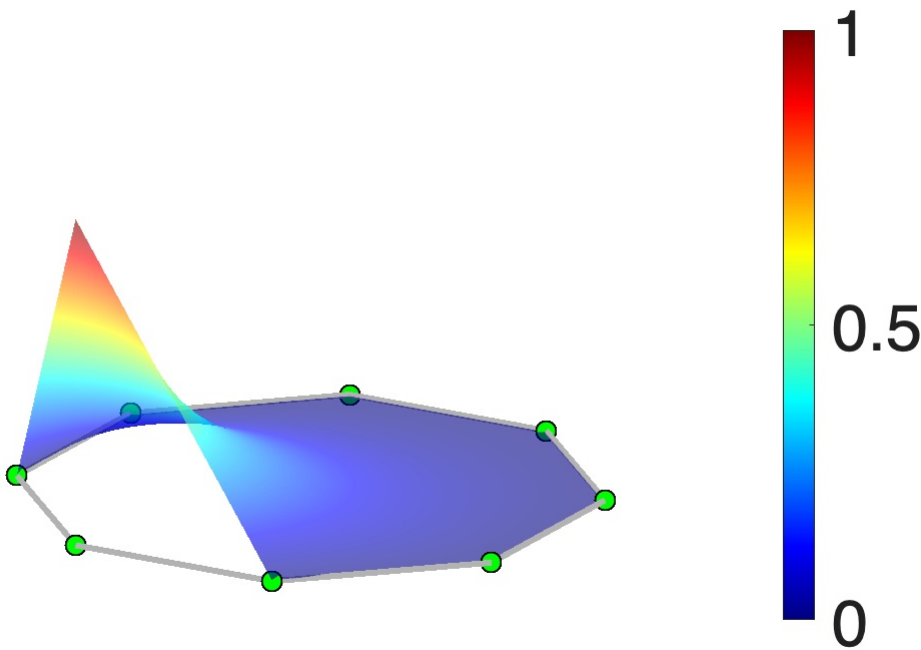}
\includegraphics[width=0.24\textwidth]{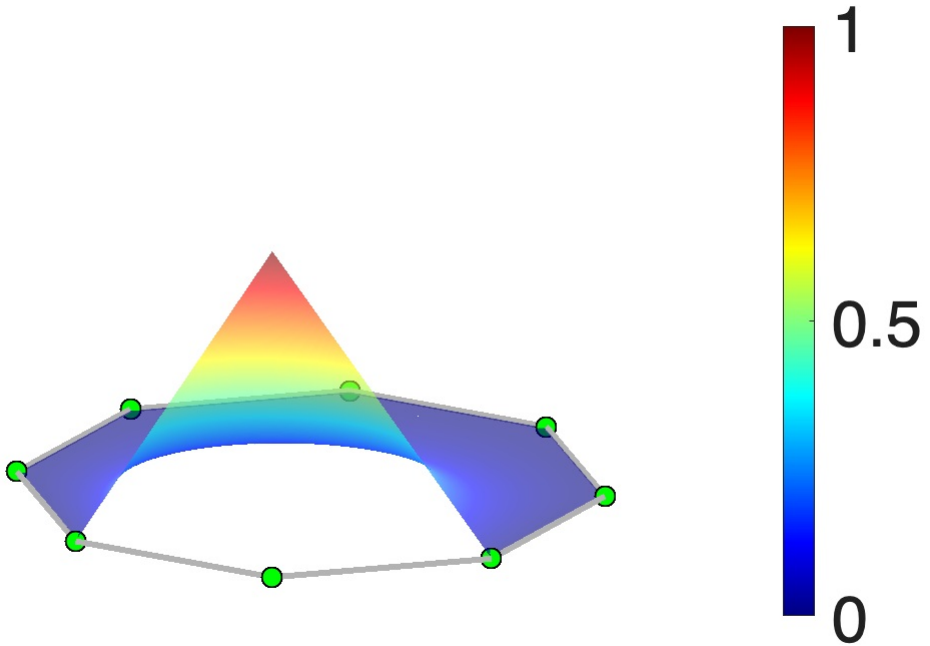}
\includegraphics[width=0.24\textwidth]{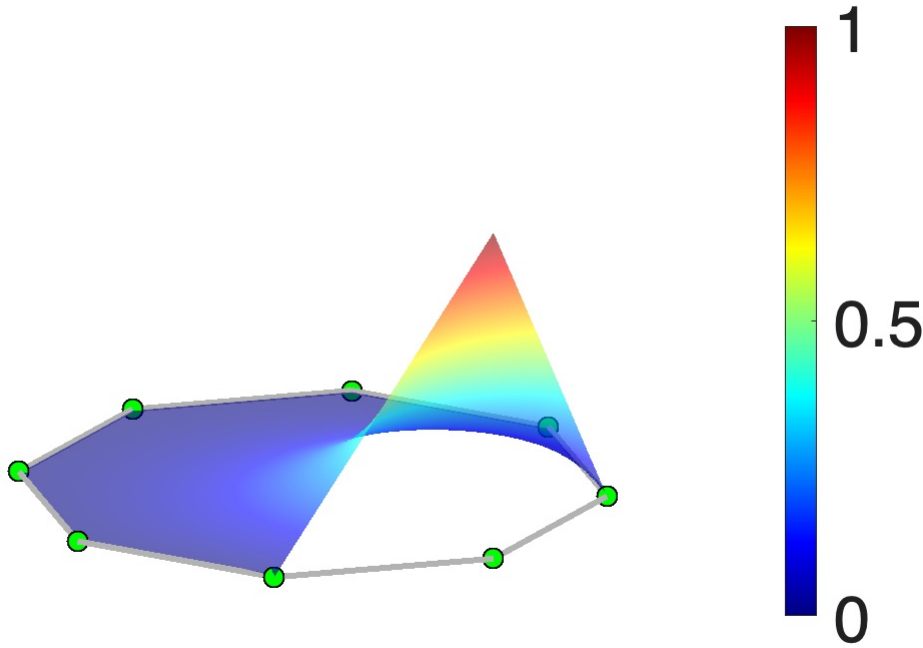}
\includegraphics[width=0.24\textwidth]{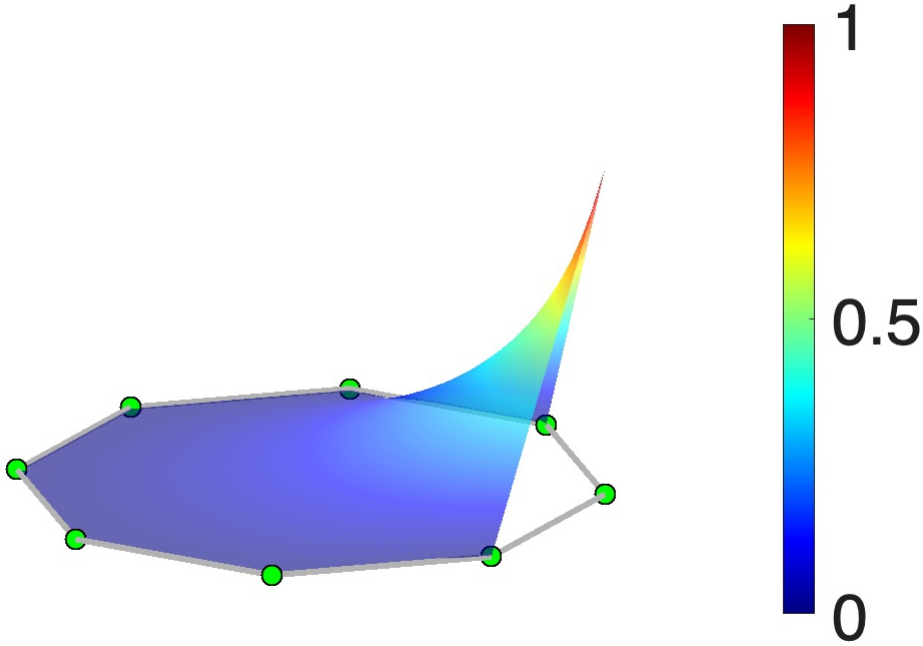}
}
\caption{Plots of Wachspress coordinates for a regular octagon.} \label{fig:wsp_octagon_12345678}
\end{figure}

In Figs.~\ref{fig:wsp_quad_1234}--\ref{fig:wsp_octagon_12345678},
plots of Wachspress coordinates on convex polygons (quadrilateral,  pentagon and octagon) are presented. Observe that on the edges of an element these coordinates share the same properties as finite element shape functions on triangles and quadrilaterals.  The Wachspress coordinate $\lambda_i(\vx)$ for vertex $v_i$ is unity at $\vx_i$
and zero on edges that do not contain $v_i$. Consequently,  only 
$\lambda_{i}$ and $\lambda_{i+1}$ (both are affine functions) are nonzero on edge $e_i$; in addition,
$\lambda_{i} + \lambda_{i+1} = 1$ on $e_i$.

\section{Formulation}\label{sec:formulation}
To fix ideas, we begin by outlining the 
construction of the neural network trial function in one 
dimension. We then present the construction of
the bilinear Coons transfinite interpolant~\citep{Coons:1967:SCA} over the unit square.
This is followed by the derivation of 
$g$, the Wachspress-based transfinite interpolant---viewed as a lifting 
operator---that extends the Dirichlet boundary function to the interior of a convex polygonal domain~\citep{Randrianarivony:2011:OTI}.  
Numerical computations are
performed to verify
that transfinite interpolation on polygonal
domains is met. Triangular, square, quadrilateral, pentagonal, and octagonal domains are considered. 
Wachspress
coordinates $\vm{\lambda}$ are used
in the input (geometric) layer of the neural network.  Let $\N(\vm{\lambda};\vm{\theta})$ be the neural network's output. 
The transfinite trial function in PINNs is constructed by first
subtracting from $\N$ the extension of the boundary restriction
of $\N$ into the interior of the polygonal domain, and then 
adding $g$ to it.  
For clarity, superscripts ADF and TFI are used to 
distinguish between the
neural network trial function based on ADFs~\citep{Sukumar:2022:EIB} from
that constructed with the
transfinite interpolant in this work.

\subsection{Univariate trial function}
\label{subsec:trialfunction_1D}
Let $N_\theta(x;\vm{\theta})$ be the
neural network's output in one dimension.
Consider the closed interval 
$\bar{\Omega} = [-1,1]$ shown in~\fref{fig:1D_domain},
and assume that Dirichlet boundary conditions
$u(-1) = \bar{u}_1$ and $u(1) = \bar{u}_2$ are 
prescribed from a function $u(x)$.  On using barycentric interpolation (convex
combination), define the following affine
function $g(x)$ that meets these
boundary conditions:
\begin{equation}\label{eq:lift_1D}
g(x) := \liftTFI[u] = 
\left( \frac{1-x}{2} \right) \bar{u}_1 +  \left( \frac{1+x}{2} \right) \bar{u}_2 ,
\end{equation}
where $\liftTFI$ is a lifting
operator, and~\eqref{eq:lift_1D} extends boundary data to $\bar \Omega$. In one dimension, $g$
is also the
harmonic extension---solution
of the Laplace equation with boundary data as
Dirichlet boundary conditions.
\begin{figure}[!tbh]
\centering
\begin{tikzpicture}[scale=1]
    \coordinate (A) at (0,0); \filldraw[blue] (A) circle (0.1cm);
    
    \node[below] at (0,-0.2) {$x = -1$};
    \node[above] at (0, 0.2) {$u(-1) = \bar{u}_1$};
    
    \coordinate (B) at (8,0); \filldraw[blue] (B) circle (0.1cm);
    \node[below] at (8,-0.2) {$x = 1$};
    \node[above] at (8, 0.2) {$u(1) = \bar{u}_2$};
   
    % Draw the line
   \draw[line width=1pt, blue] (A) -- (B);

   \draw[red, -{Stealth[length=3mm, width=3mm]}] (8.25,0.0) -- (9.25,0.0);
   \node at (9.5,0) {$x$};
\end{tikzpicture}
\caption{One-dimensional interval.}\label{fig:1D_domain}
\end{figure}

Since $g(x)$ in~\eqref{eq:lift_1D}
meets the Dirichlet boundary conditions, to form the ansatz we must add a contribution that vanishes on the boundary. To this end, we subtract from $\N(x;\vm{\theta})$ its lifting from the boundary of the interval $[-1,1]$ (affine 
interpolant), and then add $g$ to it. 
On doing so, we can write the trial function for PINNs
and deep Ritz as:
\begin{equation}\label{eq:TFI_1D}
\begin{split}
u_{\theta}^{\textrm{TFI}}(x;\vm{\theta}) 
&= g(x) + \{ \N(x;\vm{\theta}) -
   \liftTFI [ \N(x;\vm{\theta}) ] \} \\
&= g(x) + \N(x;\vm{\theta})
   - \left( \frac{1-x}{2} \right) \N(-1;\vm{\theta}) 
   - \left( \frac{1+x}{2} \right) \N (1;\vm{\theta}) ,
\end{split}
\end{equation}
which is identical to the expression in TFC that 
is derived in~\citet{Mortari:2017:TOC}. 
If interpolation of both value and derivative is desired
(for example, in clamped Euler--Bernoulli beam problems), then $C^1$ cubic Hermite finite 
element shape functions are suitable to construct 
$g(x)$.\footnote{To satisfy value-periodic boundary conditions in
one dimension (suffices for deep Ritz), one 
can set $g = \N(-1;\vm{\theta})$ in~\eqref{eq:TFI_1D}. 
Instead of linear
interpolation in~\eqref{eq:TFI_1D}, can 
adapt univariate $C^1$ Hermite shape functions 
to enforce both value- and
derivative-periodic boundary conditions on
$[-1, 1]$.  Refer
to~\citep{Dong:2021:MRP} for alternative approaches
to impose periodic boundary conditions in
neural networks.} 
It is readily verified from~\eqref{eq:TFI_1D} that the Dirichlet boundary conditions are satisfied
at $x = \pm 1$.
In contrast, the trial function in PINNs using approximate
distance functions with R-equivalence (order of normalization, $m = 1$) is~\citep{Sukumar:2022:EIB}:
\begin{equation}\label{eq:ADF_1D}
\uADF (x;\vm{\theta}) = g(x) + \phi(x)
                \N(x;\vm{\theta}), \quad
                \phi(x) = \left(\frac{1}{1+x} + \frac{1}{1-x} \right)^{-1} = \frac{1-x^2}{2},
\end{equation}
where $\phi(x)$ (ADF to the boundary) vanishes at
$x = \pm 1$ to 
ensure that~\eqref{eq:ADF_1D} matches
the Dirichlet boundary data.  
We now draw a few inferences on the properties of the 
two trial functions (TFI and ADF):
\begin{enumerate}
\item The trial functions in~\eqref{eq:TFI_1D} 
and~\eqref{eq:ADF_1D}---despite their differing forms---strictly 
enforce the Dirichlet boundary conditions.
\item On
inspection, when compared to~\eqref{eq:ADF_1D} that has
a multiplicative structure, 
the trial function $\uTFI (x;\vm{\theta})$
in~\eqref{eq:TFI_1D} 
incorporates
$\N(x;\vm{\theta})$ additively, which is simpler and it
enhances the properties of the trial space. 
If the exact solution is present in the 
span of functions that are contained in the 
neural network space, then it is desirable
that the form of the trial function chosen
in PINNs is able to represent the exact solution. 
The TFI trial function in~\eqref{eq:TFI_1D} satisfies
this property but it cannot be met by
the trial function with ADFs in~\eqref{eq:ADF_1D}. 
This is
confirmed on observing that if
$\N(x;\cdot) := u(x)$ is set in~\eqref{eq:TFI_1D}, then
$\uTFI(x;\cdot) = u(x)$ is obtained. 
\item The trial function with ADF in~\eqref{eq:ADF_1D} couples 
the behavior on the boundary and the interior due to the presence  of the product term $\phi \N$, which is undesirable. By construction,
$\phi$ satisfies
$\phi(-1)= 0$ and $\phi'(-1) = 1$. For simplicity, assume $g = 0$. Then, the spatial derivative of $\uADF$ in the vicinity of the left boundary ($\phi \to 0$, $\phi^\prime \to 1$) is
$\phi \mathbb{N}_\theta^\prime + \phi^\prime \N  \approx \N$ 
(note that $\mathbb{N}_\theta^\prime$ does not appear)
and $\nabla_{\vm{\theta}} ( \uADF ) 
= \phi \nabla_{\vm{\theta}} ( \N ) $---which can enhance stiffness and impede gradient flow (training) in residual minimization (especially, if there is fine-scale physics near the Dirichlet boundary such as a boundary layer). 
It can also lead to loss of coercivity 
in the vicinity of the boundary when deep Ritz is used: 
gradient term is dominant in the potential energy functional (loss function) 
and since 
the stiffness $\phi^2$ (multiplies $(\mathbb{N}_\theta^\prime)^2$ in the strain energy
term) goes to zero, high-frequency noise can appear in $\N$ near the boundary that will not be considered in computing the loss function.
On the other hand, the spatial derivative of $\uTFI$
is $\mathbb{N}_\theta^\prime$ plus a term that is independent of $x$, and so the potential energy functional treats 
$\mathbb{N}_\theta^\prime$ with equal importance for all points in $[-1,1]$ (as in standard PINN). 
\item The con of~\eqref{eq:TFI_1D} is that complexity increases since for
every point $x$, the points $x = -1$ and $x = 1$ must also 
be provided in the forward pass during network training.
\end{enumerate}

\subsection{Bilinear Coons transfinite interpolant over the unit square}
\label{subsec:coons}
The bilinear Coons surface patch~\citep{Coons:1967:SCA} provides continuous blending of boundary curves that form a curvilinear rectangle, or
equivalently a transfinite interpolant to the boundary functions on a rectangle.  This is achieved by 
using the Boolean sum operator that results in the convex
combination of functions that are prescribed on opposite edges of
a rectangle and a correction term containing the 
vertex values of the four boundary functions.  We now present the bilinear Coons interpolant over the unit square.

Consider a unit square domain with vertices ordered in counterclockwise orientation. On the boundary, Dirichlet functions ${\talpha}_1(x)$, $\talpha_2(y)$, 
$\talpha_3(x)$ and $\talpha_4(y)$ are prescribed
(see~\fref{fig:bilinear_Coons}), which are continuous at the vertices.  The vertex values of the boundary functions are chosen as $\talpha_1(0),\,\talpha_1(1),\,\talpha_3(1)$ and $\talpha_3(0)$. The bilinear transfinite operator 
is constructed using the Boolean sum of two linear projection operators:
\begin{equation}\label{eq:G}
    G = P_x \oplus P_y = P_x + P_y - P_x P_y, 
\end{equation}
where $P_x[\vm{\talpha}]$ and $P_y[\vm{\talpha}]$ represent functions that are convex
combinations of opposite boundary functions, and
$P_x P_y [\vm{\talpha}]$ is the bilinear interpolant of the four vertex values. In~\eqref{eq:G}, the vertex values are accounted for twice due to $P_x[\vm{\talpha}]$ and $P_y [\vm{\talpha}]$, and therefore the contribution from 
$P_xP_y[\vm{\talpha}]$ must be subtracted.
The expressions for $P_x(x,y)$, $P_y(x,y)$ and $P_x P_y (x,y)$ 
are:
\begin{subequations}
\begin{align}
    P_x(x,y) &=  P_x[\tilde{\vm{\alpha}}] = (1-x) \talpha_4(y) 
   + x  \talpha_2(y), \quad
    P_y(x,y) = P_y[\tilde{\vm{\alpha}}] = (1-y) \talpha_1(x) + y
     \talpha_3(x), \\
    P_x P_y(x,y) &= P_x P_y [\tilde{\vm{\alpha}}] = 
    (1-x)(1-y) \talpha_1(0) + 
    x(1-y)\talpha_1(1) + xy \talpha_3(1) + y(1-x) \talpha_3(0),
\end{align}
\end{subequations}
and therefore $G(x,y)$ from~\eqref{eq:G} is given by
\begin{equation}\label{eq:G_Coons}
\begin{split}
    G(x,y) = G[\tilde{\vm{\alpha}}] =& \
    (1-x) \talpha_4(y) 
   + x  \alpha_2(y) + (1-y) \talpha_1(x) + y
     \alpha_3(x) \\
    & \ - (1-x)(1-y) \talpha_1(0) -
    x(1-y)\talpha_1(1) - xy \talpha_3(1) - y(1-x) \talpha_3(0).
\end{split}
\end{equation}
On noting that $1 = 1 - x + x = 1 - y + y$, we can equivalently express $G(x,y)$ as a weighted sum of the bilinear interpolating functions:
\begin{equation}\label{eq:G_FEbasis}
\begin{split}
    G(x,y) = & \ (1-x)(1-y) \left[ \talpha_1 (x) + \talpha_4 (y) 
                 - \talpha_1( 0) \right] 
    + x(1-y) \left[ \talpha_2 (y) + \talpha_1 (x) 
                 - \talpha_1 (1)  \right] \\
    & \ + xy  \left[ \talpha_3 (x) + \talpha_2( y) 
                 - \talpha_3 (1)  \right] 
    + y(1-x) \left[ \talpha_4 (y) + \talpha_3 (x) 
                 - \talpha_3( 0) \right] ,
\end{split}
\end{equation}
which provides the
key insight underlying the construction of the transfinite interpolant over convex polygons~\citep{Randrianarivony:2011:OTI} that is presented next.
\begin{figure}
\centering
\begin{tikzpicture}[scale=1]
    % Define the square's corners
    \coordinate (A) at (0,0); \filldraw[blue] (A) circle (0.1cm);
    \node[below] at (A) {\large $(0,0)$};
   
    \coordinate (B) at (5,0); \filldraw[blue] (B) circle (0.1cm);
    \node[below] at (B) {\large $(1,0)$};
   
    \coordinate (C) at (5,5); \filldraw[blue] (C) circle (0.1cm);
    \node[above] at (C) {\large $(1,1)$};
   
    \coordinate (D) at (0,5); \filldraw[blue] (D) circle (0.1cm);
    \node[above] at (D) {\large $(0,1)$};
   
    % Draw the blue square
   \draw[line width=1pt, blue] (A) -- (B) -- (C) -- (D) -- cycle;

    % Draw red arrows below each edge
    % Bottom edge
    \draw[red, -{Stealth[length=3mm, width=3mm]}] (2.0,-0.25) -- (3.0,-0.25);
    \node at (2.5,-0.7) {{\large $\talpha_1(x)$}};
    % Right edge
    \draw[red, -{Stealth[length=3mm, width=3mm]}] (5.35,2.0) -- (5.35,3.0);
    \node at (6.05,2.5) {{\large $\talpha_2(y)$}};
    % Top edge
    \draw[red, -{Stealth[length=3mm, width=3mm]}] (2.0,5.3) -- (3.0,5.3);
    \node at (2.5,5.75) {{\large $\talpha_3(x)$}};
    % Left edge
    \draw[red, -{Stealth[length=3mm, width=3mm]}] (-0.35,2.0) -- (-0.35,3.0);
    \node at (-1.05,2.5) {{\large $\talpha_4(y)$}};
    % write vertex IDs: 1,2,3,4
    \node at (-0.3,0.1) {\large $1$};
    \node at (5.3,0.1)  {\large $2$};
    \node at (5.3,5.0)  {\large $3$};
    \node at (-0.3,5.0) {\large $4$};
\end{tikzpicture}
\caption{Bilinear Coons interpolation on the unit square.}
\label{fig:bilinear_Coons}
\end{figure}

\subsection{Wachspress-based transfinite interpolant}
\label{subsec:TFI_2D}
\citet{Randrianarivony:2011:OTI} showed that bilinear Coons interpolation can be seen from a topological perspective using blending functions that are generalized barycentric coordinates, with a formula that involves a projection onto the faces of a polytope (faces, edges and vertices for three-dimensional convex polyhedra, and edges and vertices for two-dimensional convex polygons). 
The transfinite interpolant of the Dirichlet 
boundary condition, $g$, lifts functions from the polygonal boundary to its interior.

We present the main elements of the derivation of the transfinite interpolant~\citep{Randrianarivony:2011:OTI}.  We consider the transfinite
formula for a quadrilateral from which its extension to general convex polygons becomes evident. As a particular choice of generalized barycentric coordinates, we adopt Wachspress
coordinates~\citep{Wachspress:2016:RBA}. This choice is 
guided by the
fact that Wachspress coordinates are efficient to compute (closed-form expressions),
are $C^\infty$ smooth over a polygon, 
and have a bounded Laplacian; a drawback
is that they are restricted to convex polygons. Mean value 
coordinates~\citep{Floater:2003:MVC,Hormann:2006:MVC} 
are well-defined
for convex and nonconvex planar polygons; however, they are 
$C^0$ at the vertices of a polygon and even on the square 
their Laplacian is unbounded at its
vertices (see~\fref{fig:gfunction_mvc_b} in~\sref{subsec:TFI_verification_polygons}). 
Harmonic coordinates~\citep{Joshi:2007:HCF}
satisfy Laplace equation and have
many positive attributes on arbitrary (convex and nonconvex) planar polygons, but they
are not known in closed form nor simple formulas are available to compute them.

Let $P$ be an open, polygonal domain 
with $n$ vertices. Its boundary is $\partial P$ and the 
closure of $P$ is denoted by $\bar P$. Consider a piecewise continuous 
boundary function $\B := \B(x,y)$, which
is prescribed on the boundary of the polygon.
The restriction of $\B$ to each edge is a
bivariate
function of $x$ and $y$ over a convex polygon, with $x$ and $y$ being
affinely related on each edge. Since 
$x$ and $y$ can be expressed in terms of $\vm{\lambda}$ through~\eqref{eq:reproducing}, these boundary functions can be equivalently expressed 
in terms of $\vm{\lambda}$.  Furthermore, as noted earlier
in~\sref{sec:wsp},
since on $e_i$ ($i$-th edge) only $\lambda_i$ and $\lambda_{i+1}$ are nonzero with
$\lambda_i + \lambda_{i+1} = 1$ on $e_i$, each boundary function on an edge 
can be parametrized as a function of a single variable in $[0,1]$.  Let the $n$ parametric functions
$\{\alpha_i\}_{i=1}^n : [0,1]\to \mathbb{R}$
(counterclockwise orientation starting at the bottom edge; see~\fref{fig:TFI_quad} for the case $n = 4$) 
be prescribed on $\partial P$.  With some abuse of
notation, we retain the use of $\B = \B\bigl(x(\vm{\lambda}),y(\vm{\lambda})\bigr)$ to represent the boundary function that is now a function of
Wachspress coordinates. 

On a quadrilateral, the 
Dirichlet function 
$\B(\lambda_1,\lambda_2,\lambda_3,\lambda_4) : \partial P \to \Re$ is defined 
on the boundary such that~\citep{Randrianarivony:2011:OTI}
\begin{subequations}\label{eq:B}
\begin{align}
\B(\lambda_1,\lambda_2,0,0) &= \alpha_1(\lambda_2), \ \ \lambda_1 + \lambda_2 = 1 \textrm{ on } e_1 ,\\
\B(0,\lambda_2,\lambda_3,0) &= \alpha_2(\lambda_3), \ \ \lambda_2 + \lambda_3 = 1 \textrm{ on } e_2 , \\
\B(0,0,\lambda_3,\lambda_4) &= \alpha_3(\lambda_4), \ \ \lambda_3 + \lambda_4 = 1 \textrm{ on } e_3 , \\
\B(\lambda_1,0,0,\lambda_4) &= \alpha_4(\lambda_1), \ \ \lambda_4 + \lambda_1 = 1 \textrm{ on } e_4  ,
\end{align}
\end{subequations}
and due to continuity at the
vertices, $\alpha_1,\,\alpha_2,\,\alpha_3,\,\alpha_4$ must satisfy the
following matching conditions:
\begin{equation} \label{eq:matchingconditions}
\alpha_1( 1) = \alpha_2 (0), \quad
\alpha_2( 1) = \alpha_3 (0), \quad 
\alpha_3 (1) = \alpha_4 (0), \quad 
\alpha_4 (1) = \alpha_1 (0).
\end{equation}
On the unit square (\fref{fig:TFI_quad_a}), 
Wachspress coordinates
are identical to finite element shape functions:
\begin{equation}\label{eq:lambda_square}
\lambda_1(\vx) = (1-x)(1-y), \quad
\lambda_2(\vx) = x(1-y), \quad
\lambda_3(\vx) = xy, \quad
\lambda_4(\vx) = y(1-x) ,
\end{equation}
whereas on the convex quadrilateral shown in~\fref{fig:TFI_quad_b}, Wachspress
coordinates are
rational functions that are a quotient of a bivariate quadratic and an affine function. Solving~\eqref{eq:wsp_quad} provides the exact solution for $\vm{\lambda}$ over a
convex quadrilateral. 
\begin{figure}[!hbt]
\centering
\begin{subfigure}{0.46\textwidth}
\begin{tikzpicture}[scale=0.83]
    % Define the square's corners
    \coordinate (A) at (0,0); \filldraw[blue] (A) circle (0.1cm);
    \node[below] at (A) {\Large $\lambda_1$};
    
    \coordinate (B) at (5,0); \filldraw[blue] (B) circle (0.1cm);
    \node[below] at (B) {\Large $\lambda_2$};
    
    \coordinate (C) at (5,5); \filldraw[blue] (C) circle (0.1cm);
    \node[above] at (C) {\Large $\lambda_3$};
    
    \coordinate (D) at (0,5); \filldraw[blue] (D) circle (0.1cm);
    \node[above] at (D) {\Large $\lambda_4$};
    
    % Draw the blue square
   \draw[line width=1pt, blue] (A) -- (B) -- (C) -- (D) -- cycle;

    % Draw red arrows below each edge
    % Bottom edge
    \draw[red, -{Stealth[length=3mm, width=3mm]}] (2.0,-0.35) -- (3.0,-0.35);
    \node at (2.5,-1) {{\Large $\alpha_1(\lambda_2)$}};
    % Right edge
    \draw[red, -{Stealth[length=3mm, width=3mm]}] (5.35,2.0) -- (5.35,3.0);
    \node at (6.5,2.5) {{\Large $\alpha_2(\lambda_3)$}};
    % Top edge
    \draw[red, -{Stealth[length=3mm, width=3mm]}] (3.0,5.35) -- (2.0,5.35);
    \node at (2.5,6.0) {{\Large $\alpha_3(\lambda_4)$}};
    % Left edge
    \draw[red, -{Stealth[length=3mm, width=3mm]}] (-0.35,3.0) -- (-0.35,2.0);
    \node at (-1.5,2.5) {{\Large $\alpha_4(\lambda_1)$}};
   
    % write P
    \node at (2.5,2.5) {{\Large $P$}};

    % write e_1 e_2 e_3 e_4
    \node at (2.5,0.3) {\large $e_1$};
    \node at (4.6,2.5) {\large $e_2$};
    \node at (2.5,4.65) {\large $e_3$};
    \node at (0.4,2.5){\large $e_4$};

    % write vertex IDs: 1,2,3,4
    \node at (-0.3,0.1) {$1$};
    \node at (5.3,0.1)  {$2$};
    \node at (5.3,5.0)  {$3$};
    \node at (-0.3,5.0) {$4$};
\end{tikzpicture}
\subcaption{}\label{fig:TFI_quad_a}
\end{subfigure} \hfill
\begin{subfigure}{0.46\textwidth}
\begin{tikzpicture}[scale=0.88]
% Define the quadrilateral's corners
    \coordinate (A) at (0,0); \filldraw[blue] (A) circle (0.1cm);
    \node[below] at (A) {\Large $\lambda_1$};
    
    \coordinate (B) at (5,0.5); \filldraw[blue] (B) circle (0.1cm);
    \node[below] at (B) {\Large $\lambda_2$};
    
    \coordinate (C) at (4,5.5); \filldraw[blue] (C) circle (0.1cm);
    \node[above] at (C) {\Large $\lambda_3$};
    
    \coordinate (D) at (1,4); \filldraw[blue] (D) circle (0.1cm);
    \node[above,yshift=0.08cm] at (D) {\Large $\lambda_4$};
    
    % Draw the blue quadrilateral 
   \draw[line width=1pt, blue] (A) -- (B) -- (C) -- (D) -- cycle;

    % Draw red arrows below each edge
    % Bottom edge
    \draw[red, -{Stealth[length=3mm, width=3mm]}] (2.0,-0.2) -- (3.0,-0.1);
    \node at (2.5,-0.75) {{\Large $\alpha_1(\lambda_2)$}};
    % 2nd edge
    \draw[red, -{Stealth[length=3mm, width=3mm]}] (5.0,2.4) -- (4.8,3.4);
    \node at (5.95,2.95) {{\Large $\alpha_2(\lambda_3)$}};
    % 3rd edge
    \draw[red, -{Stealth[length=3mm, width=3mm]}] (3.0,5.4)
    -- (2.0,4.9);
    \node at (2.35,5.9) {{\Large $\alpha_3(\lambda_4)$}};
    % 4th edge
    \draw[red, -{Stealth[length=3mm, width=3mm]}] (0.3,2.5) -- (0.05,1.5); 
    \node at (-0.9,2.1) {{\Large $\alpha_4(\lambda_1)$}};
    % write P
    \node at (2.5,2.5) {{\Large $P$}};
\end{tikzpicture}
\subcaption{}\label{fig:TFI_quad_b}
\end{subfigure}
\caption{Boundary functions for Wachspress-based transfinite interpolation
on (a) square and (b) quadrilateral domains.} \label{fig:TFI_quad}
\end{figure}

We now present the transfinite formula. Let $\Pi$ denote a face (dimension $0$ or $1$) of the polygon. Consider an arbitrary face $\Pi$ of $\bar P$ and let $I$ be the set that contains the
indices of the vertices that are incident upon $\Pi$. 
The set $\Pi$ contains edges and vertices, and
$I$ will contain either 2 vertices if $\Pi$ is an edge or 
$1$ vertex ($I = \Pi$) if $\Pi$ is a vertex.
Let $v_q$ be any vertex such that $q \in I$. 
Now, introduce a projection
that depends on both $\Pi$ and $v_q$~\citep{Randrianarivony:2011:OTI}:
\begin{equation}\label{proj_faces}
\proj := \proj_{\Pi,v_q} : \bar P \to \bar \Pi ,
\end{equation}
where $\proj$ takes the barycentric coordinates of a point 
in $\bar P$ and projects it onto $\bar \Pi$.
If $\vm{\lambda}$ are the barycentric coordinates
for a point $\vx \in \bar P$. then
the formula for transfinite interpolation over the 
polygonal domain
$\bar P$ is~\citep{Randrianarivony:2011:OTI}:
\begin{subequations}\label{eq:TFI_formula}
\begin{align}
g(\vm{\lambda}) := \liftTFI[\B] &= \sum_{i=1}^n 
\lambda_i \sum_{\Pi \in {\cal G}(i)}
(-1)^{ \mathrm{dim} (\Pi)+1 } \, 
 \B \circ \mathbb{P}_{\Pi,v_i}
( \vm{\lambda} ),  \label{eq:TFI_formula_a} \\
\intertext{where 
$\liftTFI$ (appeared earlier in~\sref{subsec:trialfunction_1D}) is a lifting operator that extends the boundary function
to the interior of the domain, 
$\circ$ is the composition operator, 
${\cal G}(i)$ is the set of topological entities
(edges and vertices of the polygon) such that the vertex $v_i$ is incident upon them and 
$\proj$ is a projection such that for
$\vm{\lambda} = (\lambda_1, \dots, \lambda_n)$, the
barycentric coordinates of the image $\vm{\mu} := \vm{\mu}(\vm{\lambda})
= (\mu_1, \dots, \mu_n)$ are defined as:}
\mu_j &:= \begin{cases}
0 & \textrm{if } v_j \not \subseteq \Pi \\
\lambda_j & \textrm{if } v_j \subseteq \Pi \ \textrm{and } j \neq i \\
1 - \sum\limits_{ j \in I\backslash \{ i \} }  \lambda_j & \textrm{otherwise}
\end{cases} \quad \label{eq:TFI_formula_b} .
\end{align} 
\end{subequations}
Equation~\eqref{eq:TFI_formula_a} for $g$ is a transfinite blended interpolant, where Wachspress coordinates 
$\{\lambda_i\}_{i=1}^n$ serve as blending functions. Each
Wachspress coordinate that is associated with a vertex $v_i$ is multiplied by three terms (two arise from the edges connected to $v_i$ and one from $v_i$ itself) that are a function of the boundary function 
$\B(\lambda_1,\lambda_2,\lambda_3,\lambda_4)$. To further clarify, we  refer 
to~\eqref{eq:TFI_formula_a} and obtain the three terms that are associated with the first vertex. For $i = 1$, the set
${\cal G} (1) := \{ \Pi_1^1 , \, \Pi_1^2 , \, \Pi_1^3 \} 
= \{ e_1, e_4, v_1 \}$ ($e_1$ and $e_4$ are edges; $v_1$ is
a vertex) and $\mathrm{dim} \, \bigl({\cal G} (1) \bigr) = [1 , 1 , 0]$ are the dimensions of the
entities in ${\cal G}(1)$.  
In addition, let $I_1^j$
be the index set for $\Pi_1^j$ and
$\vm{\mu}_1^j  := \proj_{\Pi_1^j,v_1} ( \vm{\lambda} ) $. 
Then, we obtain
$I_1^1 = \{1,\, 2\}$, $I_1^2 = \{1,\, 4\}$,
$I_1^3 = \{1\}$; and
$\vm{\mu}_1^1 = ( 1 - \lambda_2,\,\lambda_2, \, 0, \, 0 ) $,
$\vm{\mu}_1^2 = ( 1 - \lambda_4,\, 0, \, 0 , \, \lambda_4 ) $,
$\vm{\mu}_1^3 = ( 1,\, 0, \, 0, \, 0 )$. 
The image $\vm{\mu}$ is the projection of the generalized barycentric coordinates of a point in the polygon onto either an edge or a vertex. Due to the facet-reducing
property of nonnegative (convex)
Wachspress coordinates, on an edge $e_k$ 
only $\lambda_k$ and $\lambda_{k+1}$ are nonzero and they sum to unity, whereas at a vertex $v_k$ we obtain
$\lambda_k = 1$ and $\lambda_j = 0$
for $j \ne k$. 
Using these,
we can write
$\B \circ \vm{\mu}_1^1 = \B(1-\lambda_2,\lambda_2,0,0)$,
$\B \circ \vm{\mu}_1^2  = \B(1-\lambda_4,0,0,\lambda_4)$ and
$\B \circ \vm{\mu}_1^3  = \B(1,0,0,0)$ 
for the terms that appear
in~\eqref{eq:TFI_formula_a}.  
The contribution
from $i = 1$ in~\eqref{eq:TFI_formula_a} is
$\lambda_1 [  \B(1-\lambda_2,\lambda_2,0,0)
+ \B(1-\lambda_4,0,0,\lambda_4) - \B(1,0,0,0) ] $. 
The negative contribution in this weighted sum
is when $\Pi$ is a vertex; the same structure also appears in the bilinear Coons transfinite interpolant in~\eqref{eq:G_FEbasis}.
On following similar steps for 
$i = 2$, $i = 3$ and $i = 4$, we can express
$g(\vm{\lambda})$ for the square and quadrilateral
domains shown in~\fref{fig:TFI_quad} as~\citep{Randrianarivony:2011:OTI}:
\begin{equation}\label{eq:g_quad}
\begin{split}
g(\vm{\lambda}) =
\ & \lambda_1 \bigl[ \B(1-\lambda_2,\lambda_2,0,0)
+ \B(1-\lambda_4,0,0,\lambda_4) - \B(1,0,0,0) \bigr] \\
& + \lambda_2 \bigl[ \B(0,1-\lambda_3,\lambda_3,0)  
+  \B(\lambda_1,1-\lambda_1,0,0) - \B(0,1,0,0) \bigr] \\
& + \lambda_3 \bigl[  \B(0,0,1-\lambda_4,\lambda_4) 
+ \B(0,\lambda_2,1-\lambda_2,0) - \B(0,0,1,0) \bigr] \\
& + \lambda_4 \bigl[ \B(\lambda_1,0,0,1-\lambda_1) 
+  \B(0,0,\lambda_3,1-\lambda_3) - \B(0,0,0,1) \bigr].
\end{split}
\end{equation}
The expression for $g$
in~\eqref{eq:g_quad} holds on the convex quadrilateral in~\fref{fig:TFI_quad_b}, where
the Wachspress coordinates, $\vm{\lambda}$, are rational
functions in Cartesian coordinates~\citep{Wachspress:2016:RBA}.
On using~\eqref{eq:B}, \eqref{eq:g_quad} becomes
\begin{equation}\label{eq:g_quad_alpha}
\begin{split}
g(\vm{\lambda}) =
\ & \lambda_1 \bigl[ \alpha_1(\lambda_2) + \alpha_4 (1-\lambda_4) - \alpha_1(0) \bigr]
+ \lambda_2  \bigl[ \alpha_2(\lambda_3) + \alpha_1 (1 - \lambda_1) - \alpha_2(0) \bigr] \\ 
& + \lambda_3 \bigl[ \alpha_3 (\lambda_4) + \alpha_2 (1 - \lambda_2) - \alpha_3 (0) \bigr] + 
\lambda_4 \bigl[ \alpha_4 (\lambda_1) + \alpha_3(1- \lambda_3) - \alpha_4(0) \bigr] .
\end{split}
\end{equation}
On using $\vm{\lambda}$ from~\eqref{eq:lambda_square} for the unit square, 
one observes that the above equation is similar in form to~\eqref{eq:G_FEbasis}, which in part provides
a connection for the developments 
in~\citep{Randrianarivony:2011:OTI}. However, the Wachspress-based
transfinite interpolant
in~\eqref{eq:g_quad_alpha} 
and the bilinear Coons transfinite interpolant in~\eqref{eq:G_FEbasis} are
distinct over the unit square. 
An illustrative example is presented 
in~\sref{subsec:TFI_verification_polygons}.
Bilinear Coons transfinite interpolation uses affine 
blending functions to form a convex combination of the boundary Dirichlet functions, where the boundary functions remain unmodified, which is 
observable in~\eqref{eq:G_FEbasis}. In contrast, with the 
Wachspress-based transfinite interpolant given 
in~\eqref{eq:g_quad_alpha}, the argument of the boundary functions 
$\alpha_i(\vm{\lambda})$ depend on the 
Wachspress coordinates, which are bilinear functions on the square.

On observing the pattern in~\eqref{eq:g_quad}, it can
be generalized to any convex polygonal domain, such as a triangular domain
or a polygonal one with more than four edges. 
For example,
the expression for $g$ on a triangular domain
becomes
\begin{equation}\label{eq:g_triangle}
\begin{split}
g(\vm{\lambda}) =
\ & \lambda_1 \bigl[ \B(1-\lambda_2,\lambda_2,0)
+ \B(1-\lambda_3,0,\lambda_3) - \B(1,0,0) \bigr] \\
& + \lambda_2 \bigl[ \B(0,1-\lambda_3,\lambda_3)  
+ \B(\lambda_1,1-\lambda_1,0) - \B(0,1,0) \bigr] \\
& + \lambda_3 \bigl[ \B(\lambda_1,0,1-\lambda_1)
+ \B(0,\lambda_2,1 - \lambda_2) - \B(0,0,1) \bigr],
\end{split}
\end{equation}
which can be shown to agree with the Coons triangular transfinite interpolant~\citep{Barnhill:1973:SIT}.  
 
Figure~\ref{fig:TFI_pentagon} shows the Dirichlet boundary functions on a pentagonal domain. 
On this pentagonal domain, $g$ is:
\begin{figure}
\centering
\begin{tikzpicture}[scale=1]
% Define the pentagon's corners
    \coordinate (A) at (0.1,0.4); \filldraw[blue] (A) circle (0.1cm);
    \node[below] at (A) {\Large $\lambda_1$};
    
    \coordinate (B) at (5,0.5); \filldraw[blue] (B) circle (0.1cm);
    \node[below] at (B) {\Large $\lambda_2$};

     \coordinate (C) at (5.75,3.0); \filldraw[blue] (C) circle (0.1cm);
    \node[right] at (C) {\Large $\lambda_3$};
    
    \coordinate (D) at (4,5.5); \filldraw[blue] (D) circle (0.1cm);
    \node[above,yshift=0.05cm] at (D) {\Large $\lambda_4$};
    
    \coordinate (E) at (1,4); \filldraw[blue] (E) circle (0.1cm);
    \node[above,yshift=0.09cm] at (E) {\Large $\lambda_5$};
    
    % Draw the blue pentagon
   \draw[line width=1pt, blue] (A) -- (B) -- (C) -- (D) -- (E) -- cycle;

    % Draw red arrows below each edge
    % Bottom (1st) edge
    \draw[red, -{Stealth[length=3mm, width=3mm]}] (2.0,0.1) -- (3.0,0.120408);
    \node at (2.5,-0.44) {{\Large $\alpha_1(\lambda_2)$}};
    % 2nd edge
    \draw[red, -{Stealth[length=3mm, width=3mm]}] (5.63,1.2) -- (5.9,2.1);
    \node at (6.7,1.55) {{\Large $\alpha_2(\lambda_3)$}};
    % 3rd edge
    \draw[red, -{Stealth[length=3mm, width=3mm]}] (5.6,3.85) -- (5.005,4.7);
    \node at (6.3,4.4) {{\Large $\alpha_3(\lambda_4)$}};
    % 4th  edge
    \draw[red, -{Stealth[length=3mm, width=3mm]}] (2.8,5.35)
    -- (1.9,4.9);
    \node at (2.2,5.8) {{\Large $\alpha_4(\lambda_5)$}};
    % 5th edge
    \draw[red, -{Stealth[length=3mm, width=3mm]}] (0.28,2.6) -- (0.03,1.6); 
    \node at (-0.8,2.2) {{\Large $\alpha_5(\lambda_1)$}};
    % write P
    \node at (3.02,2.68) {{\Large $P$}};
\end{tikzpicture}
\caption{Boundary functions for Wachspress-based transfinite interpolation 
on a pentagonal domain.}\label{fig:TFI_pentagon}
\end{figure}
\begin{equation*}
\begin{split}
g(\vm{\lambda}) =
\ & \lambda_1 \bigl[  \B(1-\lambda_2,\lambda_2,0,0,0)
+ \B(1-\lambda_5,0,0,0,\lambda_5) - \B(1,0,0,0,0) \bigr] \\
& + \lambda_2 \bigl[ \B(0,1-\lambda_3,\lambda_3,0,0)  
+ \B(\lambda_1,1-\lambda_1,0,0,0) - \B(0,1,0,0,0) \bigr] \\
& + \lambda_3 \bigl[  \B(0,0,1-\lambda_4,\lambda_4,0)
+ \B(0,\lambda_2,1-\lambda_2,0,0) - \B(0,0,1,0,0) \bigr] \\
& + \lambda_4 \bigl[  \B(0,0,0,1-\lambda_5,\lambda_5)
+  \B(0,0,\lambda_3,1-\lambda_3,0) - \B(0,0,0,1,0) \bigr] \\
& + \lambda_5 \bigl[ \B(\lambda_1,0,0,0,1-\lambda_1)
+ \B(0,0,0,\lambda_4,1-\lambda_4)
 - \B(0,0,0,0,1) \bigr] ,
\end{split}
\end{equation*}
which in terms of $\{\alpha_i(\cdot)\}_{i=1}^5$ is given by
\begin{equation}\label{eq:g_pentagon}
\begin{split}
g(\vm{\lambda}) =
\ & \lambda_1 \bigl[ \alpha_1 (\lambda_2) + 
\alpha_5 (1-\lambda_5) - \alpha_1(0) \bigr]
 + \lambda_2 \bigl[ \alpha_2 (\lambda_3) 
+ \alpha_1(1-\lambda_1) - \alpha_2(0) \bigr] \\
& + \lambda_3 \bigl[ \alpha_3 (\lambda_4)
+ \alpha_2(1-\lambda_2) - \alpha_3(0) \bigr]
 + \lambda_4 \bigl[ \alpha_4( \lambda_5) 
+  \alpha_3(1-\lambda_3) - \alpha_4(0) \bigr] \\
& + \lambda_5 \bigl[ \alpha_5(\lambda_1)
+ \alpha_4(1-\lambda_4) - \alpha_5(0) \bigr] .
\end{split}
\end{equation}
For a domain that is bounded by a 
convex $n$-gon, $g$ 
can be expressed in 
compact form as: 
\begin{equation}\label{eq:g_polygon}
g(\vm{\lambda}) = \liftTFI [
\vm{\alpha}
(\vm{\lambda} ) ]  =
\sum_{i = 1}^n \lambda_i \Bigl[ 
\alpha_i(\lambda_{i+1}) + 
\alpha_{i-1}(1 \!-\!\lambda_{i-1}) - \alpha_i(0) \Bigr] .
\end{equation}
\subsection{Numerical computations of the transfinite interpolant}
\label{subsec:TFI_verification_polygons}
\begin{figure}[!tbh]
\centering
\begin{subfigure}{0.48\textwidth}
\includegraphics[width=\textwidth]{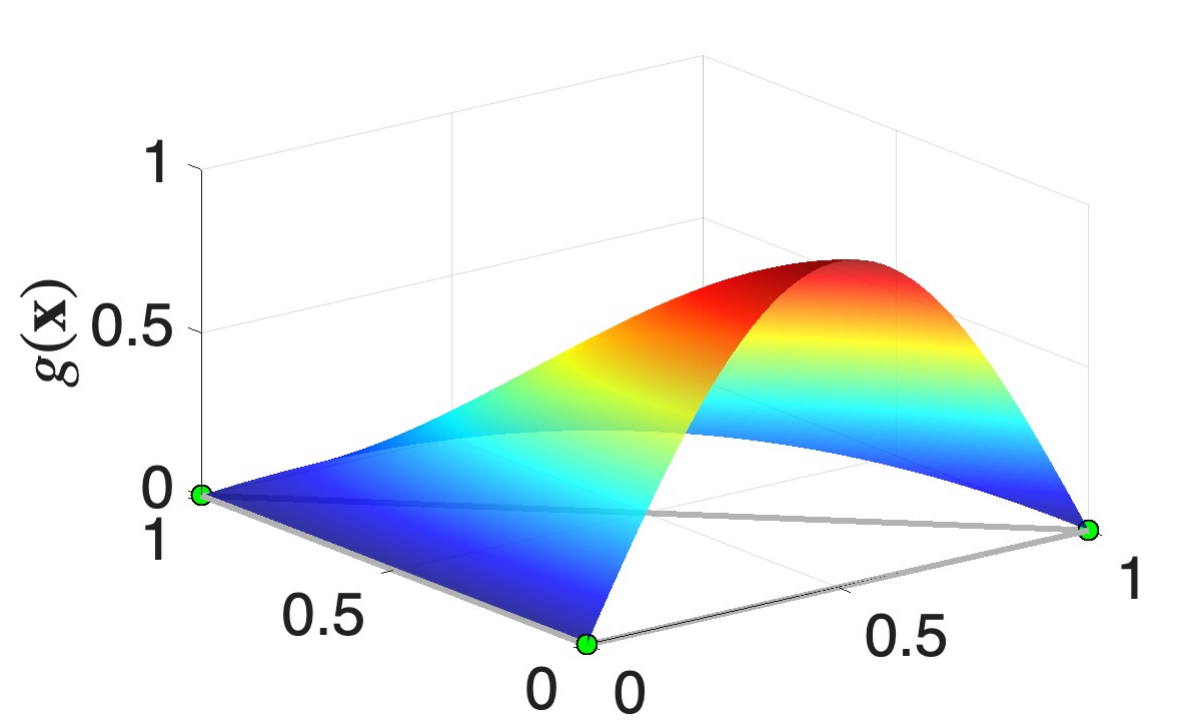}
\subcaption{}\label{fig:gfunction_tri}
\end{subfigure} \hfill
\begin{subfigure}{0.48\textwidth}
\includegraphics[width=\textwidth]{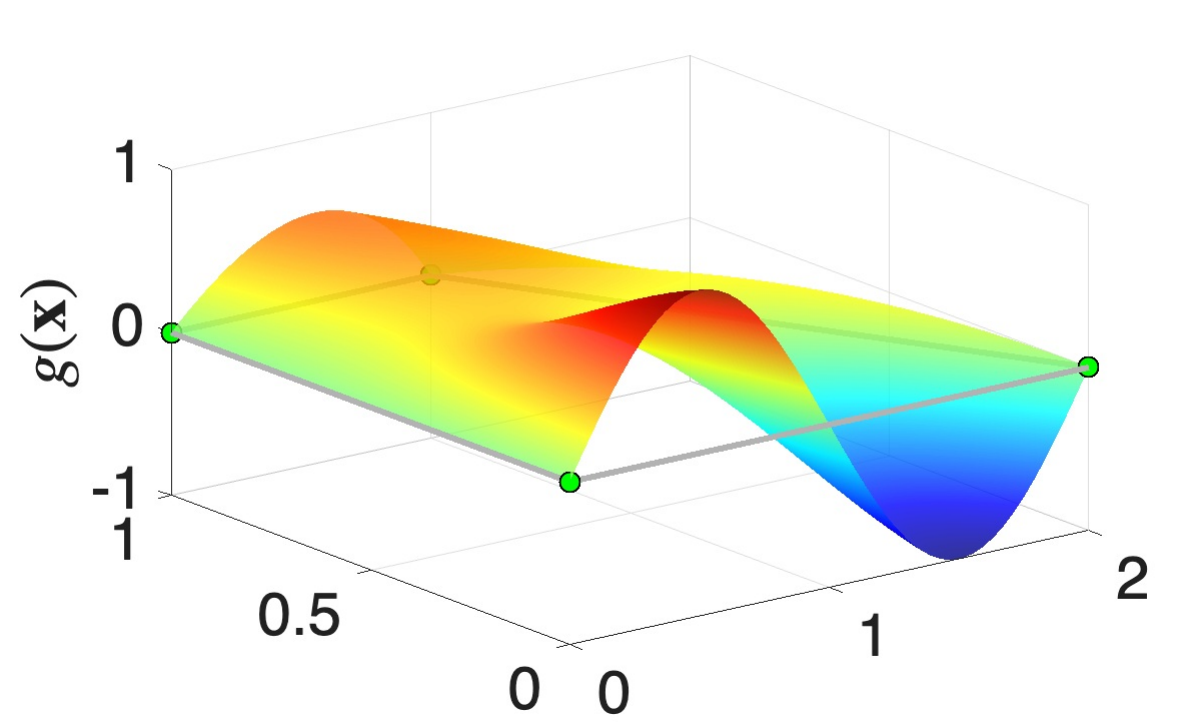}
\subcaption{}\label{fig:gfunction_quad}
\end{subfigure} \hfill
\begin{subfigure}{0.48\textwidth}
\includegraphics[width=\textwidth]{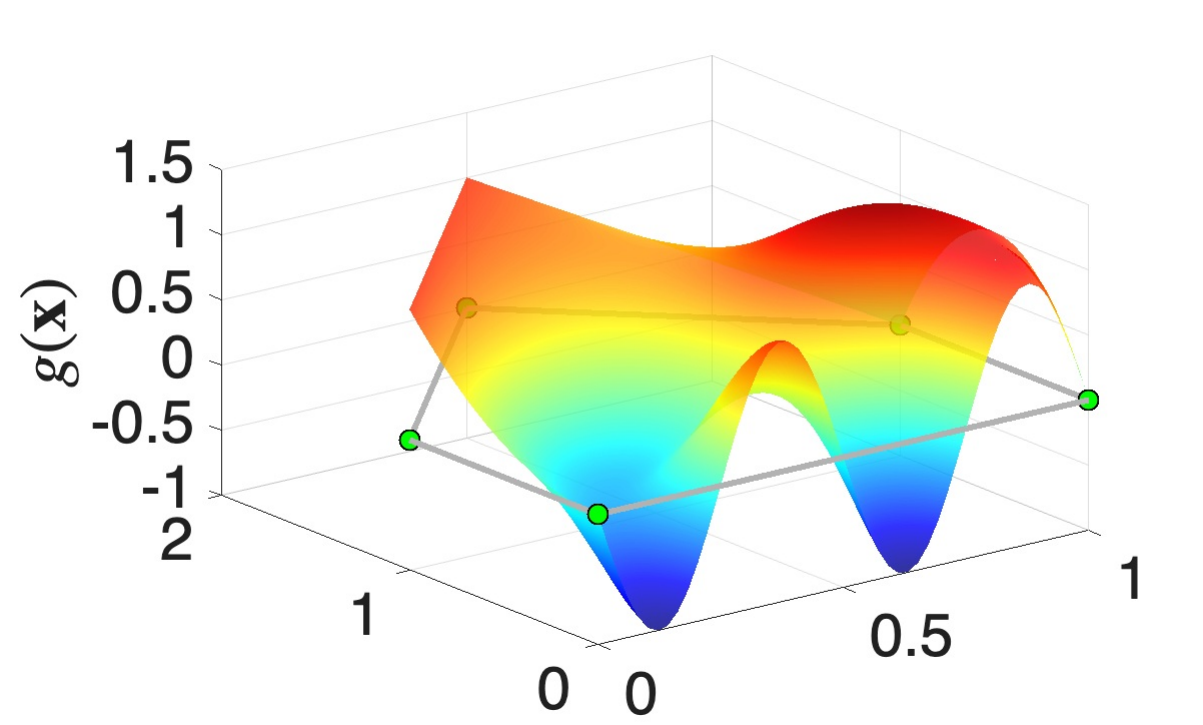}
\subcaption{}\label{fig:gfunction_pentagon}
\end{subfigure} \hfill
\begin{subfigure}{0.48\textwidth}
\includegraphics[width=\textwidth]{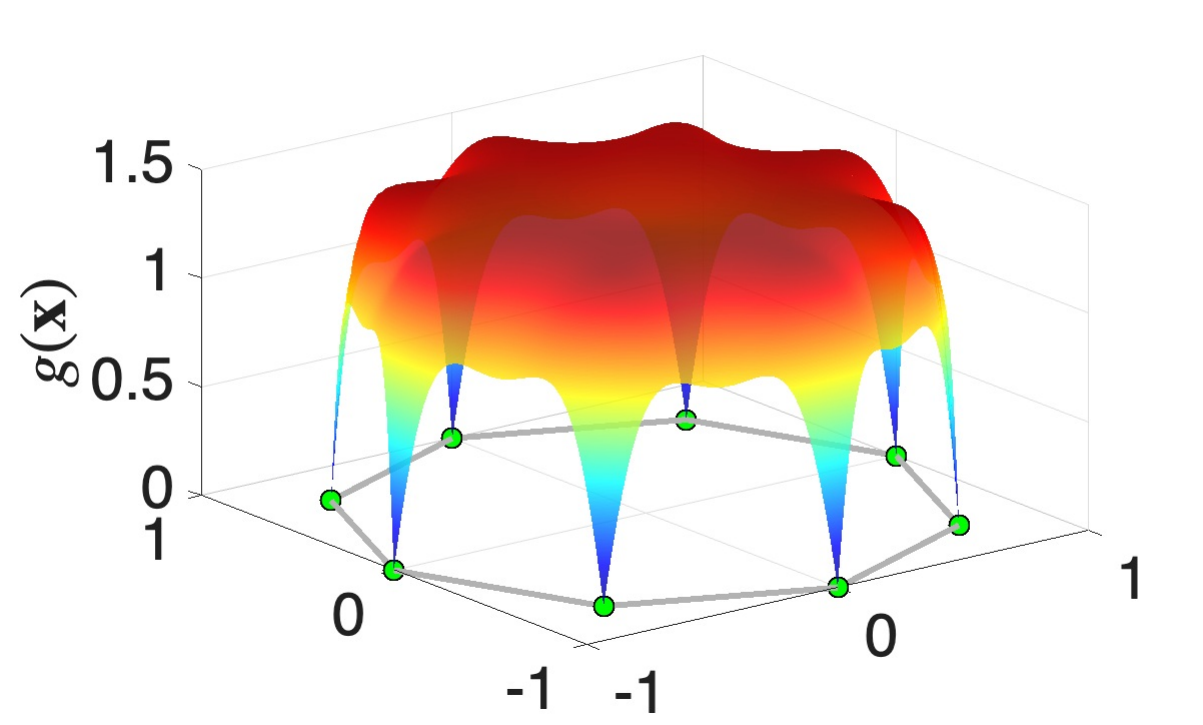}
\subcaption{}\label{fig:gfunction_octagon}
\end{subfigure} \hfill
\caption{Plots of Wachspress-based transfinite interpolant on
(a) triangular, (b) quadrilateral, (c) pentagonal and (d) octagonal domains.}
\label{fig:gfunction}
\end{figure}

Wachspress coordinates are used to compute the transfinite 
interpolant over triangular, quadrilateral, pentagonal, and octagonal domains.
The vertices of the triangle are:
$(0,0)$, $(1,0)$ and $(0,1)$. The Dirichlet boundary functions,
$\{\talpha_i\}_{i=1}^3$,
on the boundary edges 
of the triangle are chosen as:
\begin{equation*}
\talpha_1 (\vx) = \sin (\pi x), \quad
\talpha_2 (\vx) = y(1-y), \quad \talpha_3(\vx) = 0.
\end{equation*}
On using~\eqref{eq:g_triangle}, the transfinite interpolant 
over the triangle is determined:
\begin{equation*}
g(\vx) = x \sin [ \pi (x + y) ] 
- (x + y - 1) \sin ( \pi x ) - x y(x + y - 2) .  
\end{equation*}
The vertices of the quadrilateral are:
$(0,0)$, $(2,0)$, $(1,1)$ and $(0,1)$, and 
$\{\talpha_i\}_{i=1}^4$ are selected as: 
\begin{equation*}
\talpha_1 (\vx) = \sin (\pi x), \quad 
\talpha_2 (\vx) = (x - 1)(2 - x), \quad 
\talpha_3 (\vx) = x \bigl( e - e^x \bigr), \quad
\talpha_4 (\vx) = 0.
\end{equation*}
The exact formula for Wachspress coordinates
on the quadrilateral is obtained by solving~\eqref{eq:wsp_quad}:
\begin{equation}\label{eq:wsp_quad_coordinates}
\vm{\lambda} (\vx) = 
\begin{Bmatrix}
 \dfrac{(1-y) (2 - x - y)}{2-y}, & &\dfrac{x (1-y)}{2-y}, & &\dfrac{xy}{2-y}, && \dfrac{y (2 - x - y)}{2 - y} 
\end{Bmatrix}^\top ,
~\end{equation}                 
and plots of these Wachspress coordinates are presented in~\fref{fig:wsp_quad_1234}. Equation~\eqref{eq:g_quad_alpha} is 
used to determine the transfinite interpolant.
The vertices of the pentagon are:
$(0,0)$, $(1,0)$, $(1,1), (1/2,2), (0,1)$. The Dirichlet boundary functions, $\{\talpha_i(\vx)\}_{i=1}^5$, are chosen as:
\begin{equation*}
\tilde{\alpha}_1(\vx) = -\sin (4 \pi x), \quad
\tilde{\alpha}_2(\vx) = 4y(1-y), \quad
\tilde{\alpha}_3(\vx) = y-1, \quad 
\tilde{\alpha}_4(\vx) = 1, \quad 
\tilde{\alpha}_5(\vx) = y^2.
\end{equation*}
Lastly, we compute the
transfinite interpolant for a regular octagon that is 
inscribed in a unit circle with center
at the origin.  The function 
$- \sin [ 4 \pi ( x^2 + y^2) ] $
is prescribed on the boundary.  
For the pentagon and regular octagon, the 
Wachspress coordinates are computed
using~\eqref{eq:wsp}, and their
plots are shown in Figs.~\ref{fig:wsp_pentagon_12345}
and~\ref{fig:wsp_octagon_12345678}, respectively.
On using~\eqref{eq:g_pentagon} for the pentagon
and~\eqref{eq:g_polygon} for the octagon, the transfinite 
interpolants are computed. 
The transfinite interpolants over the triangular,
quadrilateral, pentagonal and octagonal domains
are plotted in~\fref{fig:gfunction}, and we observe that the Dirichlet boundary
functions are captured.

\begin{figure}[!tbh]
\centering
\begin{subfigure}{0.48\textwidth}
\includegraphics[width=\textwidth]{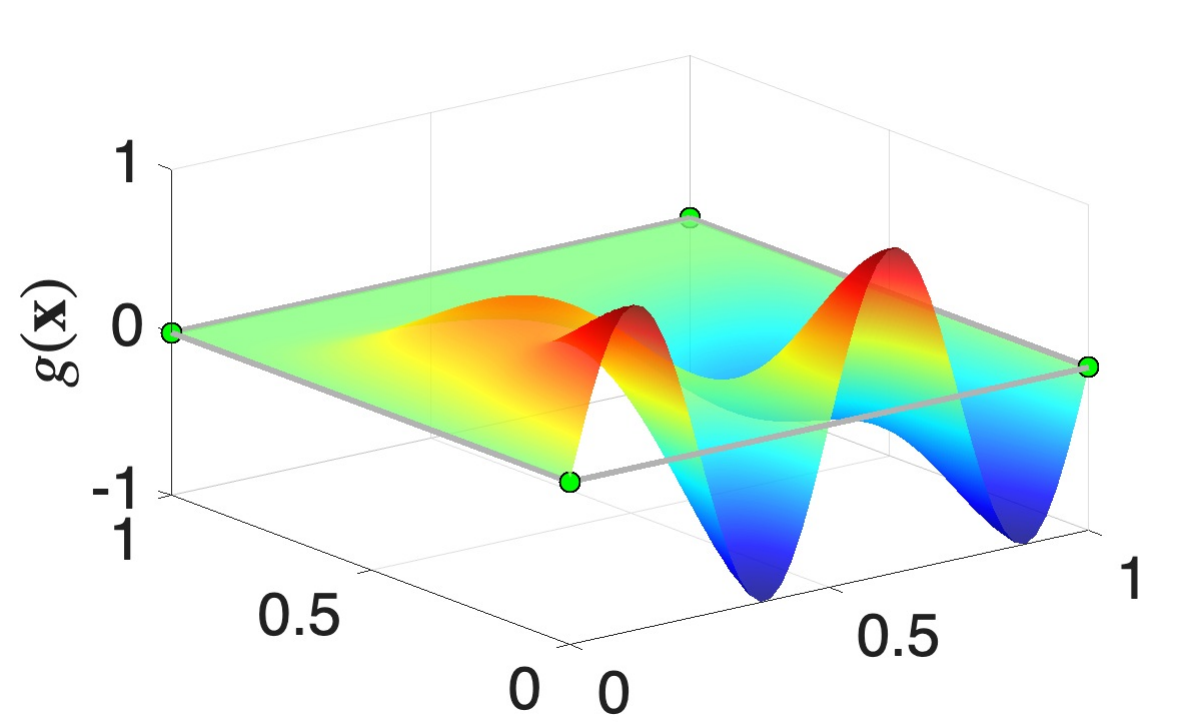}
\subcaption{}\label{fig:gfunction_mvc_a}
\end{subfigure} \hfill
\begin{subfigure}{0.48\textwidth}
\includegraphics[width=\textwidth]{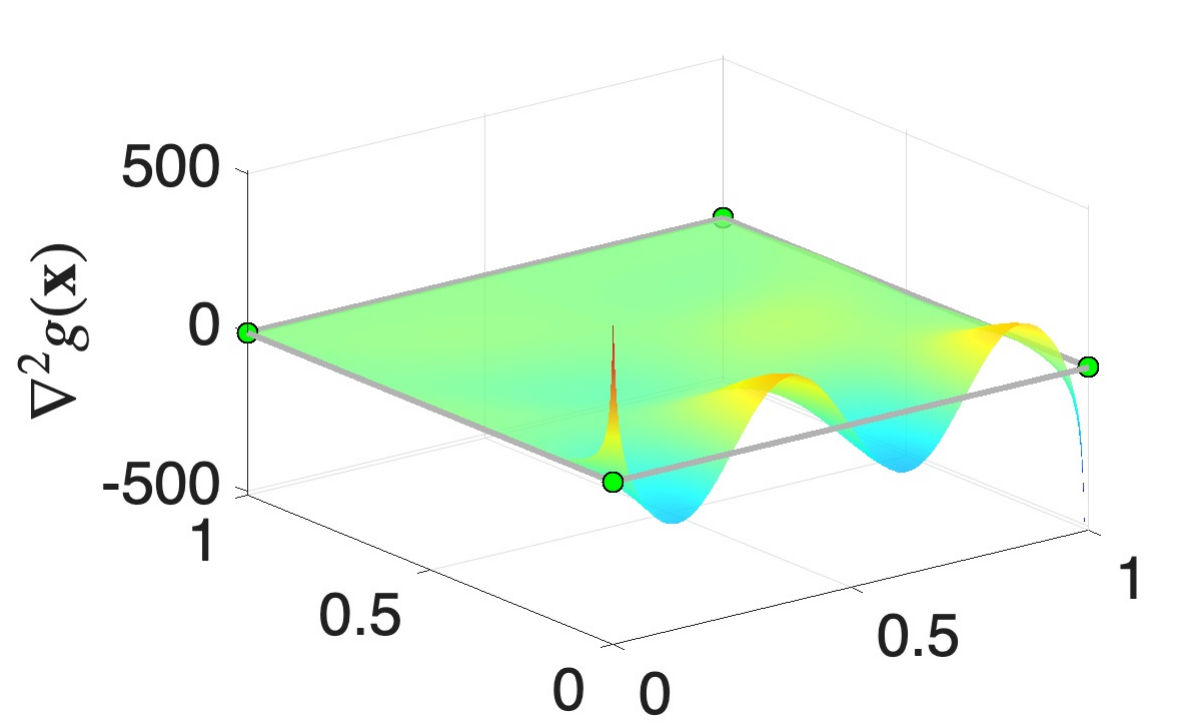}
\caption{}\label{fig:gfunction_mvc_b}
\end{subfigure} 
\caption{Behavior of transfinite interpolant using
mean value coordinates on the unit square. (a) Interpolant $g$ and (b)  Laplacian of $g$. The function $\sin(4 \pi x)$ is prescribed on the bottom edge and zero is imposed on the other edges. The interpolant is well-behaved, but 
$\nabla^2 g$ diverges at two
of the vertices that belong to the edge on which
nonzero Dirichlet boundary conditions are prescribed.}
\label{fig:gfunction_mvc}
\end{figure}

The behavior of the transfinite interpolant with mean value coordinates over a square is assessed. On
a square, mean value
coordinates are obtained by solving~\eqref{eq:wsp_quad} with
$\rho_i(\vx) = h_{i}(\vx)$~\citep{Dieci:2023:MC}.
Plots of
$g$ over a square are presented in~\fref{fig:gfunction_mvc}.
Mean value coordinates are $C^0$ at the vertices of a planar polygon~\citep{Floater:2015:GBC}, and we find that 
their Laplacian diverges at the vertices.
The interpolant $g$ is well-behaved in~\fref{fig:gfunction_mvc_a}, but 
$| \nabla^2 g | \to \infty$ at the vertices of the edge with nonzero Dirichlet conditions (see~\fref{fig:gfunction_mvc_b}).  Since the neural 
network trial function that is constructed with the
Wachspress-based transfinite formulation has a 
bounded Laplacian, it is used to solve PDEs in~\sref{sec:numerical_experiments}.

\subsection{Transfinite neural network trial function}\label{subsec:trial_2D}
From~\sref{subsec:TFI_2D}, we recall that
the transfinite interpolant of the Dirichlet 
boundary condition lifts functions from the polygonal boundary to its interior. On following the rationale and steps described in one dimension
(see~\sref{subsec:trialfunction_1D}), we express the transfinite neural network trial function in 2D as the difference between the neural network's output and 
the extension of its boundary restriction into the interior of the domain, 
with $g$ added to it. 
On using
the lifting operator defined in~\eqref{eq:TFI_formula_a} and
the formula for $g$ in~\eqref{eq:g_quad} for a quadrilateral domain,
the lifting operator applied
to the neural network's output yields
\begin{equation} \label{eq:PN_quad}
\begin{split}
\liftTFI [\N (\vm{\lambda};\vm{\theta})] = \ & 
 \lambda_1 \Bigl[ \N(1-\lambda_2,\lambda_2,0,0;\vm{\theta}) +
\N(1-\lambda_4,0,0,\lambda_4; \vm{\theta}) - \N(1,0,0,0; \vm{\theta}) \Bigr] \\
& + \lambda_2 \Bigl[ \N(0,1-\lambda_3,\lambda_3,0; \vm{\theta} )
+ \N(\lambda_1,1-\lambda_1,0,0; \vm{\theta} )
 - \N(0,1,0,0; \vm{\theta} ) \Bigr] \\
& +  \lambda_3 \Bigl[ \N(0,0,1-\lambda_4,\lambda_4; \vm{\theta} )
  + \N(0,\lambda_2,1-\lambda_2,0; \vm{\theta} )
  - \N(0,0,1,0; \vm{\theta} ) \Bigr]  \\
 & + \lambda_4 \Bigl[ \N(\lambda_1,0,0,1-\lambda_1; \vm{\theta} )
   + \N(0,0,\lambda_3,1-\lambda_3; \vm{\theta} )
   - \N(0,0,0,1; \vm{\theta} ) \Bigr] .
\end{split}
\end{equation}
Hence, the transfinite neural network
trial function on a quadrilateral domain
is written as:
\begin{equation}\label{eq:TFI_trial}
\uTFI(\vm{\lambda};\vm{\theta})
= g(\vm{\lambda}) +
\N(\vm{\lambda};\vm{\theta}) - \liftTFI [ \N(\vm{\lambda};\vm{\theta}) ],
\end{equation}
where we point out that 
by construction 
$\N(\vm{\lambda};\vm{\theta}) - \liftTFI [ \N(\vm{\lambda};\vm{\theta})] = 0 $
on $\partial P$. For 
the triangular and pentagonal
domains, $g$ is provided 
in~\eqref{eq:g_triangle} and~\eqref{eq:g_pentagon},
respectively, and the corresponding $\liftTFI$ operator is used
to compute $\liftTFI[\N]$. If homogeneous
Dirichlet boundary conditions are prescribed
on the boundary, then $g(\vm{\lambda}) = 0$ and the
trial function reduces to
\begin{equation}\label{eq:TFI_trial_zero}
\uTFI(\vm{\lambda};\vm{\theta})
= \N(\vm{\lambda};\vm{\theta}) - \liftTFI [ \N(\vm{\lambda};\vm{\theta}) ].
\end{equation}
The Laplacian of~\eqref{eq:TFI_trial} and~\eqref{eq:TFI_trial_zero} 
are bounded at the vertices of a convex polygonal domain
since each Wachspress coordinate $\lambda_i(\vx)$ is $C^\infty$ at the vertices.
Contrast this to the trial function with ADFs,
\begin{equation}\label{eq:ADF_trial_zero}
\uADF(\vx;\vm{\theta}) = \phi(\vx) \Nxt ,
\end{equation}
which appears
in~\sref{sec:unbounded_Laplacian}. On a square domain,
the Laplacian of the ADF to the boundary of the square, 
namely $\nabla^2 \phi$, diverges at the vertices 
(see~\fref{fig:Laplacian-b}) and consequently the Laplacian 
of the trial function in~\eqref{eq:ADF_trial_zero}
is unbounded at the vertices. In this contribution, 
we use the
trial function in~\eqref{eq:TFI_trial} to overcome 
this limitation of the trial function that is based on ADFs~\citep{Sukumar:2022:EIB}.

\subsection{Python implementation of the transfinite 
interpolant}\label{subsec:python}
This section presents 
standalone Python functions for the implementation of the
action of the lifting operator on the boundary function that yields $g$ in~\eqref{eq:g_polygon}. These functions are displayed under 
Listing~\ref{lst:python_script}. For simplicity, a
unit square is chosen on which Wachspress coordinates are computed; the
function \texttt{compute\_wachpress\_square} needs to be replaced for other 
polygons.  For a convex quadrilateral, use the exact formula given 
in~\eqref{eq:wsp_quad} and implement~\eqref{eq:wsp_global} 
for any other convex polygon. 
The transfinite interpolant is computed in function \texttt{compute\_g\_function}, which
is applicable for any convex polygonal domain. 
The Dirichlet boundary conditions are specified in Cartesian coordinates on the edges of the polygon in the function \texttt{boundary\_function}.
As boundary conditions, we have chosen $u = \sin (\pi x)$ on the top edge and zero on
the remaining three edges. 
Since $g$ is a function of the Wachspress coordinates and does not
involve any contribution from the neural network's output, the computation of the Wachspress coordinates and the $g$ function are a one-time operation that is done prior to model training. 
The Python implementation of the lifting operator acting on the neural network's output, which is required to form the transfinite 
trial function $\uTFI(\vm{\lambda};\vm{\theta})$ given in~\eqref{eq:TFI_trial},
mirrors the implementation for $g$.

\begin{lstlisting}[style=mystyle,float=t,language=Python,caption={Python functions to compute the transfinite interpolant $g$.},label={lst:python_script}]
import numpy as np

def compute_wachspress_square(x_pts, y_pts, n_verts):
    """
    Compute bilinear Wachspress coordinates on the unit square
    Vertices: (0,0), (1,0), (1,1), (0,1)
    """
    lambda_wsp = np.zeros((n_verts, x_pts.shape[0]))
    lambda_wsp[0,:] = (1.0 - x_pts) * (1.0 - y_pts)
    lambda_wsp[1,:] = x_pts * (1.0 - y_pts)
    lambda_wsp[2,:] = x_pts * y_pts
    lambda_wsp[3,:] = (1.0 - x_pts) * y_pts    
    return lambda_wsp
 
def compute_g_function(x_verts, y_verts, lambda_wsp):
    """
    Compute transfinite interpolant of Dirichlet boundary functions
    """
    n_verts = x_verts.shape[0]
    g = np.zeros(lambda_wsp.shape[1])
    for i in range(n_verts):
        """
        get the other vertices for edges attached to vertex i
        """  
        ip1 = (i + 1) % n_verts  # next vertex
        im1 = (i - 1) % n_verts  # previous vertex

        """
        (x,y) coordinates where the boundary function will be evaluated
        """
        x_ip1 = (1.0 - lambda_wsp[ip1,:]) * x_verts[i] + lambda_wsp[ip1,:] * x_verts[ip1]
        y_ip1 = (1.0 - lambda_wsp[ip1,:]) * y_verts[i] + lambda_wsp[ip1,:] * y_verts[ip1]
        x_im1 = (1.0 - lambda_wsp[im1,:]) * x_verts[i] + lambda_wsp[im1,:] * x_verts[im1]
        y_im1 = (1.0 - lambda_wsp[im1,:]) * y_verts[i] + lambda_wsp[im1,:] * y_verts[im1]

        """
        Contributions to the transfinite boundary interpolant from this vertex
        """
        g += lambda_wsp[i,:] * (boundary_function(i, x_ip1, y_ip1) \
                                + boundary_function(im1, x_im1, y_im1) \
                                - boundary_function(i, x_verts[i], y_verts[i]))
    return g
 
def boundary_function(edge_idx, x, y):
    """
    Define the boundary function for any given edge
    """
    return np.select([edge_idx == 2], [np.sin(np.pi * x)], default=0)
\end{lstlisting}

\section{Network architecture and model training}\label{sec:model_training}
\begin{figure}
\centering
\includegraphics[width=0.7\textwidth]{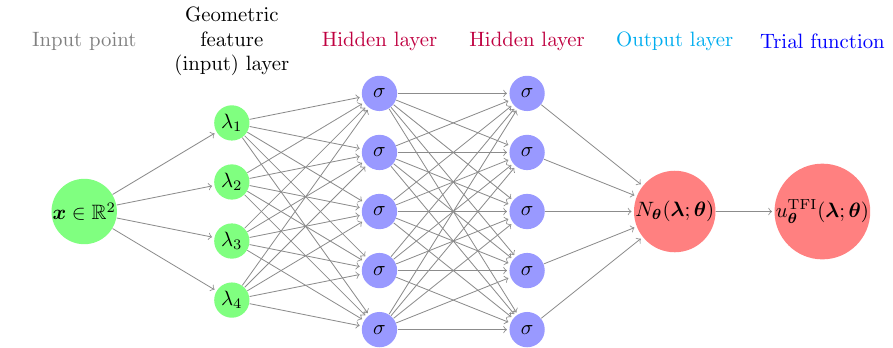}
\medskip
\caption{Schematic of neural network architecture 
$4$--$5$--$5$--$1$ on a 
quadrilateral domain ($n = 4$). On
solving~\eqref{eq:wsp_quad}, 
an input point $\vx \in \bar \Omega$ 
is mapped to $\vm{\lambda} \in [0,1]^4$, which are the
Wachspress coordinates on a quadrilateral.
The $\vm{\lambda}$ appear in the geometric
feature (map) layer of the neural network. On using
$\NLt$ and $g(\vm{\lambda})$ in~\eqref{eq:TFI_trial}, the trial
function $\uLtTFI$ is computed.}
\label{fig:dnn}
\end{figure}

In this paper, we use deep neural networks
to solve linear and nonlinear
Poisson problems over polygonal domains.
The network architecture is a fully connected 
feed-forward multilayer perceptron. The neural network consists of an input layer, $L$ hidden layers with 
${\cal N}_\ell$
neurons in the $\ell$-th hidden layer, and an
output layer with one neuron.  For all the numerical
experiments, we choose
the tanh activation function, 
$\sigma(z) = \tanh(z)$.  For a polygonal domain $\Omega$,
a point $\vx \in \bar \Omega$ is mapped 
to Wachspress coordinates, 
$\vm{\lambda}: \bar \Omega \to [0,1]^n$, which
serve as a geometric feature map and reside in
the geometric feature (input) layer. 
Then, $\vm{\lambda}$ is passed to the next (hidden) layer in the
neural network.
Let $\NLt$ be the neural
network output, with
$\vm{\theta} := \{ \vm{W}, \vm{b} \}$ the unknown parameter vector that consists of
weights $\vm{W}_\ell \in \Re^{ {\cal N}_\ell \times { \cal N}_{\ell-1}}$ and 
biases $\vm{b}_\ell \in \Re^{ {\cal N}_\ell }$. 
We express $\NLt$ via the composition of $T^{(\ell)}$
($\ell=1,2,\ldots,L)$ and a linear map ${\cal G}$ as:
\begin{equation}\label{eq:composition}
\NLt = {\cal G} \circ T^{(L)} \circ T^{(L-1)} \circ \ldots \circ 
T^{(1)}(\vm{x} \mapsto \vm{\lambda}),
\end{equation}
where ${\cal G}: \Re^{{\cal N}_L} \rightarrow \Re$ is the linear mapping
for the output layer and in each hidden layer ($\ell = 1,2,\ldots,L$), the
nonlinear mapping is:
\begin{equation}\label{eq:Tmap}
T^{(\ell)}(\vm{z}) = \sigma(\vm{W}_\ell \cdot \vm{z} + \vm{b}_\ell ),
\end{equation}
where $\vm{z} \in \Re^{{\cal N}_{\ell-1}}$. A schematic of the network architecture for a quadrilateral
domain ($n = 4$) is shown
in~\fref{fig:dnn}.

The formulation to strongly enforce
Dirichlet boundary conditions that is
described in~\sref{sec:formulation} has been implemented for PINNs using PyTorch~\citep{Paszke:2016:PYT}.  
All simulations are executed on an NVIDIA H100 GPU.
For all problems, two to four hidden layers are used
with either 20, 30 or 40 neurons in each hidden layer. 
For improved accuracy in PINNs,
use of self-adaptive weight training~\citep{McClenny:2023:SAP} has been
adopted in many prior studies using PINNs. However, very recently, algorithms have emerged (based on 
quasi-Newton methods and natural gradient)
that report significantly better accuracy in
neural network training,
cf.~\citep{Rathore:2024b:CTP,Urban:2025:UOP,Jnini:2025:DNG,Kiyani:2025:OOP,Wang:2025:GAP}. Most of these 
approaches advocate an
initial training phase with the Adam optimizer, 
subsequent switch to L-BFGS and  Nystr{\"o}m Newton conjugate gradient or use of self-scaled Broyden algorithm to 
further decrease the loss.  The rationale for this sequence stems
from the properties of stochastic gradient-based and Newton-based
algorithms: the former is able to `escape' 
saddle points in a complex landscape, whereas the latter can be attracted to saddle points. Once the loss has attained sufficient decrease with Adam 
in the neighborhood of a good local minima, L-BFGS and other
quasi-Newton
methods can enable (when Hessian of loss function is 
well-conditioned) convergence to a much smaller loss. In this paper, we adopt the simple modification proposed in~\citep{Urban:2025:UOP} of using the $\log$ of the loss with L-BFGS, where it has been shown to
deliver accuracies that are comparable to that obtained using SSBroyden. In the computations that follow, 
a learning rate of 
$10^{-3}$ 
is used
for Adam. The strong Wolfe condition
is used in the line
search for L-BFGS. 
For the highly oscillatory harmonic
problem in~\sref{subsec:laplace}, mini-batch training is used; single batch training is used for all other problems. Apart from selecting the number of epochs for 
Adam and L-BFGS, we do not perform any 
hyperparameter tuning to produce the PINN and deep Ritz solutions.

\section{Numerical experiments}\label{sec:numerical_experiments}

\subsection{Harmonic problems on the unit square}\label{subsec:laplace}
In~\citep{Sukumar:2022:EIB}, the ADF to the boundary of the square diverged at the vertices of the square, and therefore collocation points in PINNs had to be chosen away from the boundary of the square.  Here, we first revisit
a Laplace problem from~\citep{Sukumar:2022:EIB} to compare PINN solutions using ADFs and the Wachspress-based transfinite 
trial function that is proposed in this paper. 
We solve $\nabla^2 u = 0$ in
$\Omega = (0,1)^2$ with the
Dirichlet boundary conditions~\citep{Sukumar:2022:EIB}: 
\begin{equation}\label{eq:laplace}
u(\vx) = 0 \ \ \textrm{on } \Gamma_1, \ \Gamma_2, \ \Gamma_3, \quad
u(\vx) = \sin \pi x \ \ \textrm{on } \Gamma_4,
\end{equation}
where the boundaries edges are: 
$\Gamma_1 = \{ (x,y) : x = 0, \ 0 \le y \le 1 \}$,  
$\Gamma_2 = \{ (x,y) : 0 \le x \le 1, \ y = 0 \}$,
$\Gamma_3 = \{ (x,y) : x = 1, \ 0 \le y \le 1\}$, and
$\Gamma_4 = \{ (x,y) : 0 \le x \le 1, \ y = 1\}$.
The exact solution for this boundary-value problem is:
\begin{equation}\label{eq:laplace_exact}
u(\vx) = \frac{ \sinh( \pi y ) \sin \pi x }
              { \sinh (\pi) } .
\end{equation}
On noting that $\alpha_1(\lambda_2) 
= \alpha_2(\lambda_3) = \alpha_4(\lambda_1) = 0$, 
$\alpha_3(\lambda_4) = 
\sin [ \pi (1 - \lambda_4)]$,  and 
using~\eqref{eq:g_polygon}, we can write $g$ as:
\begin{equation}\label{eq:g_sinebc}
g(\vm{\lambda}) = \lambda_3 
\sin [ \pi  (1 -  \lambda_4 ) ] + 
 \lambda_4 \sin ( \pi \lambda_3 ) = 
\lambda_3 \sin (\pi \lambda_4) +
\lambda_4 \sin (\pi \lambda_3) ,
\end{equation}
and then the transfinite trial function from~\eqref{eq:TFI_trial}
is given by
\begin{equation}\label{eq:trial_laplace_one}
\uTFI(\vm{\lambda};\vm{\theta})
= \lambda_3 \sin (\pi \lambda_4) 
+ \lambda_4 
\sin ( \pi \lambda_3 )  + 
\N(\vm{\lambda};\vm{\theta}) - \liftTFI [ \N(\vm{\lambda};\vm{\theta}) ], 
\end{equation}
where $\liftTFI [\N(\vm{\lambda};\vm{\theta}) ]$
is given in~\eqref{eq:PN_quad}.  
On using~\eqref{eq:lambda_square}, $g$ in~\eqref{eq:g_sinebc} becomes $ xy \sin [ \pi y (1 - x) ] + y(1-x) \sin(\pi x y)$, whereas
the Coons transfinite interpolant from~\eqref{eq:G_Coons} is
$y \sin (\pi x)$.  This shows that in general the 
Wachspress-based transfinite interpolant and the Coons transfinite
interpolant are distinct on a square.

\begin{figure}
\centering
\begin{subfigure}{0.23\textwidth}
\includegraphics[width=\textwidth]{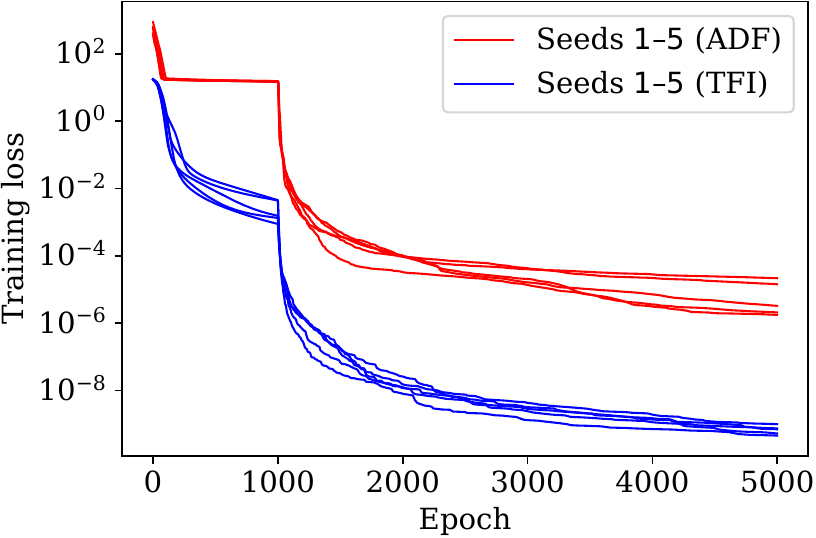}
\subcaption{$\delta = 10^{-2}$}\label{fig:harmonic_square_a}
\end{subfigure} \hfill
\begin{subfigure}{0.23\textwidth}
\includegraphics[width=\textwidth]{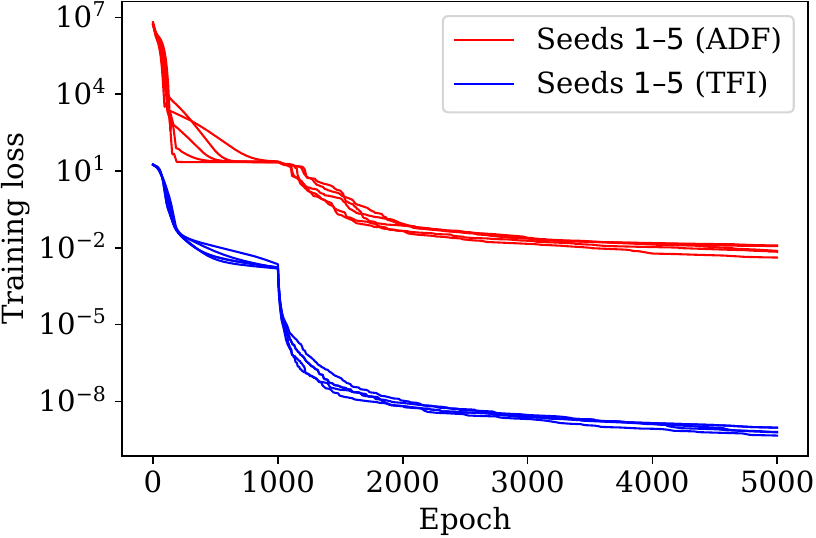}
\subcaption{$\delta = 10^{-4}$}\label{fig:harmonic_square_b}
\end{subfigure} \hfill
\begin{subfigure}{0.23\textwidth}
\includegraphics[width=\textwidth]{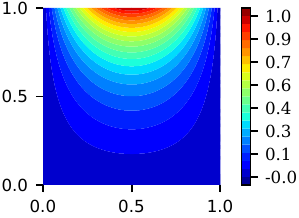}
\subcaption{$u$}\label{fig:harmonic_square_c}
\end{subfigure} \hfill
\begin{subfigure}{0.23\textwidth}
\includegraphics[width=\textwidth]{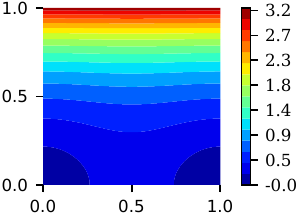}
\subcaption{$|\nabla u|$}\label{fig:harmonic_square_d}
\end{subfigure} 
\begin{subfigure}{0.23\textwidth}
\includegraphics[width=\textwidth]{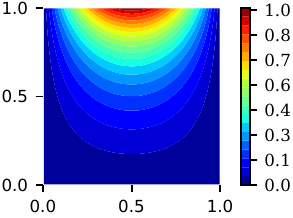}
\subcaption{$\uADF$}\label{fig:harmonic_square_e}
\end{subfigure} \hfill
\begin{subfigure}{0.23\textwidth}
\includegraphics[width=\textwidth]{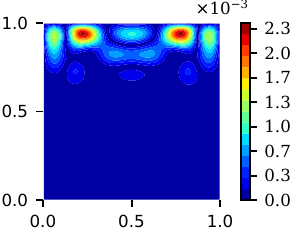}
\subcaption{ $| u - \uADF|$ } \label{fig:harmonic_square_f}
\end{subfigure} \hfill
\begin{subfigure}{0.23\textwidth}
\includegraphics[width=\textwidth]{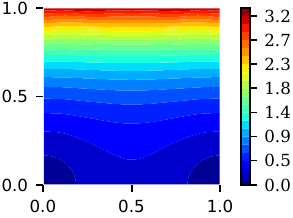}
\subcaption{$| \nabla \uADF |$}\label{fig:harmonic_square_g}
\end{subfigure} \hfill
\begin{subfigure}{0.23\textwidth}
\includegraphics[width=\textwidth]{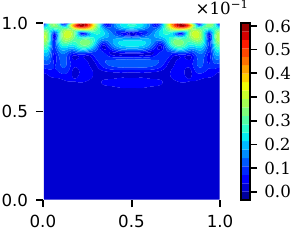}
\subcaption{ $ | \nabla ( u  - \uADF ) |    $ } \label{fig:harmonic_square_h}
\end{subfigure} 
\begin{subfigure}{0.23\textwidth}
\includegraphics[width=\textwidth]{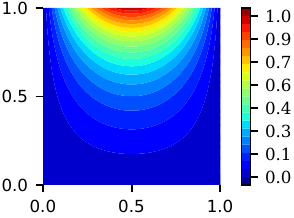}
\subcaption{$\uTFI$}\label{fig:harmonic_square_i}
\end{subfigure} \hfill
\begin{subfigure}{0.23\textwidth}
\includegraphics[width=\textwidth]{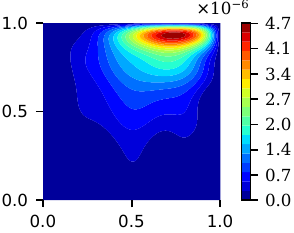}
\subcaption{ $| u - \uTFI|$ } \label{fig:harmonic_square_j}
\end{subfigure} \hfill
\begin{subfigure}{0.23\textwidth}
\includegraphics[width=\textwidth]{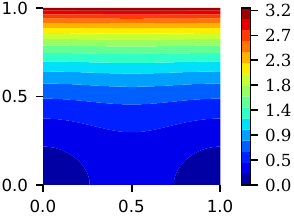}
\subcaption{$| \nabla \uTFI |$}\label{fig:harmonic_square_k}
\end{subfigure} \hfill
\begin{subfigure}{0.23\textwidth}
\includegraphics[width=\textwidth]{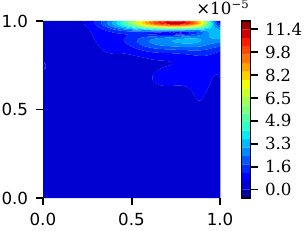}
\subcaption{ $ | \nabla ( u  - \uTFI ) |    $ } \label{fig:harmonic_square_l}
\end{subfigure} 
\begin{subfigure}{0.23\textwidth}
\includegraphics[width=\textwidth]{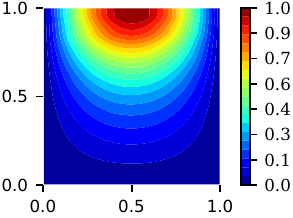}
\subcaption{$\uADF$}\label{fig:harmonic_square_m}
\end{subfigure} \hfill
\begin{subfigure}{0.23\textwidth}
\includegraphics[width=\textwidth]{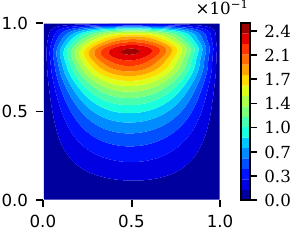}
\subcaption{ $| u - \uADF|$ } \label{fig:harmonic_square_n}
\end{subfigure} \hfill
\begin{subfigure}{0.23\textwidth}
\includegraphics[width=\textwidth]{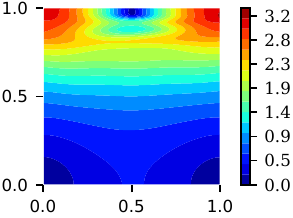}
\subcaption{$| \nabla \uADF |$}\label{fig:harmonic_square_o}
\end{subfigure} \hfill
\begin{subfigure}{0.23\textwidth}
\includegraphics[width=\textwidth]{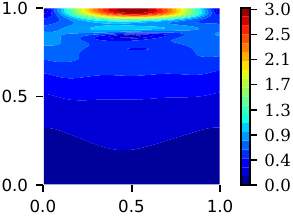}
\subcaption{ $ | \nabla ( u  - \uADF ) |    $ } \label{fig:harmonic_square_p}
\end{subfigure} 
\begin{subfigure}{0.23\textwidth}
\includegraphics[width=\textwidth]{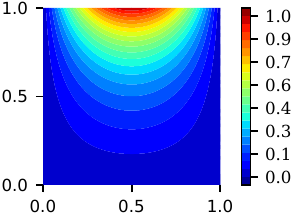}
\subcaption{$\uTFI$}\label{fig:harmonic_square_q}
\end{subfigure} \hfill
\begin{subfigure}{0.23\textwidth}
\includegraphics[width=\textwidth]{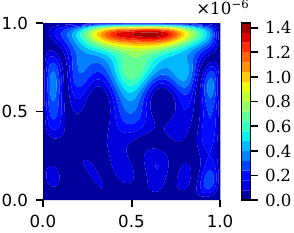}
\subcaption{ $| u - \uTFI|$ } \label{fig:harmonic_square_r}
\end{subfigure} \hfill
\begin{subfigure}{0.23\textwidth}
\includegraphics[width=\textwidth]{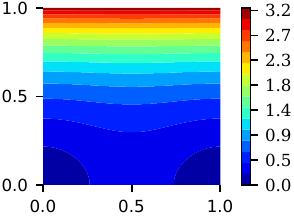}
\subcaption{$| \nabla \uTFI |$}\label{fig:harmonic_square_s}
\end{subfigure} \hfill
\begin{subfigure}{0.23\textwidth}
\includegraphics[width=\textwidth]{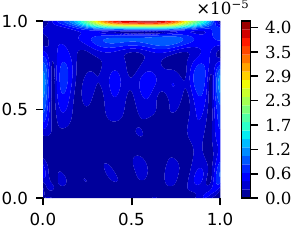}
\subcaption{ $ | \nabla ( u  - \uTFI ) |    $ } \label{fig:harmonic_square_t}
\end{subfigure} 

\caption{PINN solutions for the harmonic boundary-value problem 
in~\protect\eqref{eq:laplace} using ADF
and TFI. Network architecture is 4--20--20--1. 
Training loss using ADF~\citep{Sukumar:2022:EIB} and TFI 
for (a) $\delta = 10^{-2}$ and (b) $\delta = 10^{-4}$. 
(c), (d)
Exact solutions for $u$ and $|  \nabla u |$.
Plots in (e)--(l) are for 100 equispaced training points in 
$[0.01,0.99]^2$ ($\delta = 10^{-2}$).
Absolute 
error in $u$ and $L^2$ norm of the
gradient of the error for
(e)--(h) ADF and (i)--(l) TFI.
Plots in (m)--(t) are for 100 equispaced training points
in $[0.0001,0.9999]^2$ ($\delta = 10^{-4}$).
Absolute 
error in $u$ and $L^2$ norm of the
gradient of the error for
(m)--(p) ADF and (q)--(t) TFI.
}\label{fig:harmonic_square}
\end{figure}

The loss function for the Laplace equation with 
collocation-based PINNs is:
\begin{equation}\label{eq:loss_laplace}
    {\cal L } (\vm{\theta}) = 
    \frac{1}{M} \sum_{k=1}^M
    \left[ \nabla^2 \uTFI ( \vx_k \mapsto \vm{\lambda}^k ;\vm{\theta} )
    \right]^2,
\end{equation}
where $M$ collocation points are chosen in the domain
$[\delta, 1-\delta ]^2$, which is in the interior 
of the unit square. The transfinite trial function, 
$ \uTFI(\vm{\lambda};\vm{\theta})$
is given in~\eqref{eq:trial_laplace_one}.  
The loss function in~\eqref{eq:loss_laplace} is used
for training with the Adam optimizer, whereas its 
$\log$ is adopted for training
with the L-BFGS optimizer.
We choose 100 collocation points ($10 \times 10$ grid) for training and consider $\delta = 10^{-2}$ and 
$\delta = 10^{-4}$ in the numerical computations to compare
the solutions obtained with ADF and TFI.  
The network architecture 4--20--20--1 (two hidden
layers) is used. Network training is performed
with Adam (1,000 epochs) + L-BFGS (4,000 epochs).  The numerical results are presented in~\fref{fig:harmonic_square}. 
For the predictions, we used a refined set of $28,812$ testing points.
In Figs.~\ref{fig:harmonic_square_a} and~\ref{fig:harmonic_square_b}, we present 
the training loss curves for ADF and TFI. Training is 
conducted for five
initial seeds.  The final 
training loss for ADF and TFI are ${\cal O}(10^{-6})$
and ${\cal O}(10^{-9})$, respectively. 
Figures~\ref{fig:harmonic_square_c} and~\ref{fig:harmonic_square_d} depict contour plots of
the exact solutions for $u$
and $|\nabla u|$. 
For $\delta = 10^{-2}$, 
Figs.~\ref{fig:harmonic_square_e}--\ref{fig:harmonic_square_h}
show contour 
plots of the ADF solution, $L^2$ norm of its gradient, along
with the absolute error in $u$ and the 
$L^2$ norm of the gradient of the error. The
corresponding plots for TFI are shown 
in~Figs.~\ref{fig:harmonic_square_i}--\ref{fig:harmonic_square_l}.
In alignment with the final training loss, we observe that
the errors in $u$ as well as in $| \nabla u |$ are markedly
better for TFI when compared to ADF (see discussion
in~\sref{subsec:trialfunction_1D}).
For $\delta = 10^{-4}$, the plots for ADF and TFI are depicted
in Figs.~\ref{fig:harmonic_square_m}--\ref{fig:harmonic_square_p}
and 
Figs.~\ref{fig:harmonic_square_q}--\ref{fig:harmonic_square_t},
respectively. Compared to 
$\delta = 10^{-2}$, collocation points for
$\delta = 10^{-4}$ are much closer to the
boundary vertices, which adversely affects training
and accuracy with ADF: predictions by ADF dramatically worsen with maximum pointwise absolute error increasing from
$10^{-3}$ to $10^{-1}$. For TFI, however, the maximum pointwise absolute error
and the maximum pointwise $L^2$ norm of the gradient of the error remain at
 ${\cal O}(10^{-6})$ and
${\cal O}(10^{-5})$, respectively. In fact,
the accuracy with TFI improves with decrease in $\delta$
(compare Figs.~\ref{fig:harmonic_square_j} and~\ref{fig:harmonic_square_l} to
Figs.~\ref{fig:harmonic_square_r} and~\ref{fig:harmonic_square_t}, respectively). 
This establishes that the Wachspress-based transfinite formulation overcomes one of the main limitations of ADFs in PINNs~\citep{Sukumar:2022:EIB} (see also
the discussion in~\sref{sec:unbounded_Laplacian}).

As a second benchmark on the unit square, we consider a test problem to exemplify the benefits of
exact imposition of oscillatory Dirichlet boundary conditions. To this
end, we solve the Laplace
equation over the unit square with boundary conditions that contain
high frequencies.\footnote{The authors
thank Dr.\ Joseph
Bishop (Sandia National Laboratories)
for suggesting this problem.}  It is known that for such boundary conditions, a consequence of the smoothing properties of the Laplacian is that the solution will dampen very quickly away from the boundary edges. For very high frequencies,
the oscillations will be confined to a thin boundary layer and smoothen out sharply to zero in the bulk of the domain. We solve $\nabla^2 u = 0 $ in
$\Omega = (0,1)^2$ and impose the following
Dirichlet boundary conditions on the four edges:
\begin{equation}\label{eq:Laplace_osc_bc}
u(x,0) = u(x,1) = \sin(10 \pi x), \ \
u(0,y) = u(1,y) = \sin(10 \pi y) ,
\end{equation}
which have 5 oscillations on each edge.

\begin{figure}
\centering
\begin{subfigure}{0.28\textwidth}
\includegraphics[width=\textwidth]{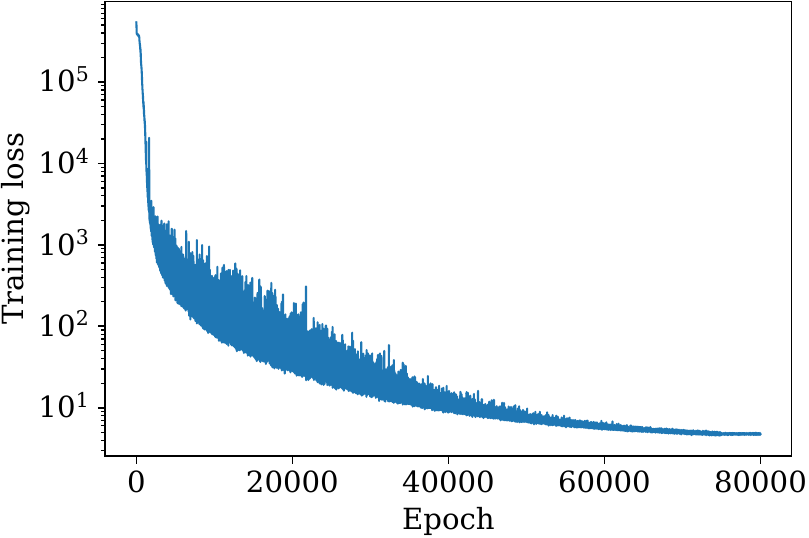}
\subcaption{}\label{fig:harmonic_osc_a}
\end{subfigure} \hfill
\begin{subfigure}{0.29\textwidth}
\includegraphics[width=\textwidth]{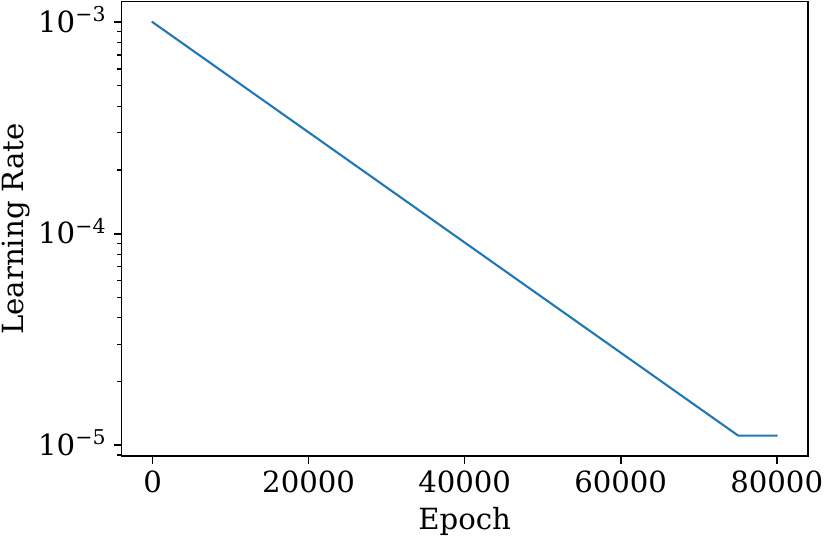}
\subcaption{}\label{fig:harmonic_osc_b}
\end{subfigure} \hfill
\begin{subfigure}{0.2\textwidth}
\includegraphics[width=\textwidth]{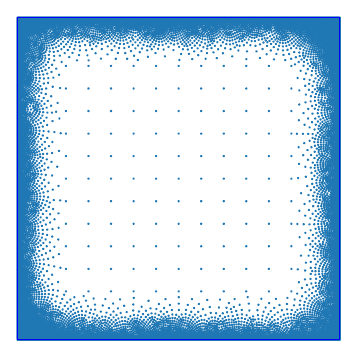}
\subcaption{}\label{fig:harmonic_osc_c}
\end{subfigure} \hfill
\begin{subfigure}{0.2\textwidth}
\includegraphics[width=\textwidth]{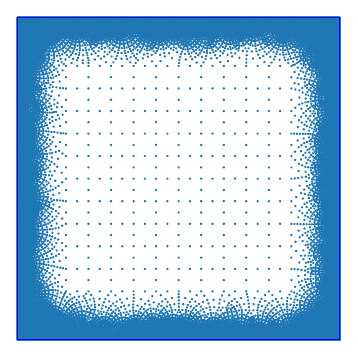}
\subcaption{}\label{fig:harmonic_osc_d}
\end{subfigure}
\begin{subfigure}{0.32\textwidth}
\includegraphics[width=\textwidth]{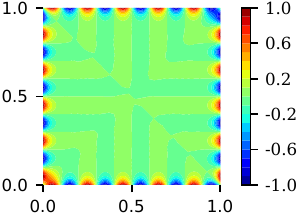}
\subcaption{$u$}\label{fig:harmonic_osc_e} 
\end{subfigure} \hfill
\begin{subfigure}{0.32\textwidth}
\includegraphics[width=\textwidth]{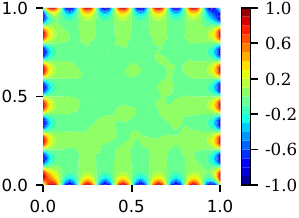}
\subcaption{$\uTFI$}\label{fig:harmonic_osc_f} 
\end{subfigure} \hfill
\begin{subfigure}{0.32\textwidth}
\includegraphics[width=\textwidth]{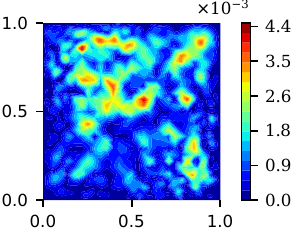}
\subcaption{$| u - \uTFI | $}\label{fig:harmonic_osc_g} 
\end{subfigure} 
\begin{subfigure}{0.3\textwidth}
\includegraphics[width=\textwidth]{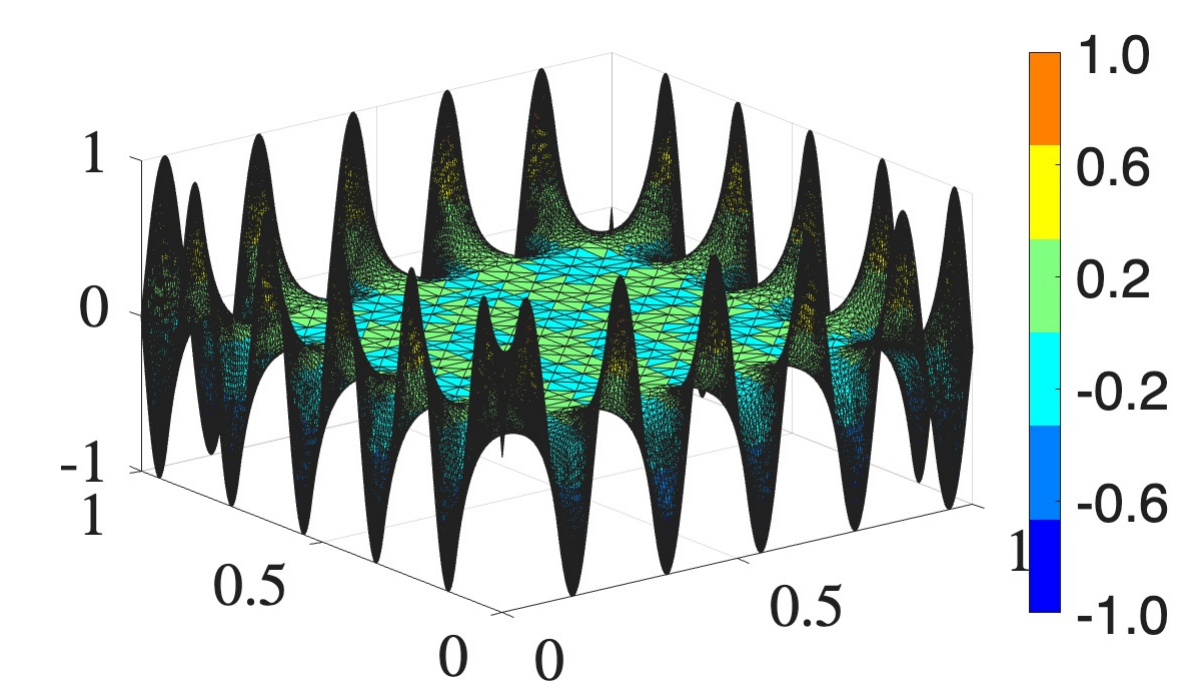}
\subcaption{$u$}\label{fig:harmonic_osc_h}
\end{subfigure} \hfill 
\begin{subfigure}{0.3\textwidth}
\includegraphics[width=\textwidth]{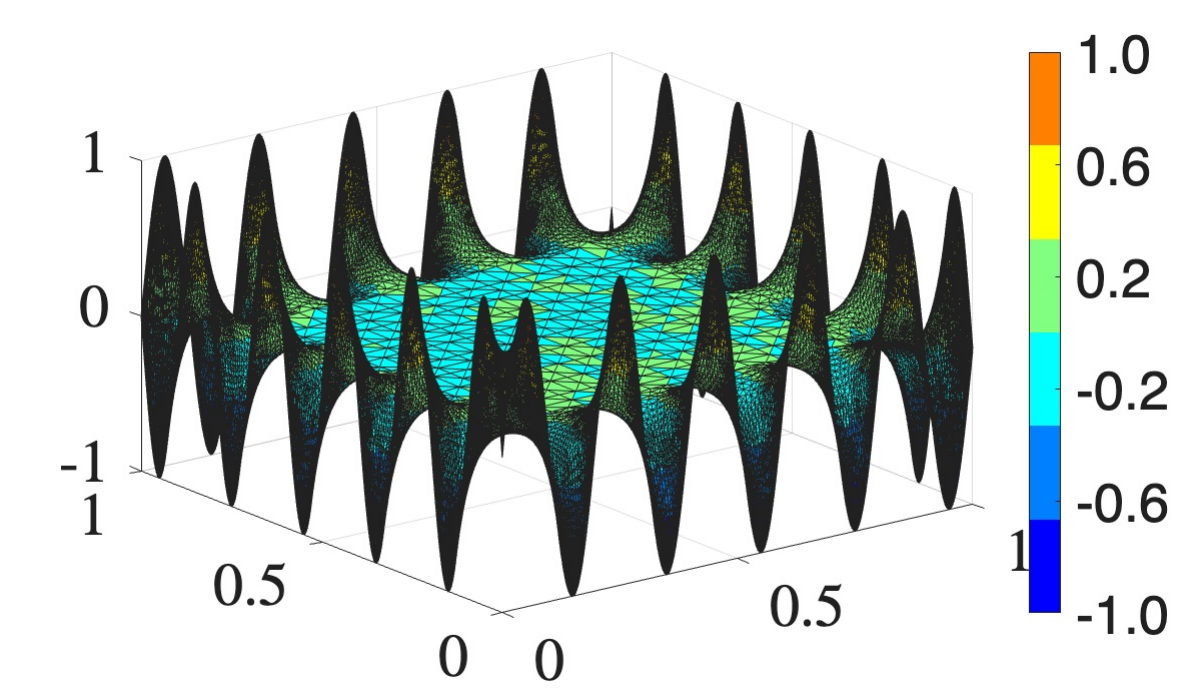}
\subcaption{$\uTFI$}\label{fig:harmonic_osc_i}
\end{subfigure} \hfill 
\begin{subfigure}{0.38\textwidth}
\includegraphics[width=\textwidth]{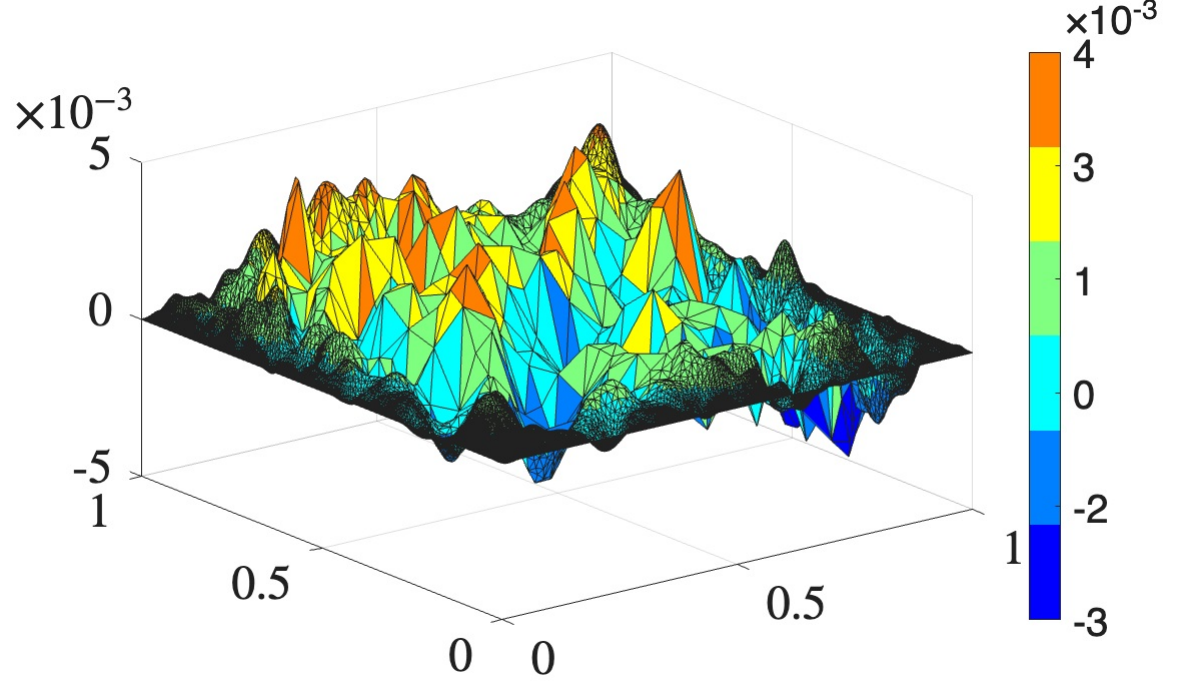}
\subcaption{$u - \uTFI$}\label{fig:harmonic_osc_j}
\end{subfigure}
\begin{subfigure}{0.24\textwidth}
\includegraphics[width=\textwidth]{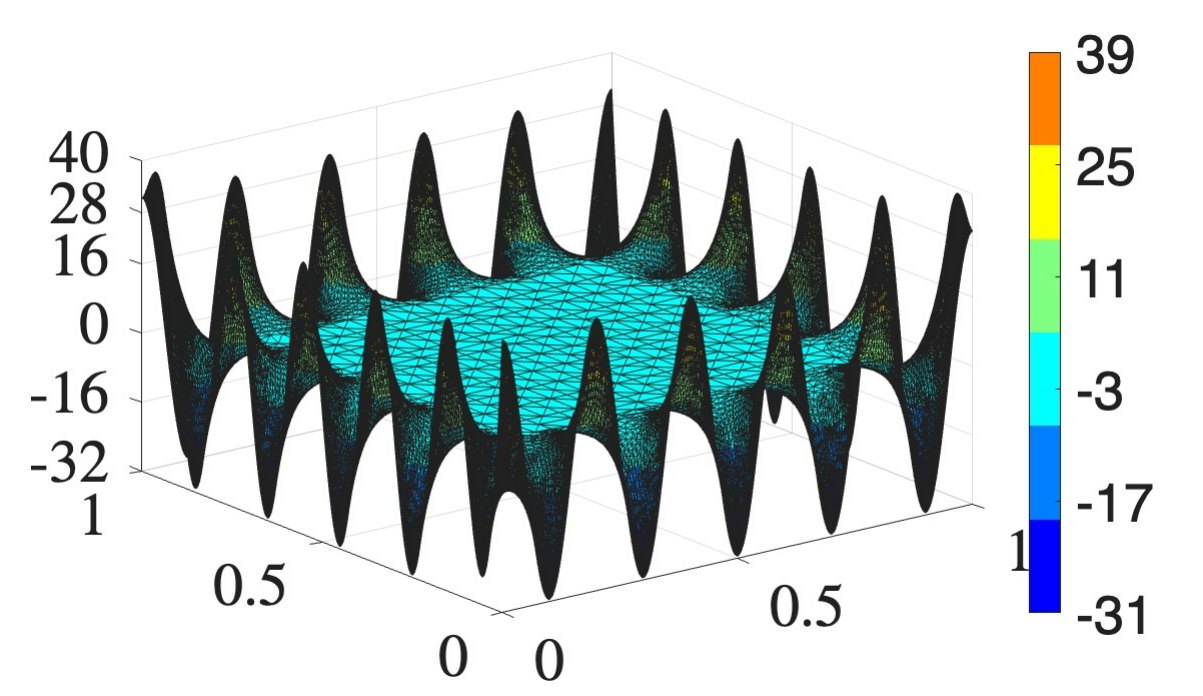}
\subcaption{$\partial \uTFI / \partial x$}\label{fig:harmonic_osc_k}
\end{subfigure} \hfill
\begin{subfigure}{0.24\textwidth}
\includegraphics[width=\textwidth]{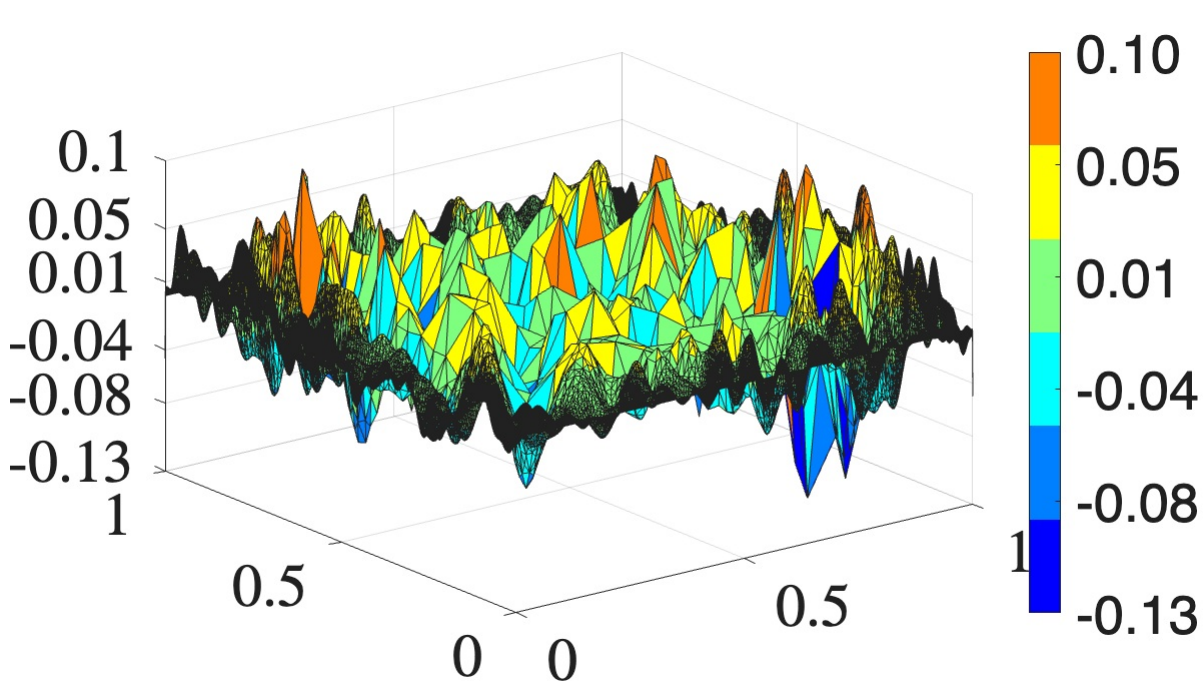}
\subcaption{$\frac{\partial}{\partial x}(u - \uTFI ) $}\label{fig:harmonic_osc_l}
\end{subfigure}  \hfill
\begin{subfigure}{0.24\textwidth}
\includegraphics[width=\textwidth]{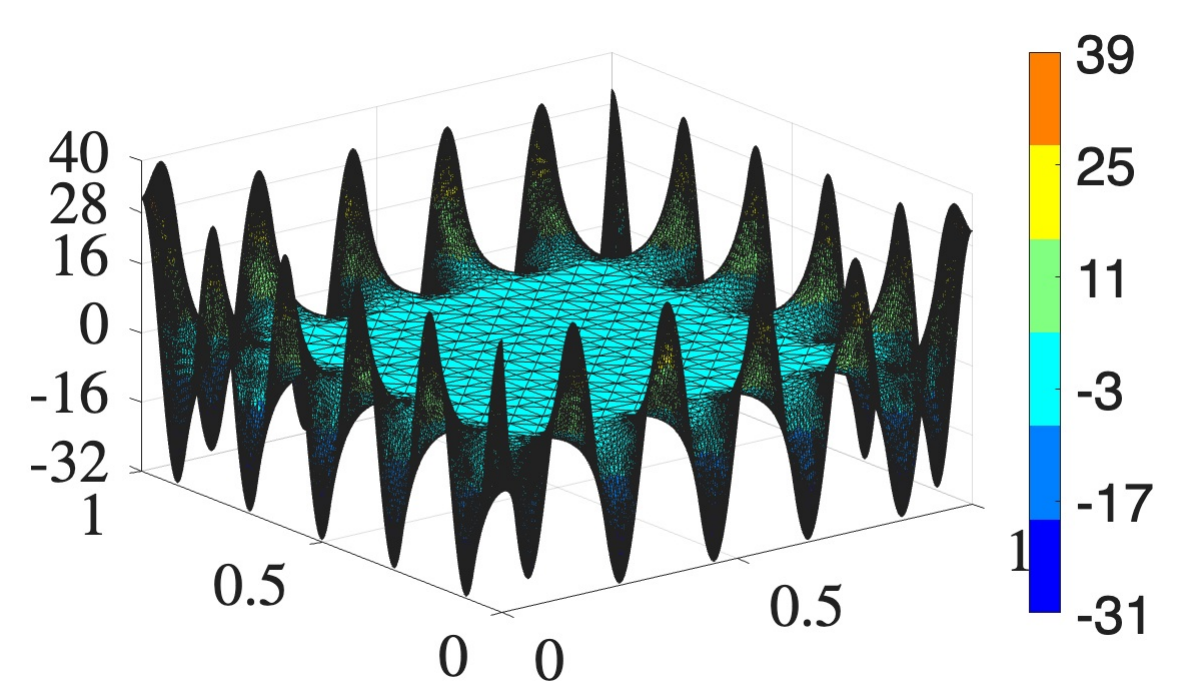}
\subcaption{$\partial \uTFI / \partial y$}\label{fig:harmonic_osc_m}
\end{subfigure}
\begin{subfigure}{0.24\textwidth}
\includegraphics[width=\textwidth]{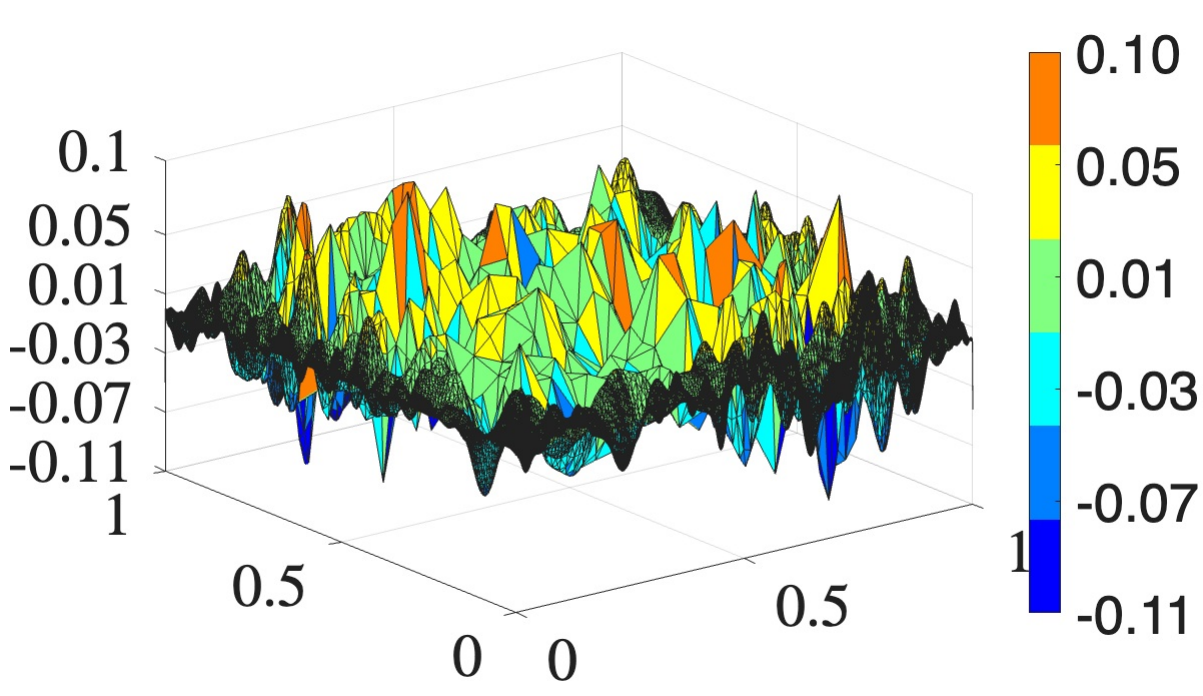}
\subcaption{$\frac{\partial}{\partial y}(u - \uTFI )$}\label{fig:harmonic_osc_n}
\end{subfigure}
\caption{PINN solution with tanh for the harmonic problem with oscillatory boundary
conditions. Reference solution is computed using a lightning Laplace solver~\citep{Trefethen:2020:LLS}.
Network architecture is 4--40--40--40--40--1. 
(a) Training loss; (b) Learning rate; (c) Collocation points for training; and
(d) Testing points for prediction. (e)--(g) Contour
plots of reference solution,
PINN solution and absolute error of PINN solution. (h)--(j)  Surface plots of reference solution,
PINN solution and error of PINN solution. 
(k)--(n) Surface plots of the partial derivatives of the PINN solution and their errors.}
\label{fig:harmonic_osc}
\end{figure}

As the reference solution, we use a `lightning' Laplace
solver~\citep{Gopal:2019:SLP,Gopal:2019:NLH};
Matlab\texttrademark\ code~\citep{Trefethen:2020:LLS} is 
used to compute the exact solution to 10-digits accuracy.\footnote{Using the method of separation of variables, the exact solution is:
$ u(\vx) = \frac{ \sin (10\pi x) 
\cosh ( 10 \pi [ y - 1/2 ] ) +
\sin (10\pi y) 
\cosh ( 10 \pi [ x - 1/2] )  }
{\cosh(5 \pi)}
$.}
Since the oscillatory solution is confined to a narrow band 
in the vicinity of the boundaries, the vertices of a Delaunay mesh that is refined in the vicinity of the
boundary of the square are used as collocation points for training the neural network. A more highly refined Delaunay mesh
is used for the predictions. Figure~\ref{fig:harmonic_osc_c}
shows the $18,000$ 
collocation points used for training. The refined Delaunay mesh consisting of 53,351 testing points is presented in 
~\fref{fig:harmonic_osc_d}. Since the 
exact solution has a high-frequency component with steep gradients that dampens very quickly with distance from the boundary, the 
loss landscape is likely to be stiff, which is characterized by sharp, narrow minima, the presence of saddle points, and large regions where the loss function changes rapidly.
Hence, to mitigate getting trapped at saddle points, network training
is done with Adam 
for $75,000$ epochs and then
L-BFGS for $5,000$ epochs. The tanh activation
function is used.
A learning rate of $10^{-3}$ is used, with
exponential decay and a decay learning rate of $0.99999$. The learning rate drops from $1 \times 10^{-3}$ to $1.1 \times 10^{-5}$ during training as shown in~\fref{fig:harmonic_osc_b}.
We use mini-batch training with $6$ mini-batches, each with a batch size of $3,000$ collocation points. 
Figure~\ref{fig:harmonic_osc_a}
displays the training loss as a function of the
number of epochs; use
of mini-batch training leads to the observed
large fluctuations in the training loss curve. The loss at the end of training drops to $4.76$, which in this instance is reasonably small given that the Laplacian of the boundary functions are of magnitude 
$100 \pi^2$. For this problem, the 
L-BFGS optimizer is not effective in further reducing the loss. 
Figures~\ref{fig:harmonic_osc_e}
and~\ref{fig:harmonic_osc_f} show
contour plots of the exact solution and the
PINN solution; corresponding surface plots are
depicted in Figs.~\ref{fig:harmonic_osc_h}
and~\ref{fig:harmonic_osc_i}. 
Figure~\ref{fig:harmonic_osc_g} shows the contour plot of the
absolute error, and the 
surface plot of the error appears 
in~\fref{fig:harmonic_osc_j}. The 
error in a narrow band of width $0.05$
from the boundaries varies from $-1.5\times 10^{-3}$ to $2 \times 10^{-3}$, while noting that the
exact solution within the same band is in the range of
$\pm 0.4$. At a distance of $0.1$ or greater from the boundaries, the maximum error is larger but in this region the exact solution is in the range
of $\pm 0.04$.  In Figs.~\ref{fig:harmonic_osc_k}--\ref{fig:harmonic_osc_n}, the partial derivatives of
the PINN solution and the errors in
$\partial u /\partial x$ and 
$\partial u / \partial y$ are presented. Steep gradients are
observed near the boundaries, which sharply decrease as one moves away from the boundaries.  The PINN solutions for
the derivatives are fairly accurate (errors are less than $0.1$),
on noting that near the boundaries the partial derivatives of $u$ lie in $(-32,39)$.
As the surface plot
in~\fref{fig:harmonic_osc_h} reveals, the
exact solution is highly oscillatory in the vicinity of the boundaries and sharply dampens away from the 
boundaries---which underscores the importance of
exactly imposing Dirichlet boundary conditions and renders this problem to be particularly 
challenging for PINNs (more so if soft imposition of Dirichlet boundary conditions is used)
considering also
the spectral bias of deep neural network training to
preferentially learn lower frequencies of a target function.

\begin{figure}
\centering
\begin{subfigure}{0.18\textwidth}
\includegraphics[width=\textwidth]{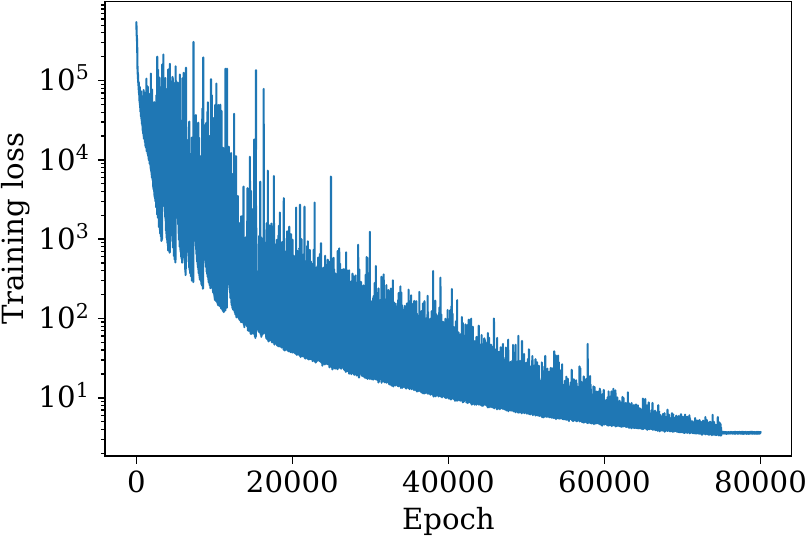}
\subcaption{}\label{fig:harmonic_osc_siren_a}
\end{subfigure} 
\begin{subfigure}{0.17\textwidth}
\includegraphics[width=\textwidth]{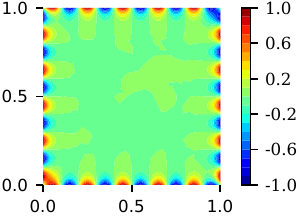}
\subcaption{$\uTFI$}\label{fig:harmonic_osc_siren_b}
\end{subfigure} 
\begin{subfigure}{0.17\textwidth}
\includegraphics[width=\textwidth]{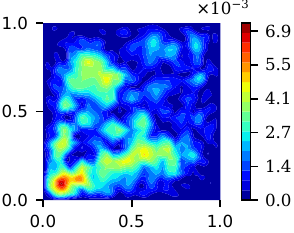}
\subcaption{ $|u - \uTFI|$ }\label{fig:harmonic_osc_siren_c}
\end{subfigure} 
\begin{subfigure}{0.22\textwidth}
\includegraphics[width=\textwidth]{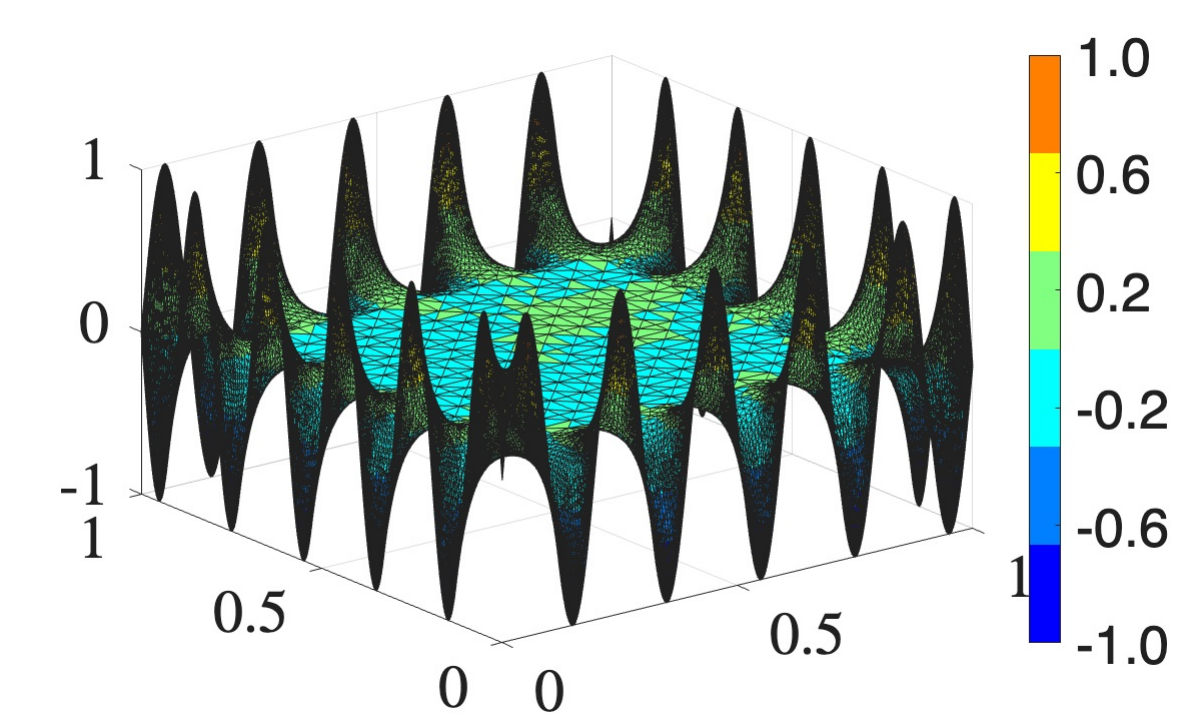}
\subcaption{$\uTFI$}\label{fig:harmonic_osc_siren_d}
\end{subfigure} 
\begin{subfigure}{0.23\textwidth}
\includegraphics[width=\textwidth]{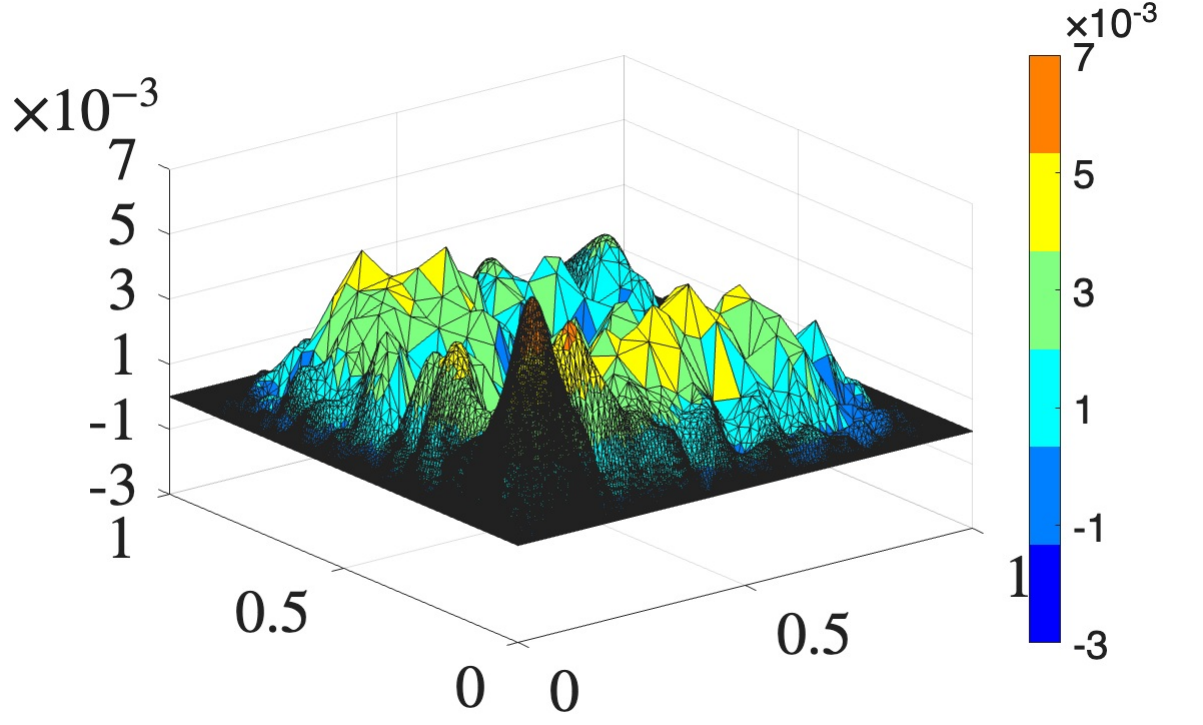}
\subcaption{$u - \uTFI$}
\label{fig:harmonic_osc_siren_e}
\end{subfigure}
\begin{subfigure}{0.24\textwidth}
\includegraphics[width=\textwidth]{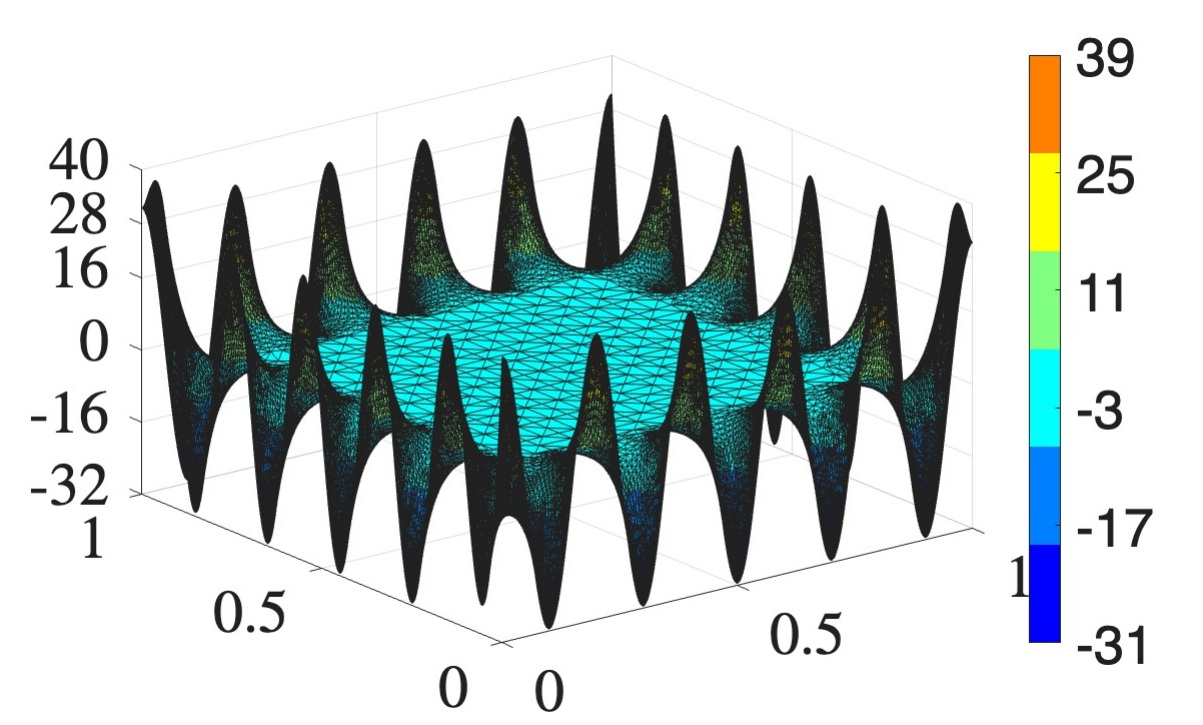}
\subcaption{$\partial \uTFI / \partial x$}
\label{fig:harmonic_osc_siren_f}
\end{subfigure} \hfill 
\begin{subfigure}{0.24\textwidth}
\includegraphics[width=\textwidth]{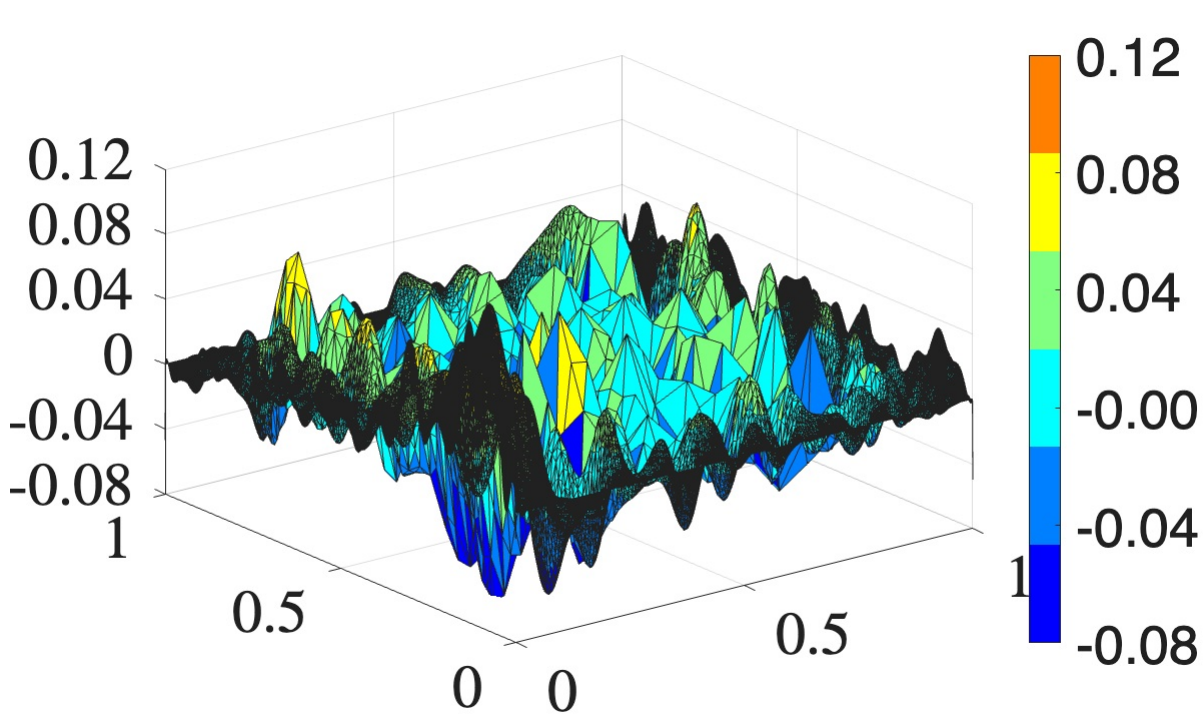}
\subcaption{$\frac{\partial}{\partial x}(u - \uTFI )$}  \label{fig:harmonic_osc_siren_g}
\end{subfigure} \hfill 
\begin{subfigure}{0.24\textwidth}
\includegraphics[width=\textwidth]{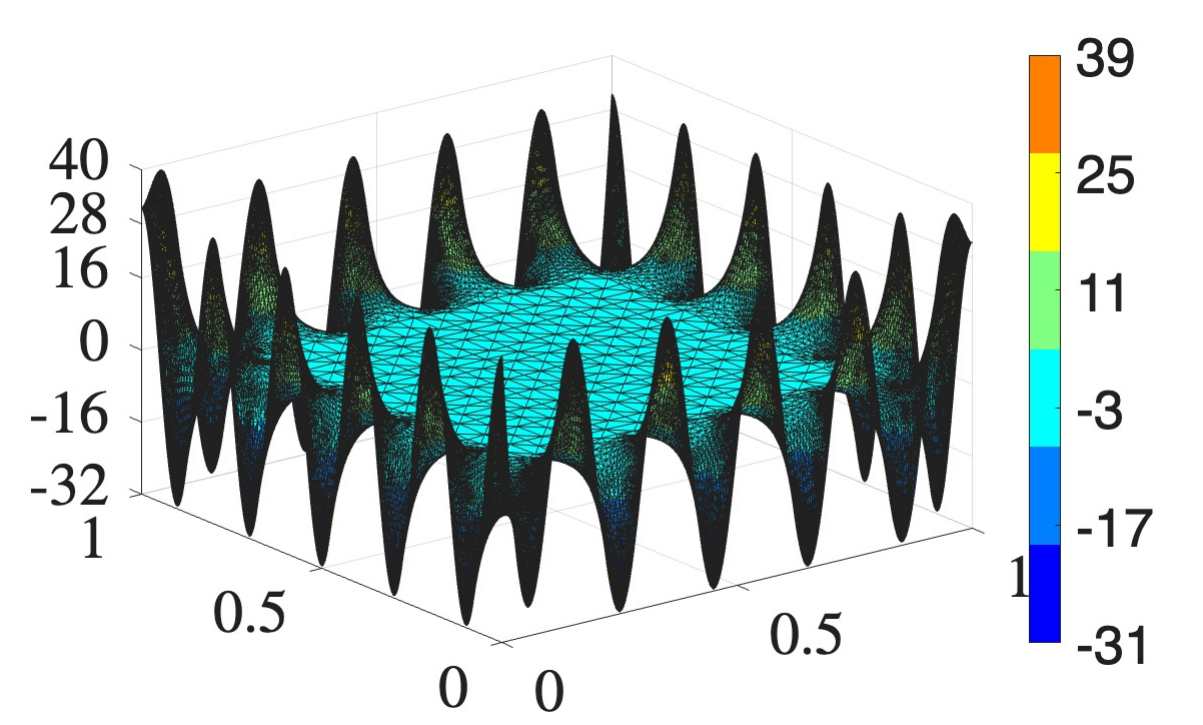}
\subcaption{$\partial \uTFI / \partial y$}\label{fig:harmonic_osc_siren_h}
\end{subfigure} \hfill 
\begin{subfigure}{0.24\textwidth}
\includegraphics[width=\textwidth]{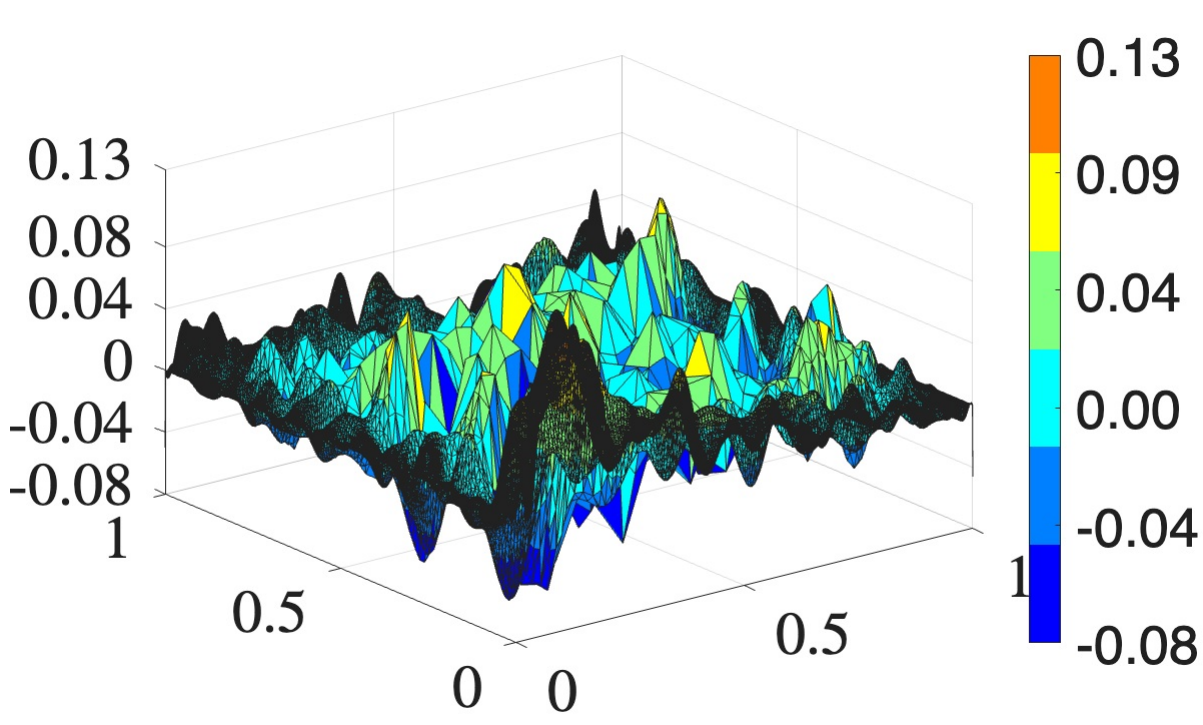}
\subcaption{ $\frac{\partial}{\partial y}(u - \uTFI )$ }
\label{fig:harmonic_osc_siren_i}
\end{subfigure}
\caption{PINN solution with SIREN for the harmonic problem with oscillatory boundary
conditions. 
Network architecture is 4--40--40--40--40--1. 
(a) Training loss; (b) PINN solution; (c)
Absolute error of PINN solution; (d),\,(e)
Surface plot of PINN solution and the
error in PINN;
(f)--(i)  Surface plots of the partial derivatives of the PINN solution and their errors.}
\label{fig:harmonic_osc_siren}
\end{figure}

A potential means to ameliorate the spectral bias is to use Fourier features or the
SIREN (Sinusoidal Representation Network)~\citep{Sitzmann:2020:INR} activation
function that uses the sine function. For this 
oscillatory problem, we adopt SIREN.
For the same set of training points, testing points and
training schedule as adopted for tanh, the numerical
results using SIREN are presented
in~\fref{fig:harmonic_osc_siren}. The loss at the end of training is $3.85$. As 
Figs.~\ref{fig:harmonic_osc_siren_b}--\ref{fig:harmonic_osc_siren_i}
reveal, the errors in $u$, $\partial u / \partial x$ and
$\partial u / \partial y$ with SIREN are 
comparable to those computed using tanh. This also 
confirms the robustness of the
Wachspress-based transfinite formulation 
across different activation functions.

\subsection{Poisson problem on a quadrilateral domain}
We consider
a quadrilateral domain $P$ with vertices
$(0, 0)$, $(1, 0)$, $(1, 1/2)$ and $(0, 1)$, and solve the Poisson problem,
\begin{equation}\label{eq:Poisson_Q}
    - \nabla^2 u  = 2 \ \ \textrm{in } P,
\end{equation}
and choose the exact solution $u$ as
\begin{equation}\label{eq:exact_Poisson_Q}
u(\vx) = (1-y)(- 2 + 2 x + y).
\end{equation}
Dirichlet boundary conditions are imposed
on $\partial P$ that are consistent with the restriction
of the exact solution in~\eqref{eq:exact_Poisson_Q} to the
boundary edges. On solving~\eqref{eq:wsp_quad}, the exact solution
for Wachspress coordinates on the quadrilateral $P$ is:
\begin{equation}\label{eq:wsp_quad_Poisson}
\vm{\lambda} (\vx) = 
\begin{Bmatrix}
  -\dfrac{\left(x-1\right)\,\left(x+2\,y-2\right)}{x-2} ,
  & & \dfrac{x\,\left(x+2\,y-2\right)}{x-2} ,
  & & -\dfrac{2\,x\,y}{x-2} ,
  & & \dfrac{2\,y\,\left(x-1\right)}{x-2} 
\end{Bmatrix}^\top .
\end{equation}        
On using~\eqref{eq:reproducing}, 
the Dirichlet boundary conditions can be cast in terms of $\{\alpha_i\}_{i=1}^4$ as:
\begin{equation*}
\alpha_1(\lambda_2) =  -2(1 - \lambda_2), \quad 
\alpha_2(\lambda_3)= \frac{\lambda_3(2 - \lambda_3)}{4}, \quad
\alpha_3(\lambda_4) = \frac{(1 - \lambda_4)(1 - 3\lambda_4)}{4}, \quad
\alpha_4(\lambda_1) = - \lambda_1 (1 + \lambda_1).
\end{equation*}
On substituting these
in~\eqref{eq:g_polygon} we obtain $g(\vx)$,
and then on using~\eqref{eq:TFI_trial}, the transfinite trial function
for PINNs or deep Ritz is obtained:
\begin{equation}\label{eq:trial_laplace_two}
\uTFI(\vm{\lambda};\vm{\theta})
= g(\vm{\lambda}) + 
\N(\vm{\lambda};\vm{\theta}) - \liftTFI [ \N(\vm{\lambda};\vm{\theta}) ], 
\end{equation}
where $\liftTFI [\N(\vm{\lambda};\vm{\theta}) ]$
is given in~\eqref{eq:PN_quad} for a quadrilateral domain.
\begin{figure}
\centering
\mbox{
\subfloat[]{\includegraphics[width=0.42\textwidth] 
{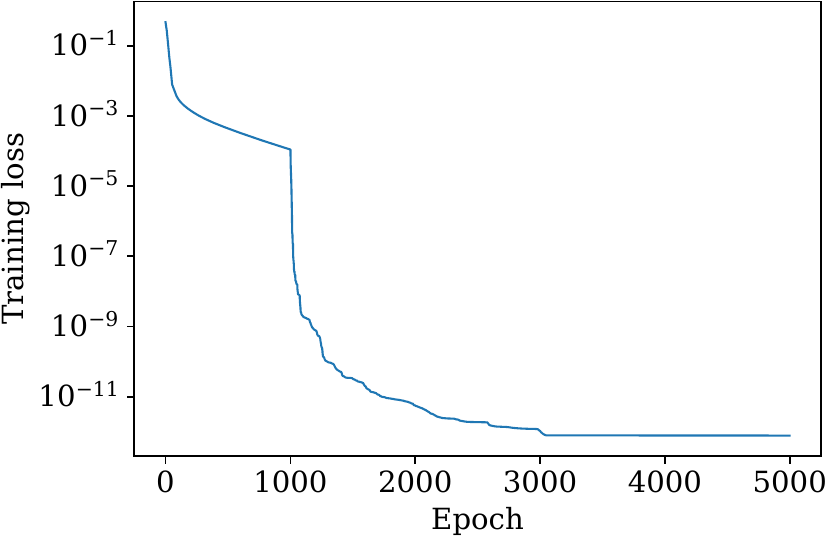}\label{fig:poisson_quad_a}} 
\hfill
\subfloat[]{\includegraphics[width=0.28\textwidth] 
{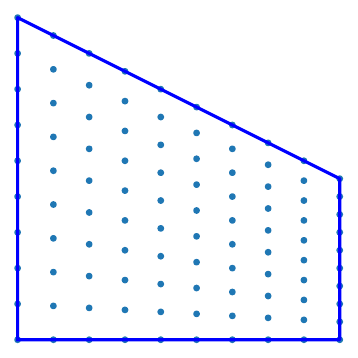}\label{fig:poisson_quad_b}} 
\hfill
\subfloat[]{\includegraphics[width=0.28\textwidth] 
{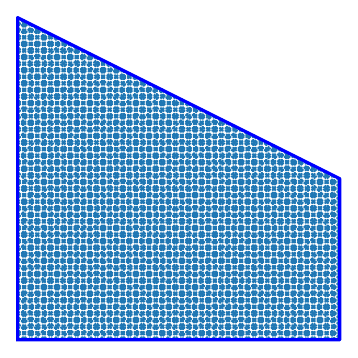}\label{fig:poisson_quad_c}} 
}
\mbox{
\subfloat[$u$]{\includegraphics[width=0.32\textwidth]{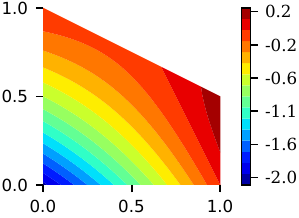}\label{fig:poisson_quad_d}}
\hspace*{0.01in} 
\subfloat[$\uTFI$]{\includegraphics[width=0.32\textwidth]{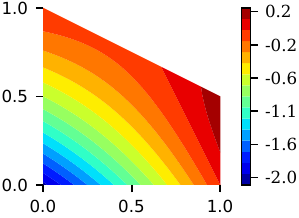}\label{fig:poisson_quad_e}}
\hspace*{0.01in} 
\subfloat[$|u - \uTFI|$]{\includegraphics[width=0.32\textwidth]{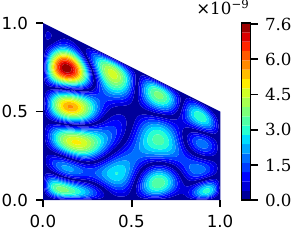}\label{fig:poisson_quad_f}}
}
\mbox{
\subfloat[$|\nabla u |$]{\includegraphics[width=0.32\textwidth]{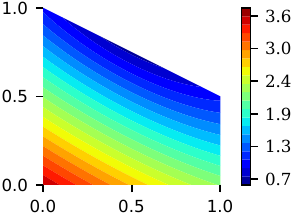}\label{fig:poisson_quad_g}}
\hspace*{0.01in}
\subfloat[$| \nabla \uTFI |$]{\includegraphics[width=0.32\textwidth]{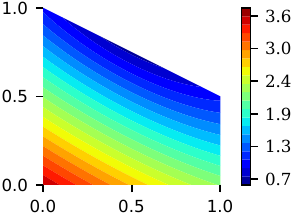}\label{fig:poisson_quad_h}}
\hspace*{0.01in} 
\subfloat[$| \nabla ( u - \uTFI) |$]{\includegraphics[width=0.32\textwidth]{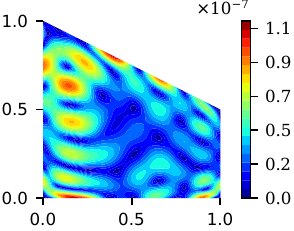}\label{fig:poisson_quad_i}}
}
\caption{PINN solution for the Poisson problem over a quadrilateral domain. Network architecture is 4--20--20--1.  (a) Training loss;
(b) Collocation points for training;
(c) Testing points for predictions;
(d) Exact solution;
(e) PINN solution; 
(f) Absolute error of PINN;
(g) $L^2$ norm of the exact gradient; 
(h) $L^2$ norm of the gradient of the PINN solution; and
(i) $L^2$ norm of the gradient of the error of PINN.
}
\label{fig:poisson_quad}
\end{figure}
For training, the network architecture 4--20--20--1 is used. A uniform 
$10 \times 10$ grid (100 collocation points) on the unit square is mapped to the quadrilateral domain via isoparametric mapping. The testing points for predictions 
are obtained using a Delaunay triangulation of the quadrilateral domain.
Predictions are made at $21,900$ points 
(see~\fref{fig:poisson_quad_c}) that are distinct from the training points 
(see~\fref{fig:poisson_quad_b}). 
Model training consists of $1,000$ epochs with the Adam
optimizer and $4,000$ epochs with the L-BFGS optimizer. The 
$\log$ of the loss is used during the L-BFGS training phase. 
Numerical results are presented in~\fref{fig:poisson_quad}.
Figure~\ref{fig:poisson_quad_a} shows a sharp drop
in training loss after the
L-BFGS optimizer; final loss after 5,000 epochs is $10^{-12}$. Comparisons of PINN and exact solutions are provided in Figs.~\ref{fig:poisson_quad_d}--\ref{fig:poisson_quad_i}. Accuracy of PINN is
${\cal O}(10^{-9})$ in
maximum pointwise absolute error (see~\fref{fig:poisson_quad_f})
and
${\cal O}(10^{-7})$  in the maximum pointwise
$L^2$ norm of the gradient of the error
(see~\fref{fig:poisson_quad_i}).  The very high accuracy in this problem stems from the exact imposition of 
Dirichlet boundary conditions, and in addition, the
simple form (quadratic) of the exact solution in~\eqref{eq:exact_Poisson_Q}.

\begin{figure} 
\centering
\mbox{
\subfloat[]{\includegraphics[width=0.52\textwidth]{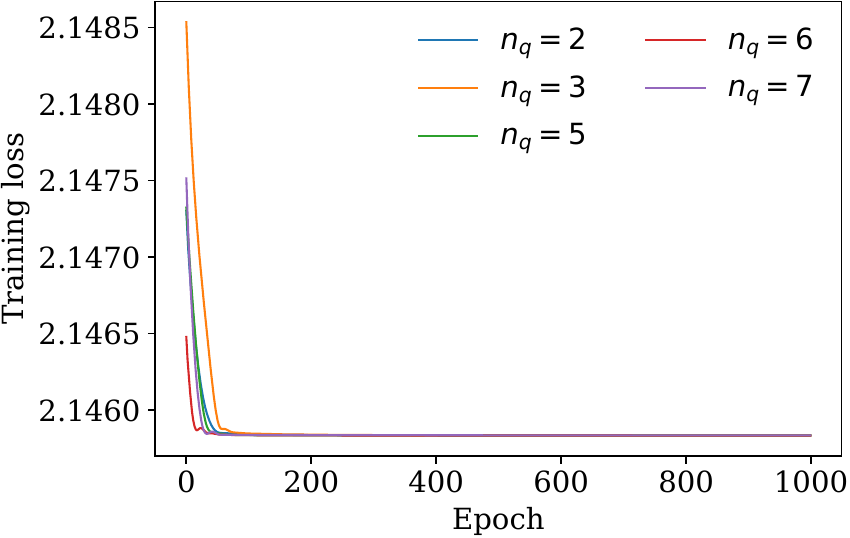}\label{fig:poisson_quad_ritz_a}}
\hspace*{0.15in}
\subfloat[]{\includegraphics[width=0.35\textwidth]{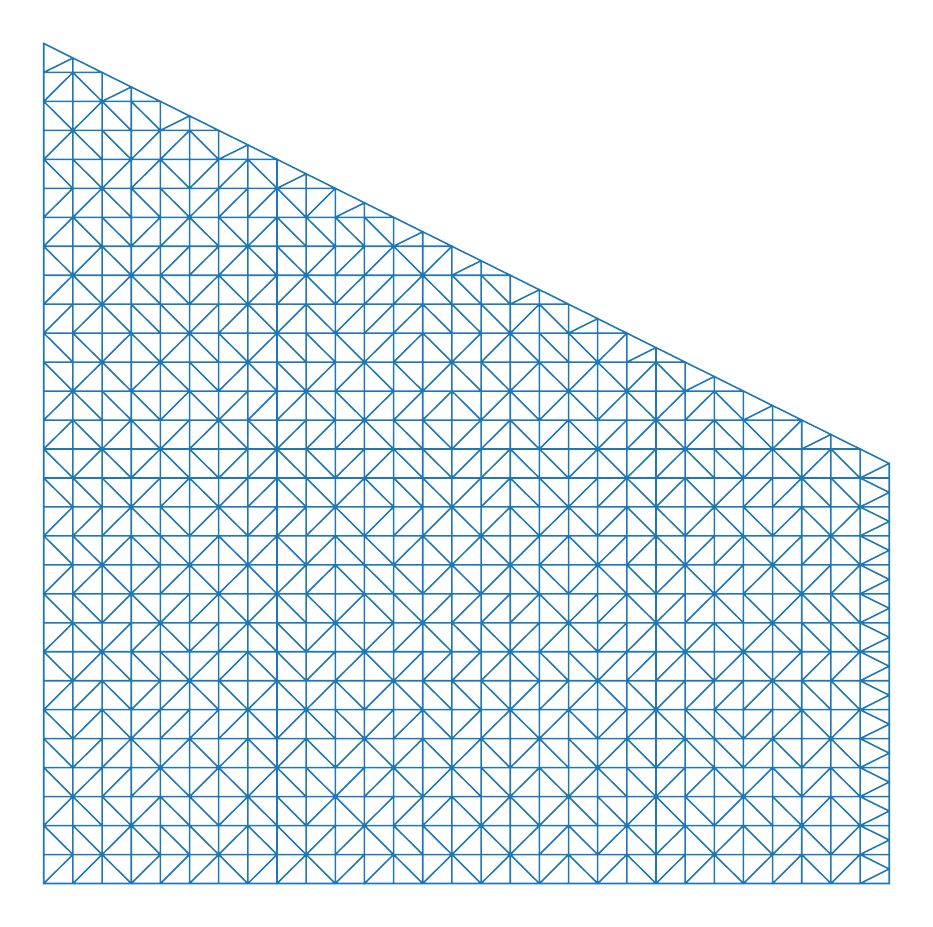}\label{fig:poisson_quad_ritz_b}}
}
\subfloat[$\uTFI$]{\includegraphics[width=0.24\textwidth]{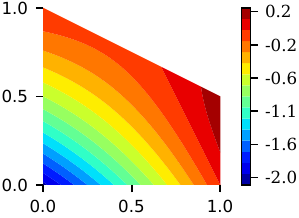}\label{fig:poisson_quad_ritz_c}}
\subfloat[$\nabla \uTFI$]{\includegraphics[width=0.24\textwidth]{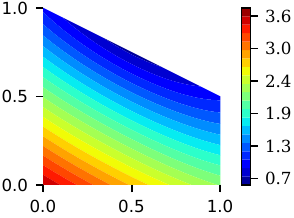}\label{fig:poisson_quad_ritz_d}}
\subfloat[$|u - \uTFI|$]{\includegraphics[width=0.24\textwidth]{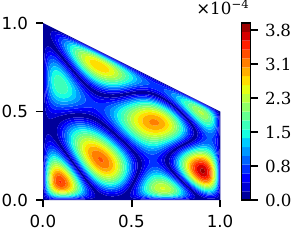}\label{fig:poisson_quad_ritz_e}}
\subfloat[$| \nabla ( u - \uTFI) |$]{\includegraphics[width=0.24\textwidth]{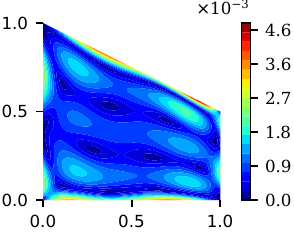}\label{fig:poisson_quad_ritz_f}}
\caption{Deep Ritz solution of the Poisson boundary-value problem 
in~\eqref{eq:Poisson_Q} over a
quadrilateral domain. Exact solution is given in~\eqref{eq:exact_Poisson_Q}. Network architecture is 4--20--20--1.
(a) Training loss;
(b) Training loss is computed using cubature over the Delaunay triangulation;
(c) PINN solution with sixth order ($n_q = 6$) cubature scheme; 
(d) $L^2$ norm of the gradient of the PINN solution;
(e) Absolute error of PINN; and
(f) $L^2$ norm of the gradient of the error of PINN.
}
\label{fig:poisson_quad_ritz}
\end{figure}

For the deep Ritz method, the loss function is set as
the potential
energy (PE) functional, $\Pi[\cdot]$:
\begin{equation}\label{eq:loss_laplace_ritz}
   {\cal L } (\vm{\theta} ) := \Pi[\uTFI(\vm{\theta})] =
   \frac{1}{2} \int_\Omega | \nabla \uxtTFI |^2 \, d\vx - \int_\Omega 2 \uxtTFI \, d\vx, 
\end{equation}
and the optimal $\vm{\theta}^\ast = \argmin\limits_{\vm{\theta}} {\cal L}(\vm{\theta})$.  Since
the kinematically admissible neural network 
trial function~\eqref{eq:trial_laplace_two} is used
in~\eqref{eq:loss_laplace_ritz}, the exact PE serves as a strict lower bound:
$\Pi[\uTFI] \ge \Pi[u]$, where $u$ is the exact solution. In
the numerical computations, the quadrilateral is partitioned into
1,290 triangular cells and a $n_q$ order
quadrature rule~\citep{Xiao:NAC}  
is used in each integration cell. 
For this problem, the potential energy of the exact solution is $2.145833333333333$.
In the numerical computations, 
we use $n_q = 2,\,3,\,5,\,6,\,7$.
Model training consisted of 500 epochs of Adam
and 500 epochs of L-BFGS optimizers, respectively. 
Numerical results are presented in~\fref{fig:poisson_quad_ritz}.
The starting training loss in~\fref{fig:poisson_quad_ritz_a} is 
relatively close to the exact PE and switching to the L-BFGS optimizer at 500 epochs does not 
make a significant difference.
The final  training losses (PE) are $2.145835, 2.145835,2.145834, 2.145834$ and $2.145835$ for integration order 2, 3, 5, 6 and 7, respectively. This yields
and error of $2 \times 10^{-6}$ in the potential energy.
Since the
exact solution is quadratic and the heat source is
a constant, the integrand in the potential energy is also a
quadratic function and so on a refined
Delaunay mesh a cubature rule of order $2$ ($n_q = 2$) suffices.
Since $\uTFI$ is kinematically
admissible, the PE (value of the loss) of the PINN solution is variational (strictly above the exact PE).  
We halt the training at 1,000 epochs since running it longer does not further reduce the loss. The model predictions are made at the same testing points that were used for collocation-based 
PINNs.  We observe from Figs.~\ref{fig:poisson_quad_ritz_e}
and~\ref{fig:poisson_quad_ritz_f} that
the maximum pointwise absolute error in the deep Ritz solution and the maximum pointwise
$L^2$ norm 
of the gradient of the error are
${\cal O}(10^{-4})$ and ${\cal O}(10^{-3})$, respectively. 
The error in the potential energy of the deep Ritz solution is $ 2 \times 10^{-6}$, and so from theory,\footnote{The difference in the potential energy (set $u_\theta = u + e$) is: 
$\Pi[u_\theta]-\Pi[u] = \Pi[u+e]-\Pi[u] =
\int_\Omega \nabla u \cdot \nabla e \, d\vx - \int_\Omega f e \, d\vx + \frac{1}{2} \int_\Omega | \nabla e|^2 \, d \vx$. Since $e$ is a kinematically admissible variation (test function), the first two terms cancel out (statement of the weak form) and we obtain the result: 
$\Pi[u_\theta]-\Pi[u] = \frac{1}{2} \int_\Omega | \nabla e|^2 \, d\vx = \frac{1}{2} |e|_{E(\Omega)}^2$, where $|\cdot|_{E(\Omega)}$ is the energy seminorm of its argument.} the
$H^1$ seminorm of the error should be 
$\sqrt{4 \times 10^{-6}} = 2\times 10^{-3}$, which is consistent with the error plot 
shown in~\fref{fig:poisson_quad_ritz_f}.
The solutions with the
deep Ritz are not as accurate as the collocation approach
(see
Figs.~\ref{fig:poisson_quad_f} and~\ref{fig:poisson_quad_i}),
which is also documented in
the literature.  
The potential energy of the collocation-based PINN solution
is $2.14583333333334$, 
which has an accuracy of ${\cal O}(10^{-14})$. Instead of directly applying the variational
principle (deep Ritz), adopting the 
variational (weak) form in PINNs is known to deliver better accuracy~\citep{Berrone:2022:VPI}. 
For all subsequent problems, we adopt the collocation approach. 

\subsection{Nonlinear Poisson problem}
\begin{figure}
\centering
\mbox{
\subfloat[]{\includegraphics[width=0.41\textwidth] 
{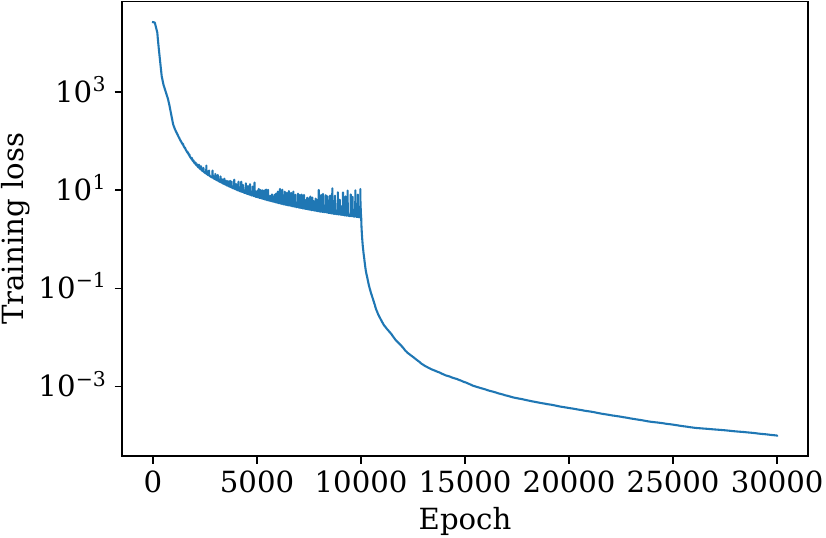}\label{fig:NLP_a}} \hspace*{0.01in}
\subfloat[]{\includegraphics[width=0.28\textwidth]
{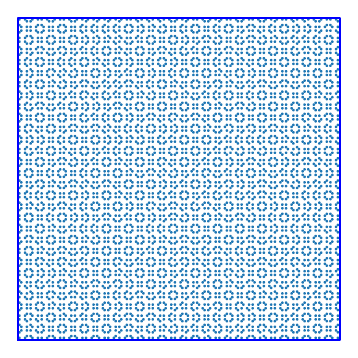}\label{fig:NLP_b}} \hspace*{0.01in}
\subfloat[]{\includegraphics[width=0.28\textwidth]
{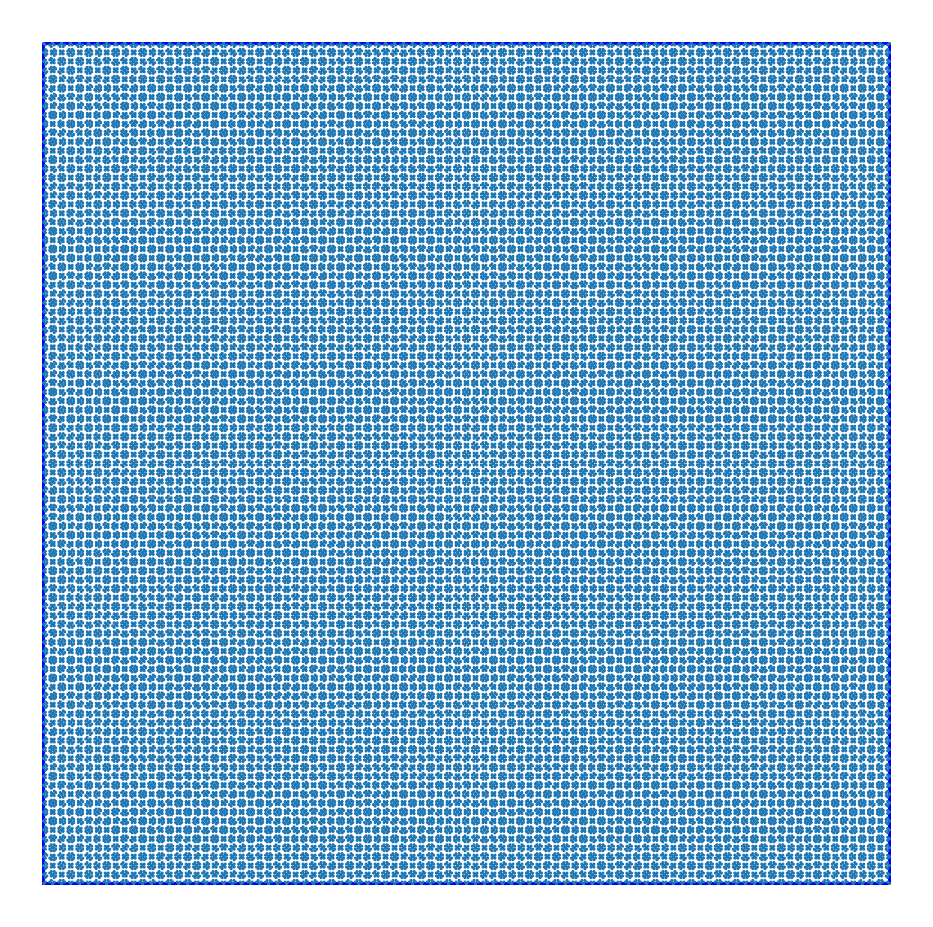}\label{fig:NLP_c}} 
}
\mbox{
\subfloat[]{\includegraphics[width=0.32\textwidth]
{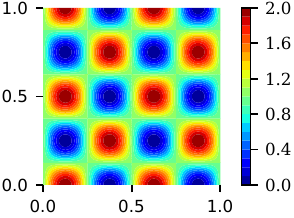}\label{fig:NLP_d}} \hspace*{0.01in}
\subfloat[]{\includegraphics[width=0.32\textwidth] 
{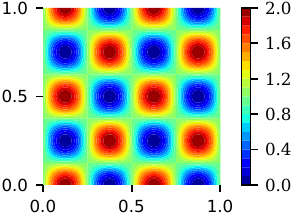}\label{fig:NLP_e}} \hspace*{0.01in}
\subfloat[]{\includegraphics[width=0.32\textwidth]
{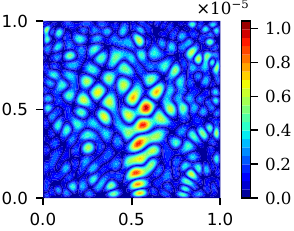}\label{fig:NLP_f}} 
}
\mbox{
\subfloat[]{\includegraphics[width=0.32\textwidth]
{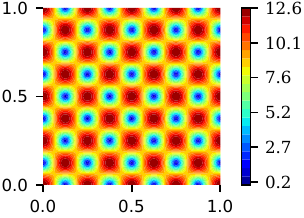}\label{fig:NLP_g}} \hspace*{0.01in}
\subfloat[]{\includegraphics[width=0.32\textwidth] 
{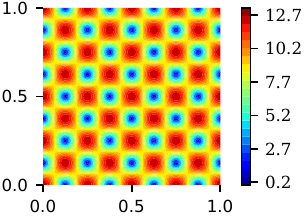}\label{fig:NLP_h}} \hspace*{0.01in}
\subfloat[]{\includegraphics[width=0.32\textwidth]
{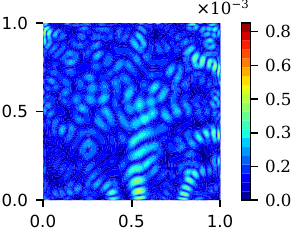}\label{fig:NLP_i}} 
}
\caption{PINN solution for the nonlinear Poisson problem over the unit square. Network architecture is 4--30--30--1.
(a) Training loss;
(b) Collocation points for training;
(c) Testing points for predictions;
(d) Exact solution, $u$; (e) PINN solution, $\uTFI$; and
(f) Absolute error of PINN, $| u - \uTFI| $;
(g) $L^2$ norm of the gradient of the exact solution, 
$| \nabla u |$; 
(h) $L^2$ norm of the gradient of the PINN solution, 
$| \nabla \uTFI |$; and
(i) $L^2$ norm of the gradient of the error of PINN, $| \nabla(u - \uTFI) | $.
}
\label{fig:NLP}
\end{figure}

We consider the
following
nonlinear Poisson problem over the unit square~\citep{Urban:2025:UOP}:
\begin{subequations}\label{eq:nonlinear_Poisson}
\begin{align}
   \nabla^2 u - \exp(u) + f &= 0 \ \ \textrm{in } \Omega = (0,1)^2, \label{eq:nonlinear_Poisson_a}
\intertext{and choose $f$ so that the exact solution $u$ 
is~\citep{Urban:2025:UOP}:}
u(\vx) &= 1 + \sin(4 \pi x ) \cos(4 \pi y) . \label{eq:nonlinear_Poisson_b}
\end{align}
\end{subequations}
The Dirichlet boundary conditions
on $\partial \Omega$ are chosen to be consistent with the exact solution in~\eqref{eq:nonlinear_Poisson_b}, i.e.,
$u = 1 + \sin(4 \pi x)$ on the bottom and top edges; and
$u = 1$ on the left and right edges.
Figure~\ref{fig:NLP_b} shows the 
$5,046$ collocation points used for training. 
A very dense grid 
of testing points, which is shown in~\fref{fig:NLP_c}, is used for predictions against the exact solution. 
The network architecture 4--30--30--1 is used.
Training was performed with Adam for $10,000$ epochs, and then with L-BFGS for $20,000$ epochs. Figure~\ref{fig:NLP_a} shows the
training loss; the loss at the end of training is $10^{-4}$. 
Figure~\ref{fig:NLP_f} depicts the
absolute error between the exact solution and the prediction by PINN, and~\fref{fig:NLP_i} shows the plot of
the $L^2$ norm of the 
gradient of the error.  
The PINN solution is accurate: maximum pointwise absolute error and maximum pointwise
$L^2$ norm of the gradient of the error are
${\cal O}(10^{-5})$ and ${\cal O}(10^{-3})$, respectively.

\subsection{Poisson problem on a parametrized geometry}
Consider a quadrilateral $Q_p$ that is parametrized by
the $y$-coordinate of the third vertex. The vertices of the quadrilateral are: 
$(0,0)$, $(1,0)$, $\bigl( 1,(1+p)/2\bigr)$ and $(0, 1)$ with
$p \in [0,1]$ so that $Q_p$ is convex.  We solve the 
following parametric (geometric) Poisson
boundary-value problem:
\begin{subequations}\label{eq:parametric_bvp}
    \begin{align}
        -\nabla^2 u(\vx,p) &= f(\vx,p) \ \textrm{in } 
        Q_p \subset \Re^2 , \label{eq:parameteric_bvp_a} \\
        u &= 0 \ \textrm{on } \partial Q_p ,\label{eq:parameteric_bvp_b} 
        \intertext{and choose $f$ so that the exact solution is:}
        u(\vx,p) &= 15 xy(1-x)(2 + px - x - 2y) .
        \label{eq:parameteric_bvp_c}
    \end{align}
\end{subequations}
The exact solution is set so that its $L^2$ norm is ${\cal O}(1)$.
For $p = 0$, $||u||_{L^2(\Omega)} = 0.525$ 
and for $p = 0$,  $||u||_{L^2(\Omega)} = 1$.
On using~\eqref{eq:wsp_quad}, the exact solution for Wachspress coordinates on $Q_p$ is obtained:
\begin{equation}\label{eq:parameteric_quad_exact}
\vm{\lambda} (\vx, p) = 
\begin{Bmatrix}
\dfrac{ (1 - x)(2 + px - x - 2y))} {2 + px - x} \  \
\dfrac{ x(2 + px - x - 2y))} {2 + px - x} \  \
\dfrac{2xy} {2 + px - x} \ \
\dfrac{2y(1-x)} {2 + px - x} 
\end{Bmatrix}^\top .
~\end{equation}        

In the computations, the input to the neural
network
is $\vm{\lambda}_p := (\vm{\lambda},p) \in \Re^5$ and
the extended convex domain $\Omega_p \subset \Re^3$ is 
bounded by the planes/surfaces given by $x = 0$, $x = 1$, $y = 0$, 
$y = (2 + px - x)/2$,
$p = 0$, and $p = 1$.  We map each quadrilateral to
the unit square (shown on the left in \fref{fig:mapping_parametric}) via the transformation:
\begin{equation}\label{eq:quad_to_square_map}
    x = \xi, \quad  y = \dfrac{\eta (2 + px - x)}{2},
    \quad J = \frac{2 - (1-p)\xi}{2},
\end{equation}
where $J$ is the Jacobian of the transformation. The Laplacian
of $u$ can be expressed in terms of the derivatives 
with respect to 
$\xi$ and $\eta$ as:
\begin{equation}\label{eq:Laplacian}
    \nabla^2 u (x,y) = u_{,\xi \xi}
    + \dfrac{2(1 - p))\eta}{A(\xi,p)}
      u_{,\xi \eta}
    + \dfrac{4 + \eta^2 (1 - p)^2}{A^2(\xi,p)}
       u_{,\eta \eta}
      + \dfrac{2 (1 - p)^2}{A^2(\xi,p)}
       u_{,\eta} ,
      \ \
      A(\xi,p) = 2 - (1-p) \xi,
\end{equation}
where a comma is used to denote the partial derivative with respect to the specified coordinate.

For collocation-based PINNs, the loss function for the Poisson equation on the parametrized geometry is:
\begin{equation}\label{eq:loss_parametric}
\begin{split}
     {\cal L}(\vm{\theta}) &= \frac{1}{|\Omega_p|}
    \int_0^1 \int_0^1 \int_0^{\frac{2+px-x}{2}} (\nabla^2 u + f )^2 \, dy \, dx \, dp = 
    \frac{1}{|\Omega_p|}
    \int_0^1 \int_0^1 \int_0^1 
     (\nabla^2 u + f )^2
     |J| 
    \, d \xi \, d \eta \, dp \\
    & \approx
    \frac{1}{M} \sum_{k=1}^M
    \bigl[ \nabla^2 \uTFI (\xi,\eta,p ;\vm{\theta} ) + f  
    \bigr]_{(\xi_k,\eta_k,p_k)}^2 \, |J(\xi_k,p_k)| , 
\end{split}
\end{equation}
where $
\uTFI (\xi,\eta,p;\vm{\theta})$ with
$(\xi,\eta) \mapsto \vm{\lambda}$,
is the neural network approximation, and
the expression for the Laplacian is provided in~\eqref{eq:Laplacian}. For training with the Adam
optimizer, the loss function in~\eqref{eq:loss_parametric} is
used, whereas its $\log$ is adopted for training with the
L-BFGS optimizer.
The training points are chosen within the unit cube shown in~\fref{fig:mapping_parametric}. Even
though the cube can be discretized using randomly chosen collocation points, for simplicity and computational efficiency, we discretize the $\xi\eta$-plane once using a set of points, $\vm{P} \in [0,1]^2$, which are uniformly spaced 
(see~\fref{fig:poisson_parametric_b}). This set of points is repeated for different values of $p \in [0,1]$ with increments of $0.1$ to fill the unit cube
with collocation points $\vm{P} \times
[0,0.1,\dots,1]$ (see the schematic shown on the right
in~\fref{fig:mapping_parametric}). In all 5,346 collocation points are used for training. This approach enables us to compute the Wachspress coordinates $\vm{\lambda}(\xi,\eta)$ 
(unit square)
just once and to reuse them for different values of $p$.
Note that on using~\eqref{eq:quad_to_square_map} in~\eqref{eq:parameteric_quad_exact}, 
$\vm{\lambda}(\xi,\eta)$ reduces to bilinear
finite element functions on a square. Training is performed just once to obtain a model that can predict the solution for the family of 
quadrilateral geometries parametrized by $p \in [0,1]$.  Predictions are made for different values of $p$ with testing points that are not in the training set. The numerical results are presented in~\fref{fig:poisson_parametric}. 
Figure~\ref{fig:poisson_parametric_a} shows the training loss as a function of the number
of epochs; the final training loss is 
 $1.4 \times 10^{-7}$. Figures~\ref{fig:poisson_parametric_d} through~\ref{fig:poisson_parametric_l} 
provide comparisons of the predicted solution with the analytical solution for
$p = 0.05$, $p = 0.55$ and $p = 0.75$, which
reveal excellent accuracy of the PINN predictions with maximum pointwise absolute errors of
${\cal O}(10^{-6})$ for the chosen values of $p$.
\begin{figure}
\centering
\begin{tikzpicture}[scale=0.91]

% -------------------------------------------------
% LEFT: Original (x,y) domain
% -------------------------------------------------
\begin{scope}[xshift=0cm,yshift=0cm]
  % polygon
  \draw[blue] (0,0) -- (3,0) -- (3,2) -- (0,3) -- cycle;

  % corner labels
  \node[below] at (0,0) {$(0,0)$};
  \node[below] at (3,0) {$(1,0)$};
  \node[above right] at (0,2.9) {$(0,1)$};
  \node[above] at (3,2) {$(1,\frac{1+p}{2})$};

  % axes
  \draw[red, -{Stealth[length=3mm, width=3mm]}] (3.2,0) -- (4.0,0);
  \node at (4.25,0) {$x$};
  \draw[red, -{Stealth[length=3mm, width=3mm]}] (0,3.2) -- (0,4.0);
  \node at (0,4.25) {$y$};

\end{scope}

% -------------------------------------------------
% Arrow
% -------------------------------------------------
\draw[->, black, very thick] (4.2,1.5) -- (5.5,1.5);

% -------------------------------------------------
% MIDDLE: Reference square (xi, eta)
% -------------------------------------------------
\begin{scope}[xshift=6.5cm,yshift=0cm]
  \draw[blue] (0,0) rectangle (3,3);

  % corner labels
  \node[below]  at (0,0) {$(0,0)$};
  \node[below] at (3,0) {$(1,0)$};
  \node[above right]  at (0,3) {$(0,1)$};
  \node[above] at (3,3) {$(1,1)$};

  % axes
  \draw[red, -{Stealth[length=3mm, width=3mm]}] (3.2,0) -- (4.0,0);
  \node at (4.25,0) {$\xi$};
  \draw[red, -{Stealth[length=3mm, width=3mm]}] (0,3.2) -- (0,4.0);
  \node at (0,4.25) {$\eta$};

\end{scope}

% -------------------------------------------------
% RIGHT: 3D cube with mid-plane
% -------------------------------------------------
\begin{scope}[xshift=13cm,yshift=1.15cm]

  % cube edges
  \draw[blue] (0,0,0) -- (3,0,0) -- (3,3,0) -- (0,3,0) -- cycle;
  \draw[blue] (0,0,3) -- (3,0,3) -- (3,3,3) -- (0,3,3) -- cycle;
  \foreach \x/\y in {0/0,3/0,3/3,0/3}
    \draw[blue] (\x,\y,0) -- (\x,\y,3);

  % mid-plane p = 0.5
  \draw[brown] (0,1.5,0) -- (3,1.5,0) -- (3,1.5,3) -- (0,1.5,3) -- cycle;

  % random dots on mid-plane
  \foreach \x/\z in {
    0.45/0.45, 1.382/0.55, 2.021/0.65,
    0.38/1.27,
    .82/1.27,
    1.32/2.38 , 2.021/2.38, 0.4/2.69,
    2.7/0.3, 2.8/2.03
  }{
    \fill[gray] (\x,1.5,\z) circle (2pt);
  }

  % labels
  \node[below] at (0,0,3) {$(0,0,0)$};
  \node[above] at (3,3,0) {$(1,1,1)$};
  \node[right] at (3.0,1.5,3) {$p\!=\!0.5$};

\end{scope}

\end{tikzpicture}
\caption{Mapping for the parametrized geometry is from
$xy$-space to $\xi\eta$-space (shown on the left). Training points are generated
in the unit cube, where $\xi, \eta, p \in [0,1]$ (shown on the right).}
\label{fig:mapping_parametric}
\end{figure}

\begin{figure}
\centering
\mbox{
\subfloat[]{\includegraphics[width=0.40\textwidth]
{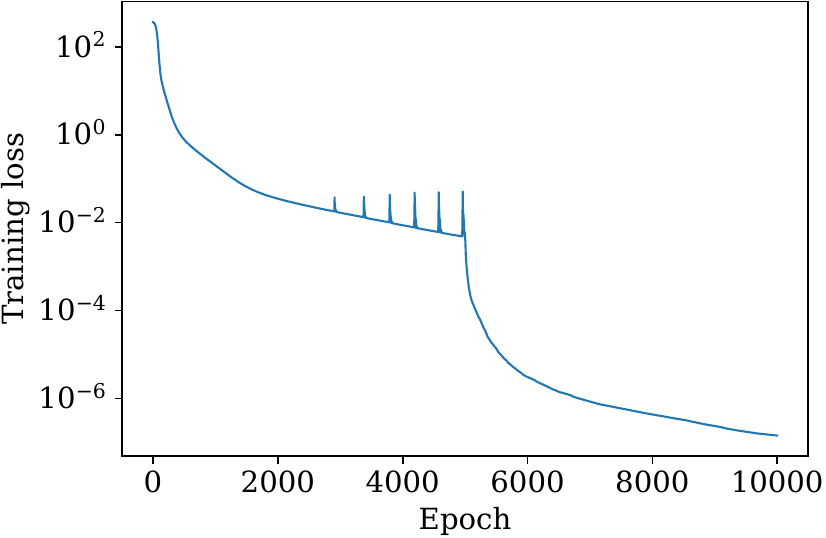}\label{fig:poisson_parametric_a}} \hspace*{0.01in}
\subfloat[]{\includegraphics[width=0.28\textwidth]
{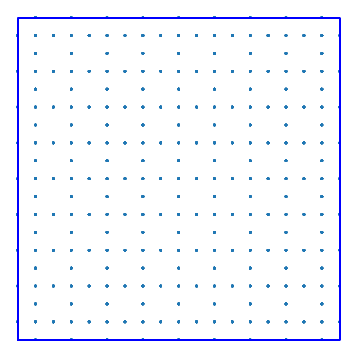}\label{fig:poisson_parametric_b}} \hspace*{0.01in} 
\subfloat[]{\includegraphics[width=0.28\textwidth]
{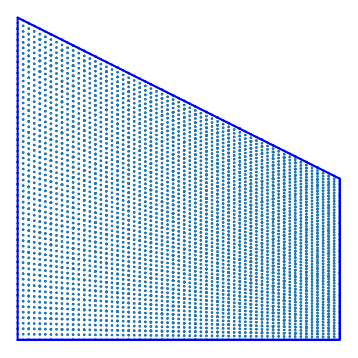}\label{fig:poisson_parametric_c}} 
}
\mbox{
\subfloat[]{\includegraphics[width=0.32\textwidth]
{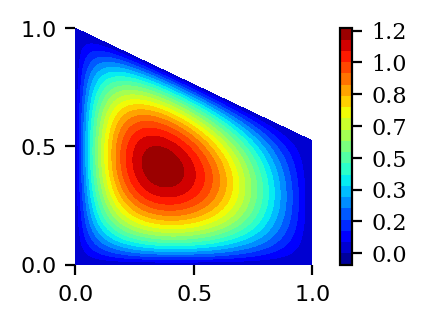}\label{fig:poisson_parametric_d}} \hspace*{0.01in} 
\subfloat[]{\includegraphics[width=0.32\textwidth]
{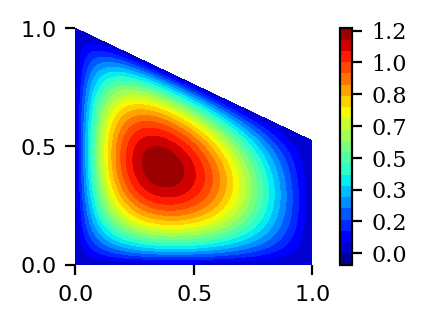}\label{fig:poisson_parametric_e}} \hspace*{0.01in} 
\subfloat[]{\includegraphics[width=0.32\textwidth]
{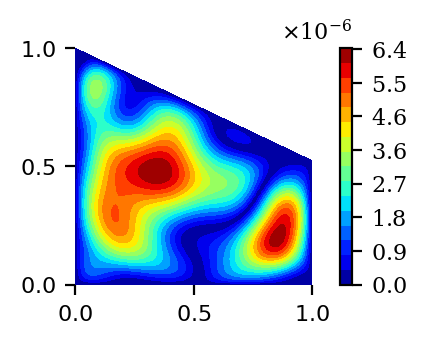}\label{fig:poisson_parametric_f}} 
}
\mbox{
\subfloat[]{\includegraphics[width=0.32\textwidth]
{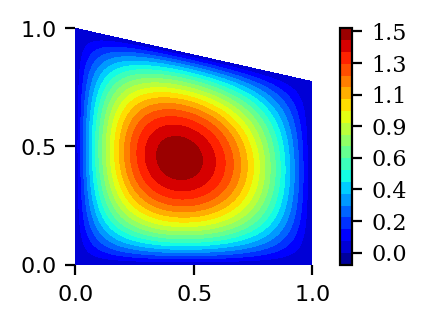}\label{fig:poisson_parametric_g}} \hspace*{0.01in}
\subfloat[]{\includegraphics[width=0.32\textwidth]
{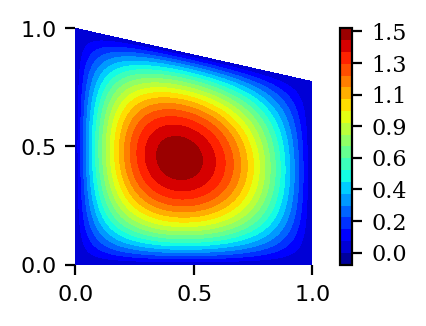}\label{fig:poisson_parametric_h}} \hspace*{0.01in} 
\subfloat[]{\includegraphics[width=0.32\textwidth]
{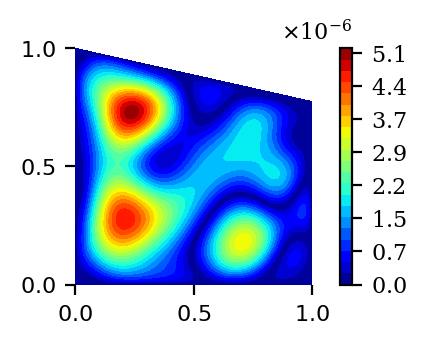}\label{fig:poisson_parametric_i}} 
}
\mbox{
\subfloat[]{\includegraphics[width=0.32\textwidth]
{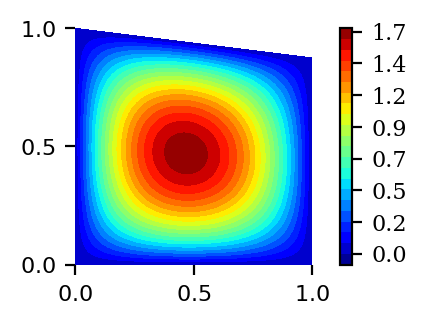}\label{fig:poisson_parametric_j}} \hspace*{0.01in} 
\subfloat[]{\includegraphics[width=0.32\textwidth]
{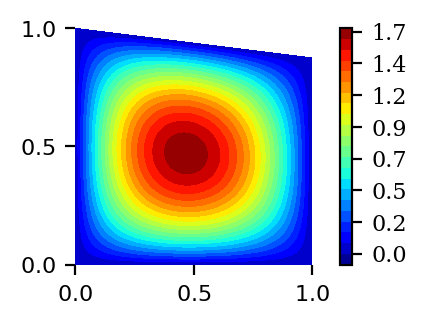}\label{fig:poisson_parametric_k}} \hspace*{0.01in} 
\subfloat[]{\includegraphics[width=0.32\textwidth]
{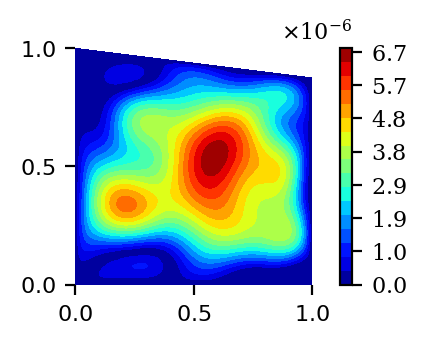}\label{fig:poisson_parametric_l}} 
}
\caption{PINN solution for the Poisson problem over family of parametrized quadrilaterals. Network architecture is 4--20--20--20--20--1. (a) Training loss. (b) Points $\vm{P} \in [0,1]^2$
are shown. Training points in the cube are given by the Cartesian product
$\vm{P} \times [0,0.1,\dots, 1]$
(see~\protect\fref{fig:mapping_parametric}). (c) For a given $p$, testing points within the quadrilateral for predictions. 
For $p = 0.05$: (d) $u$, (e) $\uTFI$, (f) $| u - \uTFI| $.
For $p = 0.55$: (g) $u$, (h) $\uTFI$, (i) $| u - \uTFI| $.
For $p = 0.75$: (j) $u$, (k) $\uTFI$, (l) $| u - \uTFI| $.
}
\label{fig:poisson_parametric}
\end{figure}

\subsection{Inverse heat conduction problem on a pentagonal domain}
\label{subsec:inverse}
\begin{figure}
\centering
\mbox{
\subfloat[]{\includegraphics[width=0.17\textwidth] 
{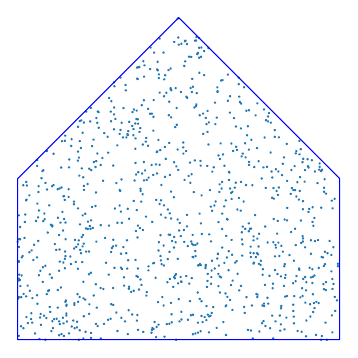}\label{fig:heat_inverse_a}} \hfill
\subfloat[]{\includegraphics[width=0.17\textwidth]
{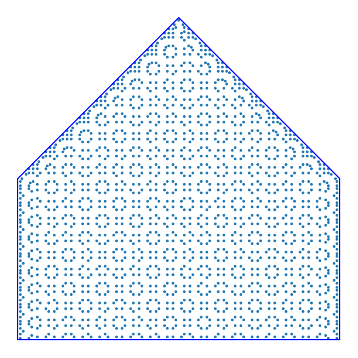}\label{fig:heat_inverse_b}} \hfill
\subfloat[]{\includegraphics[width=0.17\textwidth]
{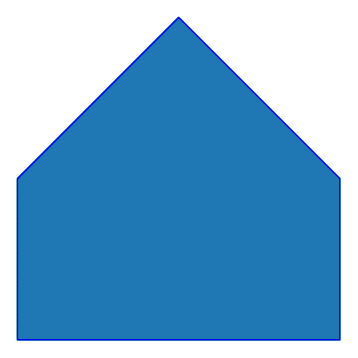}\label{fig:heat_inverse_c}} \hfill
\subfloat[]{\includegraphics[width=0.24\textwidth] 
{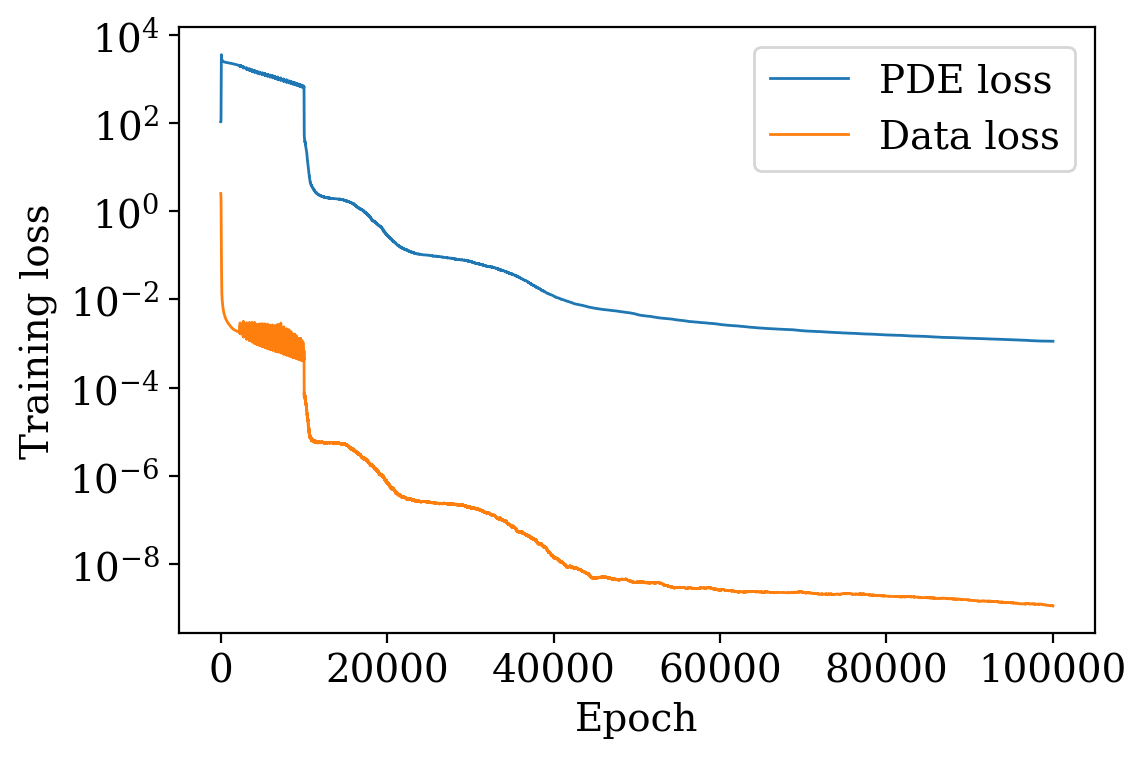}\label{fig:heat_inverse_d}} \hfill
\subfloat[]{\includegraphics[width=0.24\textwidth]
{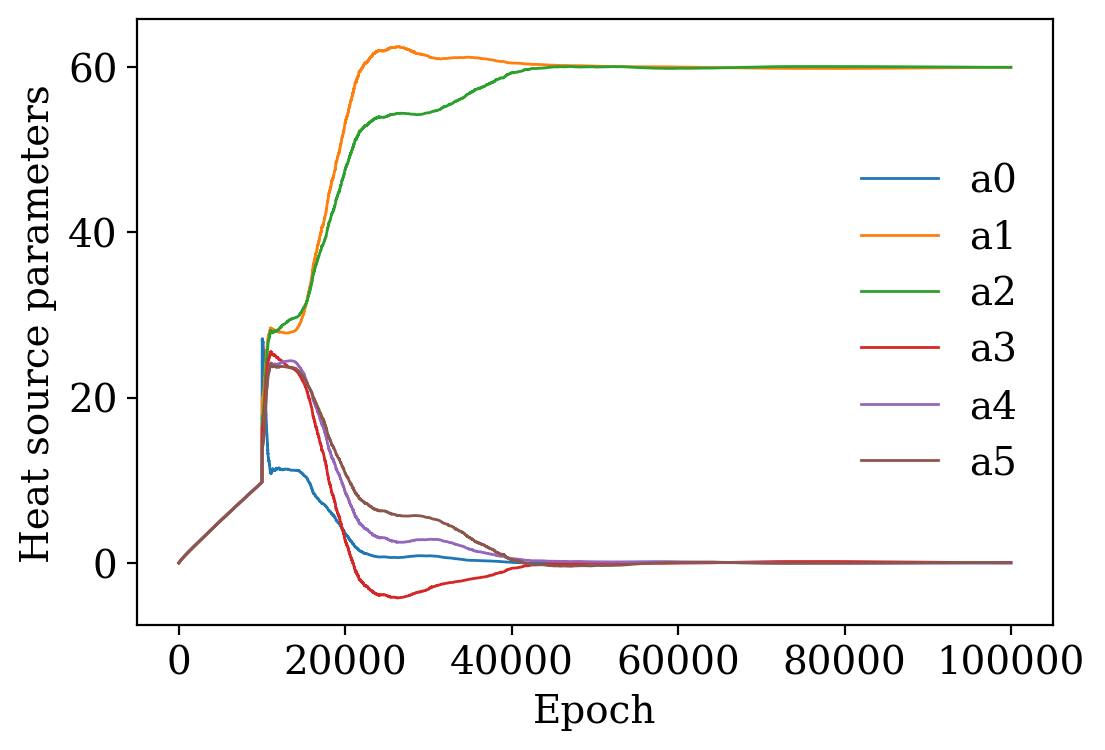}\label{fig:heat_inverse_e}}
}
\mbox{
\subfloat[]{\includegraphics[width=0.33\textwidth]
{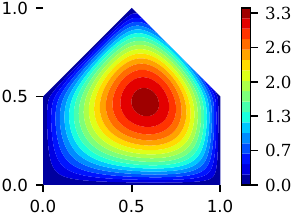}\label{fig:heat_inverse_f}} \hfill
\subfloat[]{\includegraphics[width=0.33\textwidth] 
{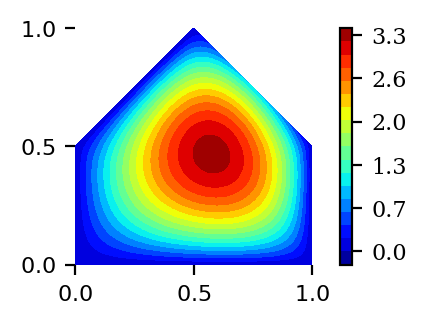}\label{fig:heat_inverse_g}} \hfill
\subfloat[]{\includegraphics[width=0.33\textwidth]
{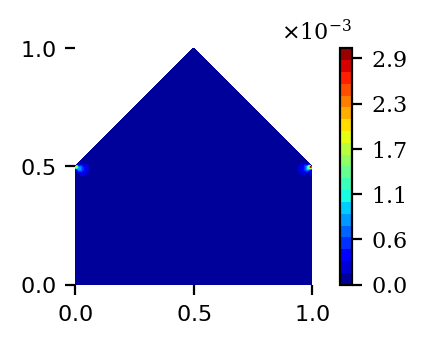}\label{fig:heat_inverse_h}} 
}
\mbox{
\subfloat[]{\includegraphics[width=0.33\textwidth]
{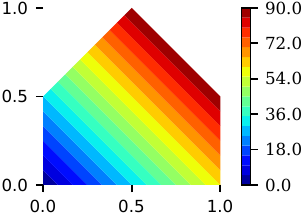}\label{fig:heat_inverse_i}} \hfill
\subfloat[]{\includegraphics[width=0.33\textwidth] 
{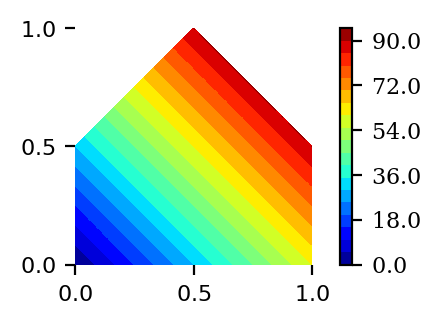}\label{fig:heat_inverse_j}} \hfill
\subfloat[]{\includegraphics[width=0.33\textwidth]
{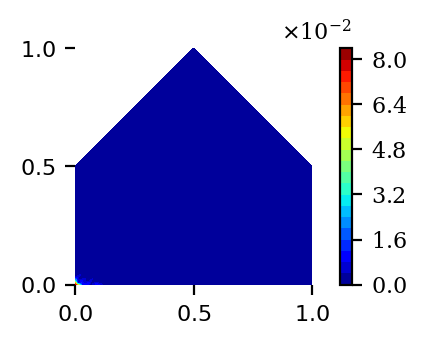}\label{fig:heat_inverse_k}} 
}
\caption{PINN solution for the inverse heat conduction problem over a pentagonal domain. Network architecture is 5--20--20--20--20--1. (a) Randomly sampled data points
for data loss;
(b) training points for PDE loss; and 
(c) testing points for predictions. 
(d) Training loss and (e) evolution of heat source 
parameters. Contour plots of
(f) FE solution, $u^h$; (g) PINN solution, $\uTFI$; and
(h) absolute error of PINN, $| u^h - \uTFI | $. 
Contour plots of 
(i) heat source in forward problem, $f$; 
(j) discovered heat source using PINN, $f_a$; and 
(k) absolute error in heat source, $|f - f_a|$. 
}
\label{fig:heat_inverse}
\end{figure}

Consider a convex pentagonal domain $P$, with vertices
$(0,0)$, $(1,0)$, $(1,1/2)$, $(1/2,1)$ and
$(0,1/2)$. Let $\{ \Gamma_i \}_{i=1}^5$ denote the edges of $P$ in counterclockwise
orientation with $\Gamma_1$ as the bottom edge.
As the forward problem, we solve the 
following steady-state heat conduction problem:
\begin{subequations}\label{eq:inverse_P}
    \begin{align}
        \nabla^2 u + f  &= 0 \ \ \textrm{in } P,
        \label{eq:inverse_P_a} \\
        u = 0 \ \ \textrm{on } \Gamma_1 \cup \Gamma_2 \cup \Gamma_4 \cup \Gamma_5, 
        \ u &= 
         \sin \bigl( \pi (2 y - 1) \bigr) 
        \ \textrm{on } \Gamma_3.
    \end{align}
\end{subequations}
We choose the heat source 
\begin{equation}\label{eq:heat_source}
f(\vx) = 60(x + y), 
\end{equation} 
and solve the forward problem
on a highly refined finite element mesh consisting of 57,454
elements ($6$-noded DC2D6 and $8$-noded DC2D8 quadratic heat transfer elements in Abaqus\texttrademark~\citep{Abaqus:2025}). The finite
element (FE) solution $u^h$ is used as the reference solution. 
As a proof-of-concept, we manufacture a simple
inverse problem to further exemplify the benefits of exactly enforcing Dirichlet boundary conditions.
For the inverse
problem, we sample the 
FE temperature at $m = 1,000$ randomly sampled points 
(see~\fref{fig:heat_inverse_a}). These points in the domain and the corresponding values of $u^h$ at these locations 
are used to construct the 
data loss for PINNs. In the PINN model, the
heat source $f_a$ is assumed to be a bivariate polynomial of degree 2: 
\begin{equation} \label{eq:heat_source_PINN}
f_a(\vx) = a_0 + a_1 x + a_2 y + a_3 x^2 + a_4 y^2 + a_5 x y,
\end{equation}
where $a_i$ $(i = 0,\dots,5)$ are unknown parameters. Now,  
our objective is to discover 
the unknown parameters
in~\eqref{eq:heat_source_PINN} with PINNs, and to assess
its accuracy against the exact heat source $f$ 
in~\eqref{eq:heat_source}.

The PDE loss depends on $f_a$, which is evolving with epochs.
In addition, the given data (primal field $u$) is the 
ground truth. 
Hence, it stands to reason that the data loss is of greater importance than the PDE loss, and so for stronger enforcement of the data constraints, the data loss term is weighted by a factor of $100,000$ compared to the PDE loss in the loss function:
\begin{equation*}
{\cal L}(\vm{\theta} ) =
\frac{1}{M} \sum_{k=1}^M
    \bigl[ \nabla^2 \uTFI ( \vx_k \mapsto \vm{\lambda}^k ;\vm{\theta} )
    + f_a(\vx_k \mapsto \vm{\lambda}^k)    \bigr]^2 +
    \frac{100000}{m} \sum_{\ell = 1}^m 
    \bigl[ u^h(\vx_\ell \mapsto \vm{\lambda}^\ell ) - \uTFI (\vx_\ell \mapsto \vm{\lambda}^\ell ;\vm{\theta}) 
    \bigr]^2 .
\end{equation*}
The network architecture 5--20--20--20--20--1 is used.
Figure~\ref{fig:heat_inverse_b} shows the 1,791 collocation points that are used to compute the PDE loss for training. Training consists of $10,000$ epochs of the Adam
optimizer, followed by $90,000$ epochs of the L-BFGS optimizer. 
Due to the higher weighting on the data loss term in ${\cal L}(\vm{\theta})$, the training process is expected to prioritize minimizing data loss over PDE loss.
Consistent with this expectation, we find the PDE and data losses  in~\fref{fig:heat_inverse_d} to be
$1.1 \times 10^{-3}$ and $1.1 \times 10^{-9}$, respectively.
It should be pointed out 
that a standard PINN setup will require a total of $7$ loss contributions with $5$ of those contributions due to the boundary losses. In contrast, our formulation has only $2$ loss terms (one from the PDE and the other one from the data). Figure~\ref{fig:heat_inverse_e} depicts the evolution of the $6$ unknown heat source parameters during training. The final values of the discovered parameters $a_0$, $a_1$ $a_2$, $a_3$, $a_4$ and $a_5$ are
$0.008007$,  $59.982616$,  $59.973508$,  $0.004112$,  $0.007594$ and  $0.046706$, respectively. Figure~\ref{fig:heat_inverse_c} 
shows the $171,465$ nodal points in the FE mesh where PINN predictions are made. Figure~\ref{fig:heat_inverse_h}
presents the contour plot of the absolute 
error of PINNs. Figures~\ref{fig:heat_inverse_i}--\ref{fig:heat_inverse_k} show
contour plots of the exact heat source $f$, discovered
heat source $f_a$ and the absolute error in the heat source
prediction.
The maximum pointwise error of the predicted heat source over the polygonal domain is about $0.5\%$.

\subsection{Eikonal problem on a pentagonal domain}
\begin{figure}
\centering
\mbox{
\subfloat[]{\includegraphics[width=0.41\textwidth] 
{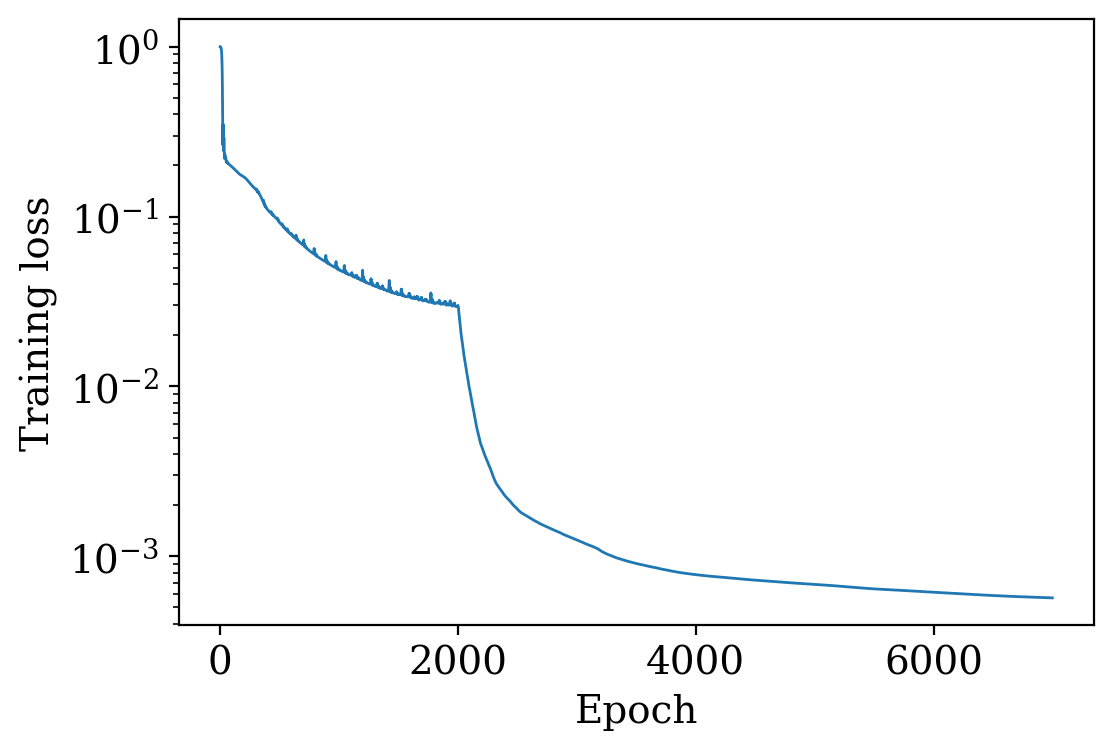}\label{fig:eikonal_a}} \hspace*{0.01in}
\subfloat[]{\includegraphics[width=0.28\textwidth]
{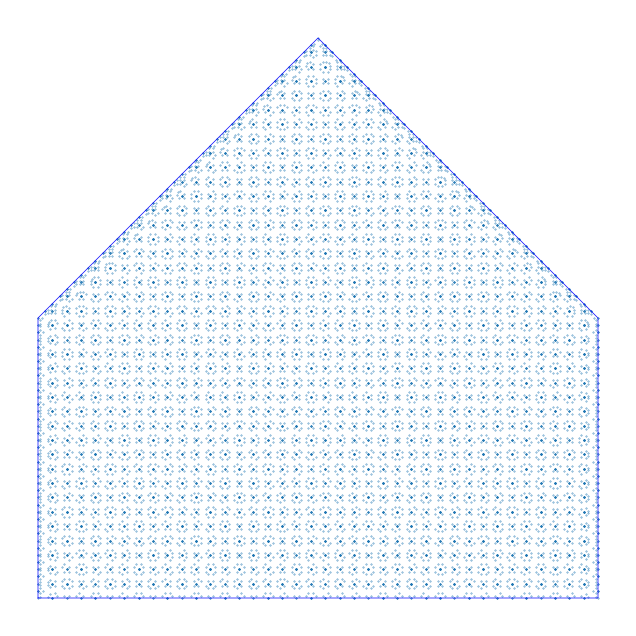}\label{fig:eikonal_b}} \hspace*{0.01in}
\subfloat[]{\includegraphics[width=0.28\textwidth]
{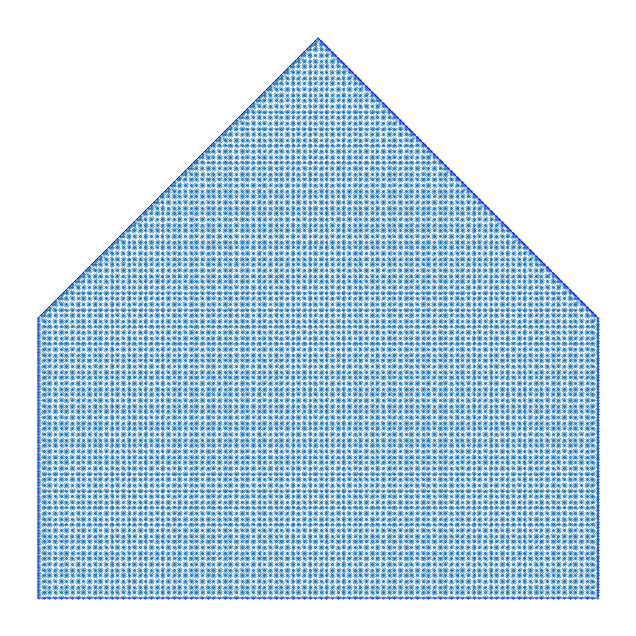}\label{fig:eikonal_c}} 
}
\mbox{
\subfloat[]{\includegraphics[width=0.32\textwidth]
{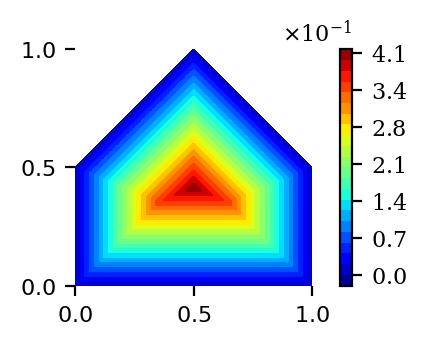}\label{fig:eikonal_d}} \hspace*{0.01in}
\subfloat[]{\includegraphics[width=0.32\textwidth] 
{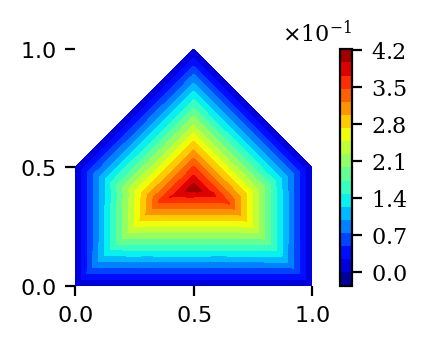}\label{fig:eikonal_e}} \hspace*{0.01in}
\subfloat[]{\includegraphics[width=0.32\textwidth]
{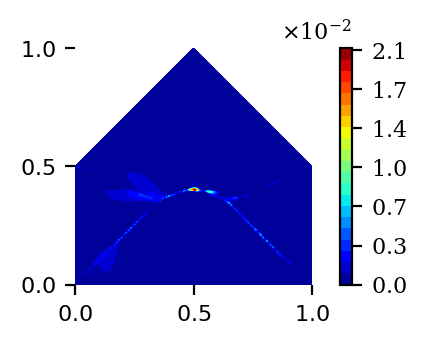}\label{fig:eikonal_f}} 
}
\caption{PINN solution for the Eikonal problem for a pentagon. Network architecture is 5--30--30--30--30--1.
(a) Training loss;
(b) Collocation points for training;
(c) Testing points for predictions;
(d) Exact solution, $u$; (e) PINN solution, $\uTFI$; and
(f) Absolute error of PINN, $| u - \uTFI| $;
}
\label{fig:Eikonal}
\end{figure}

As the final example, we solve the
Eikonal equation over the pentagonal domain $P$ that was considered in the inverse heat conduction problem. The Eikonal problem is:
\begin{subequations}\label{eq:Eikonal_problem}
\begin{align}
   \|\nabla u \| &=  1  \ \ \textrm{in } P,
        \label{eq:eikonal_a} \\
        u &= 0 \ \ \textrm{on } \Gamma = \Gamma_1 \cup \Gamma_2 \cup \Gamma_3 \cup \Gamma_4 \cup \Gamma_5 .
        \label{eq:eikonal_b}
\end{align}
\end{subequations}
The exact solution to this problem is the signed distance function to the boundary $\Gamma$. Given the ridges (along the medial axis) 
that develop in the $C^0$ (nonsmooth)
exact solution, this is a challenging problem 
for PINNs.
A refined set of collocation points is used in the model training.
Figure~\ref{fig:eikonal_b} shows the 
$14,382$ collocation points used for training. 
A very dense grid 
of testing points, which is shown in~\fref{fig:eikonal_c}, is used for predictions against the 
exact solution. 
The network architecture 5--30--30--30--30--1 is used.
Training was performed with Adam for $2,000$ epochs, and then with 
L-BFGS ($\log$ of the loss function) for $5,000$ epochs. Figure~\ref{fig:eikonal_a} shows the
training loss; the loss at the end of training is $5 \times 10^{-4}$. 
Figure~\ref{fig:eikonal_f} depicts the
absolute error between the exact solution and the prediction by PINN. 
The PINN solution is highly accurate across most of the pentagon, except near its centroid. The maximum pointwise absolute error is
${\cal O}(10^{-2})$.

\section{Conclusions}\label{sec:conclusions}
In this paper, we proposed a Wachspress-based transfinite formulation
for physics-informed neural networks 
to exactly enforce
Dirichlet boundary conditions over convex polygonal domains. This overcomes the primary limitation of 
approximate distance 
functions~\citep{Sukumar:2022:EIB}---unbounded Laplacian at the vertices of a  polygonal domain---in 
physics-informed neural networks~\citep{Raissi:2019:PIN}. 
Boolean sum operation, which is used for bilinear Coons transfinite interpolation on the square~\citep{Coons:1967:SCA}, cannot be directly extended 
to  
convex polygonal domains with four (quadrilateral) or more sides. In this work, the
transfinite interpolant was formed by using
Wachspress coordinates as the blending
functions in the formula that is 
based on a
projection onto the faces of a convex domain~\citep{Randrianarivony:2011:OTI}. This generalizes
bilinear Coons transfinite interpolation from a rectangular domain to a convex polygonal domain.
For a polygonal domain $P$ and prescribed
boundary function $\B$ on $\partial P$, the transfinite interpolant of $\B$, $g : \bar P \to C^0(\bar P)$, was viewed as a lifting that
extended the boundary function
to the interior of the domain.
To construct the transfinite trial function, we formed the difference between the neural network's output and
the extension of its boundary restriction to the interior of the polygonal domain, and
then added $g$ to it. Since the restriction
of the trial function to the boundary yielded $\B$, 
strong  enforcement of Dirichlet boundary conditions was ensured.

On convex polygonal domains, Wachspress coordinates served as a geometric feature map that encoded the boundary edges. The spatial coordinates of a point, $\vx \in \bar P$, was
mapped to Wachspress coordinates, $\vx \mapsto \vm{\lambda}$, 
which resided in the geometric feature (map) layer of the neural
network architecture.  In doing so, we showed that 
a framework emerged for solving Poisson
problems on parametrized convex geometries using
neural networks. In particular, the family of
convex quadrilaterals with one vertex at
$\bigl(1,(1+p)/2\bigr)$ $(p \in [0,1])$ 
was embedded in a curved hexahedron that was mapped to a cube. 
Besides the utility of this capability to solve 
partial differential equations using physics-informed neural networks, ideas emanating from this approach  
may also be valuable in the development of geometry-aware
neural operators~\citep{Lu:2021:DLL,Li:2020:FNO,Li:2023:FNO}.

The advance introduced in this paper permitted choosing collocation points in physics-informed neural networks that were located very close to the boundary vertices, thereby overcoming a limitation from 
previous work~\citep{Sukumar:2022:EIB}. Generalized barycentric coordinates are a natural choice as blending functions in the transfinite formula proposed in~\citep{Randrianarivony:2011:OTI}.  On convex polygons, Wachspress coordinates 
are $C^\infty$ (smooth) and their Laplacian are bounded, whereas they are unbounded for mean value coordinates. Hence, 
Wachspress coordinates were adopted in this study.
The performance of the Wachspress-based transfinite 
formulation in PINNs was successfully assessed on several linear and nonlinear problems that included a harmonic 
problem over the unit square
with highly oscillatory boundary conditions,
a Poisson problem over a parametrized quadrilateral,
an inverse heat conduction problem, and the nonlinear Eikonal equation for the distance
function to the boundary of a pentagonal domain. 
Comparisons of the PINN
predictions (absolute error in $u$ and $L^2$ norm of
the gradient of the error) were made with the exact solution (when available) or with reference finite element solutions that were computed
using the Abaqus\texttrademark\ finite element 
software package~\citep{Abaqus:2025}.
We demonstrated the sound accuracy of collocation-based PINNs and deep Ritz on convex polygonal domains, and showed that
model training over a square domain delivered accurate solutions even when interior training points were 
arbitrarily close to the boundary vertices. 

The proposed approach to exactly enforce Dirichlet boundary conditions 
is also suitable for solving 
linear eigenvalue problems over non-Cartesian (convex
polygonal) geometries with the deep Ritz 
(Rayleigh quotient) method. 
In addition, for a convex polygon $P$, if
$w(\vx) := h(\vx) - \liftTFI[h(\vx)]$ for $\vx \in \bar{P}$ 
and $0$ otherwise, we then
observe that $w(\vx) \in C^0(\Re^2)$ is compactly-supported and is a kinematically admissible test function that is suitable in domain-decomposition based variational PINNs. Here, $h(\vx)$ can be 
chosen to be any bivariate
function, including polynomials, sines and cosines, or even the neural network's output with assigned weights and biases.
Furthermore,
extensions of the proposed formulation to convex polyhedra and 
hypercubes in $\Re^n$,
and to nonconvex 
polygonal domains are also of interest. For the former, transfinite formula for convex polyhedra are provided in~\citep{Randrianarivony:2011:OTI} and Matlab\texttrademark\ code to  
compute Wachspress coordinates in three dimensions 
is available~\citep{Floater:2014:GBW}. For the 
latter, generalized barycentric coordinates on
nonconvex polygon such as metric coordinates~\citep{Malsch:2005:STD,Sukumar:2006:RAI}
and variational barycentric coordinates (uses neural fields)~\citep{Dodik:2023:VBC} are worth exploring.
Such topics offer potential directions for future research.

\section*{Acknowledgments}
NS thanks Professor Karniadakis for his generous hospitality during the author's sabbatical visit to Brown University in 2023.
Many helpful discussions with Professor Karniadakis and 
members of the \texttt{Crunch Group} 
at Brown University are also gratefully acknowledged. RR 
acknowledges the support of Dassault Syst{\`e}mes, Inc.

%%%%%%%%%%%%%%%%%%%%%%%%%%%%%%%%%%%%%%%%%%%%%%%%%%%%%%%%%%%%
%  BIBLIOGRAPHY
%%%%%%%%%%%%%%%%%%%%%%%%%%%%%%%%%%%%%%%%%%%%%%%%%%%%%%%%%%%%


\end{document}